\documentclass[12pt,reqno]{amsart}
\usepackage{amsmath}
\usepackage{amssymb}
\usepackage{graphicx}
\usepackage{amsfonts}
\usepackage[margin=0.93in, letterpaper ]{geometry}
\usepackage{lineno}
\usepackage[onehalfspacing]{setspace}
\usepackage{hyperref}

\newtheorem{theorem}{Theorem}
\theoremstyle{plain}

\newtheorem{claim}[theorem]{Claim}
\newtheorem{conclusion}[theorem]{Conclusion}

\newtheorem{corollary}[theorem]{Corollary}

\newtheorem{definition}[theorem]{Definition}
\newtheorem{example}[theorem]{Example}

\newtheorem{lemma}[theorem]{Lemma}

\newtheorem{proposition}[theorem]{Proposition}

\newtheorem{remark}[theorem]{Remark}
\newtheorem{fact}[theorem]{Fact}

\newtheorem{Facts}[theorem]{Facts}
\numberwithin{equation}{section}
\numberwithin{theorem}{subsection}
\input{tcilatex}

\begin{document}
\title{\textsf{Harmonic Analysis on }$GL_{n}$\textsf{\ over Finite Fields}}
\author{\textsf{Shamgar Gurevich}}
\address{\textit{Department of Mathematics, University of Wisconsin,
Madison, WI 53706, USA.}}
\email{shamgar@math.wisc.edu}
\author{\textsf{Roger Howe}}
\address{\textit{Department of Mathematics, Yale University, New Haven, CT
06520, USA.}}
\email{roger.howe@yale.edu}

\begin{abstract}
There are many formulas that express interesting properties of a finite
group $G$ in terms of sums over its characters. For estimating these sums,
one of the most salient quantities to understand is the \textit{character
ratio} 
\begin{equation*}
\frac{\text{trace}(\pi (g))}{\dim (\pi )},
\end{equation*}%
for an irreducible representation $\pi $ of $G$ and an element $g$ of $G$.
For example, in \cite{Diaconis-Shahshahani81} the authors stated a formula
of this type for analyzing certain random walks on $G$.

It turns out \cite{Gurevich-Howe15, Gurevich-Howe17} that for classical
groups $G$ over finite fields (which provide most examples of finite simple
groups) there are several (compatible) invariants of representations that
provide strong information on the character ratios. We call these invariants
collectively \textit{rank}.

Rank suggests a new way to organize the representations of classical groups
over finite and local fields - a way in which the building blocks are the
"smallest" representations. This is in contrast to Harish-Chandra's \textit{%
philosophy of cusp forms} that is the main organizational principle since
the 60s, and in it the building blocks are the cuspidal representations
which are, in some sense, the "largest". The philosophy of cusp forms is
well adapted to establishing the Plancherel formula for reductive groups
over local fields, and led to Lusztig's classification of the irreducible
representations of such groups over finite fields. However, analysis of
character ratios seems to benefit from a different approach.

In this note we discuss further the notion of \textit{tensor rank} for $%
GL_{n}$ over a finite field $\mathbb{F}_{q}$ and demonstrate how to get
information on representations of a given tensor rank using tools coming
from the recently studied \textit{eta correspondence},\textit{\ }as well as
the well known philosophy of cusp forms, mentioned just above.

A significant discovery so far is that although the dimensions of the
irreducible representations of a given tensor rank vary by quite a lot (they
can differ by large powers of $q$), \ for certain group elements of interest
the character ratios of these irreps are nearly equal to each other. Thus,
for purposes of this aspect of harmonic analysis, representations of a fixed
tensor rank form a natural family to study.

For clarity of exposition, we illustrate the developments with the aid of a
specific motivational example that shows how one might apply the results to
certain random walks.
\end{abstract}

\maketitle
\dedicatory{\smallskip\ \ \ \ \ \ \ \ \ \ \ \ \ \ \ \ \ \ \ \ \ \ \ \ \ \ \
\ \ \ \textrm{Dedicated to the memory of Bertram Kostant}}

\section{\textbf{Introduction}}

For a finite group $G$ we consider the set $\widehat{G}$ of (isomorphism
classes of) complex finite dimensional irreducible representations (\textit{%
irreps} for short) of $G,$ and the corresponding collection of irreducible
characters of $G$, 
\begin{equation}
\chi _{\pi },\text{ }\pi \in \widehat{G},  \label{IC}
\end{equation}%
given by $\chi _{\pi }(g)=trace(\pi (g)),$ $g\in G.$

Schur's orthogonality relations \cite{Schur1905} imply that (\ref{IC}) forms
a basis for the space of class functions on $G$. This fact gives birth to
the theory of harmonic analysis on $G$, namely the investigation of class
functions on $G$ via their expansion as a linear combination of irreducible
characters.

Starting with the work of Frobenius \cite{Frobenius1896}, through the work
of Diaconis-Shahshahani \cite{Diaconis-Shahshahani81} and others (see, e.g., 
\cite{Liebeck17, Malle14, Shalev07, Shalev17} and references there),
researchers developed explicit formulas that potentially enable one to apply
the harmonic analysis technique to many class functions that express
interesting properties of $G$.

A closer look at these formulas reveals the fact that in order to make use
of them, in many cases, one needs to have a good solution for the
following:\smallskip

\textbf{Problem (Core problem of harmonic analysis on }$G$\textbf{). }%
Estimate the character ratios%
\begin{equation}
\frac{\chi _{\pi }(g)}{\dim (\pi )},\text{ \ }\pi \in \widehat{G},\text{ }%
g\in G.  \label{CR}
\end{equation}%
We proceed to give an example.

\subsection{\textbf{Hildebrand's Random Walk Example\label{H-Ex-Sec}}}

Consider the group $G=SL_{n}(\mathbb{F}_{q})$ of $n\times n$ matrices with
entries in a finite field $\mathbb{F}_{q}$ and determinant equal to one. For
this example let us assume that $n\geq 3.$ Inside $G$ we look at the
conjugacy class $C$ of the transvection%
\begin{equation}
T\mathcal{=}%
\begin{pmatrix}
1 & 1 &  &  &  \\ 
& \ddots &  &  &  \\ 
&  & \ddots &  &  \\ 
&  &  & \ddots &  \\ 
&  &  &  & 1%
\end{pmatrix}%
,  \label{T}
\end{equation}%
with $T_{ii}=1$ for $i=1,..,n;$ $T_{12}=1,$ and $T_{ij}=0$
elsewhere.\smallskip

The following is known about $C.\smallskip $

\textbf{Fact. }We have\footnote{%
The notation $a(q)=o(b(q))$ means that $a(q)/b(q)\rightarrow 0$ as $%
q\rightarrow \infty .$}$^{\text{,}}$\footnote{%
The notation $c(q)+o(...)$ stands for $c(q)+o(c(q)).$}\smallskip ,

\begin{itemize}
\item The cardinality of $C$ is $q^{2n-2}+o(...)$ \cite{Artin57}.\smallskip

\item Every element of $G$ can be written as a product of no more than $n$
elements from $C$ \cite{Humphries80}. Moreover \footnote{%
We write $a(q)=O(b(q))$ if there is constant $A$ with $a(q)\leq A\cdot b(q)$
for all sufficiently large $q$.},%
\begin{equation}
\frac{\#(G\smallsetminus C^{<n})}{\#(G)}=1-O(\frac{1}{q}),  \label{n}
\end{equation}%
where $C^{<n}=\{g\in G;$ $g=c_{l}...c_{1}$ for $c_{i}\in C$ and $%
l<n\}.\smallskip $
\end{itemize}

Formula (\ref{n}) can be justified for example using the fact that the
elements of $G$ with all eigenvalues $\neq 1$ are outside of $C^{<n}$%
.\smallskip

\begin{figure}[h]\centering
\includegraphics
{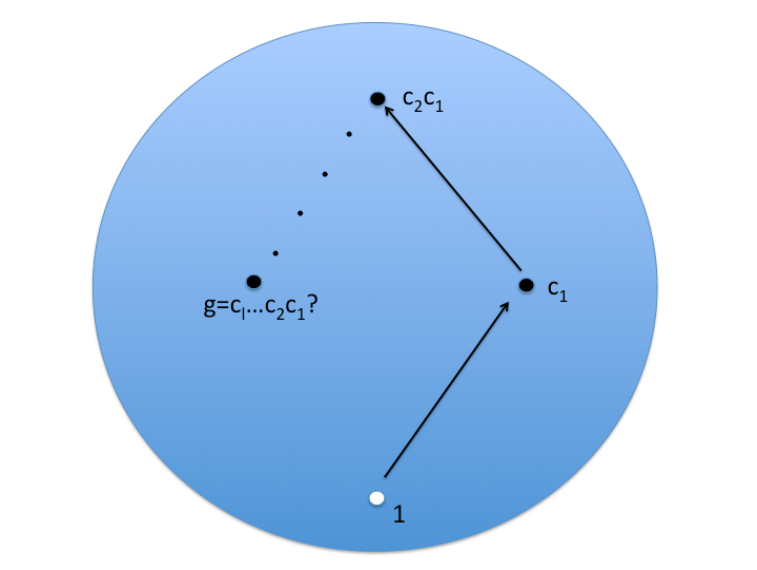}%
\caption{Random walk on $G$ using elements from the conjugacy class $C.$}%
\label{rw}%
\end{figure}

In \cite{Hildebrand92} Hildebrand looked into the problem of generating
random elements of $G$ using random elements from $C$. The mathematical
model is the following random walk on $G$---see Figure \ref{rw} for
illustration. 
We start at the identity element $1$ of $G$. Then we take element $c_{1}$
uniformly at random from $C$ and "walk" to $c_{1}.$ We can continue in this
manner and walk to $c_{2}c_{1}$, then to $c_{3}c_{2}c_{1}$ etc.

Let us denote by $P_{C}^{\ast l}(g)$ the probability that in this way after $%
l$ steps the product $c_{l}...c_{1}$ is equal to $g.$ A very general
argument \cite{Lovasz93} implies that $P_{C}^{\ast l}$ approaches the
uniform distribution $U$ on $G$ as $l\rightarrow \infty .$

To say more, \cite{Hildebrand92} consider the distance in total variation
between $P_{C}^{\ast l}$ and $U,$%
\begin{equation}
\left\Vert P_{C}^{\ast l}-U\text{ }\right\Vert _{TV}=\max_{S\subset
G}\left\vert \text{ }P_{C}^{\ast l}(S)-U(S)\right\vert .  \label{TV}
\end{equation}

It is easy to see that $\left\Vert \cdot \right\Vert _{TV}$ is equal $\frac{1%
}{2}\left\Vert \cdot \right\Vert _{L_{1}}$, i.e., half of the $L^{1}$-norm
on $G$ \cite{Diaconis-Shahshahani81}.

\begin{figure}[h]\centering
\includegraphics
{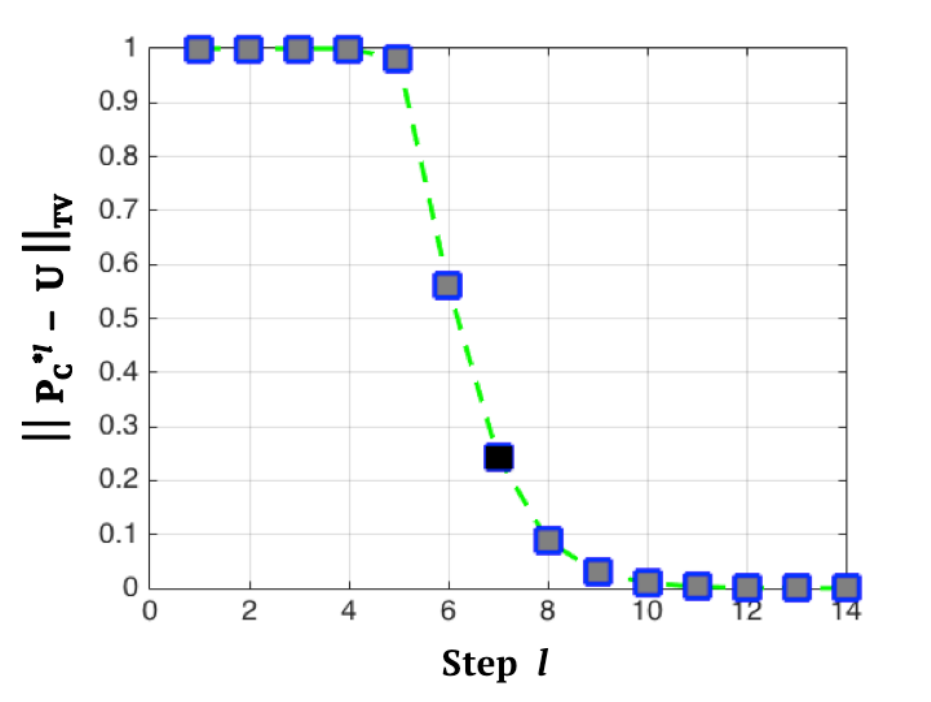}%
\caption{{}Numerics suggests that the mixing time for $G=SL_{7}(\mathbb{F}%
_{3})$ is $l_{M}\approx 7$.}\label{rw-t-sl7_3}%
\end{figure}

\begin{figure}[h]\centering
\includegraphics
{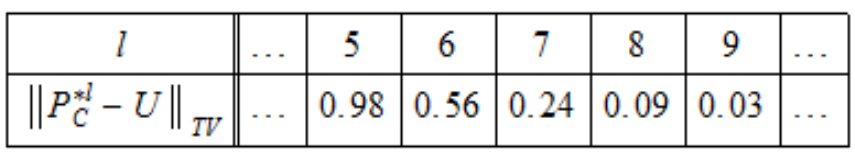}%
\caption{Numerics suggests that the mixing rate for $G=SL_{7}(\mathbb{F}_{3})
$ is $r_{M}\approx \frac{1}{3}$ .}\label{cutoffsl7_3.bmp}%
\end{figure}

The cutoff phenomenon \cite{Diaconis96} suggests that convergence to
uniformity might show a sharp cutoff, namely---see Figures \ref{rw-t-sl7_3}
and \ref{cutoffsl7_3.bmp} for illustration---the distance (\ref{TV}) stays
close to its maximum value (which is $1$) for a while, then suddenly at some
step $l_{M}$ (called \textit{mixing time})\textit{\ }drops to a quite small
value and then tends to zero exponentially fast with some exponent (called 
\textit{mixing rate}) $r_{M}$ \cite{Lovasz93}.


In our case, Formula (\ref{n}) implies that $l_{M}$ can not be less than $n$
and the numerics\footnote{%
The numerics appearing in these notes were generated with John Cannon
(Sydney) and Steve Goldstein (Madison).} that appears in Figure \ref%
{rw-t-sl7_3} illustrates, in particular, the fact that $n$ steps are
probably enough.

\begin{theorem}
\label{H-Thm}The random walk on $G=SL_{n}(\mathbb{F}_{q}),$ $n\geq 3$, using
the collection $C$ of transvections has, for sufficiently large $q$,

\begin{enumerate}
\item Mixing time $l_{M}=n$.

\item Mixing rate $r_{M}=\frac{1}{q}+O(\frac{1}{q^{n}}).$
\end{enumerate}
\end{theorem}

Theorem \ref{H-Thm} was first proved in \cite{Hildebrand92}. The results of
this note will, among other things, provide a new proof.

\subsection{\textbf{Harmonic Analysis of the Random Walk\label{S-HA-of-RW}}}

Diaconis and Shahshahani developed in \cite{Diaconis-Shahshahani81} formulas
that, in principle, enable one to estimate the mixing time $l_{M}$ and
mixing rate $r_{M}$ for random walks on finite groups. Here is the
description that is relevant for us.

The probability distribution $P_{C}^{\ast l}$ that we defined in Section \ref%
{H-Ex-Sec} is a class function on $G$, and its expansion in terms of
irreducible characters can be computed explicitly.

\begin{proposition}
We have, 
\begin{equation}
P_{C}^{\ast l}=\frac{1}{\#(G)}\sum_{\pi \in \widehat{G}}\dim (\pi )\left( 
\frac{\chi _{\pi }(T)}{\dim (\pi )}\right) ^{l}\chi _{\pi },  \label{PCl}
\end{equation}%
where $T$ is the transvection (\ref{T}).
\end{proposition}

Indeed, Formula (\ref{PCl}) can be verified using the fact that $P_{C}^{\ast
l}$ is the $l$-fold convolution of $P_{C}$ with itself, and the standard
identity for convolution of two irreducible characters.

From (\ref{PCl}) we obtain:

\begin{corollary}
\label{C-DS}For the random walk on $G$ using $C$ we have,

\begin{enumerate}
\item The total variation distance of $P_{C}^{\ast l}$ from uniformity
satisfies%
\begin{equation}
\left\Vert \text{ }P_{C}^{\ast l}-U\right\Vert _{TV}^{2}\leq \frac{1}{4}%
\underset{\mathbf{1}\neq \pi \in \widehat{G}}{\sum }\dim (\pi
)^{2}\left\vert \frac{\chi _{\pi }(T)}{\dim (\pi )}\right\vert ^{2l}.
\label{R_el}
\end{equation}

\item The mixing rate satisfies $r_{M}=\underset{\mathbf{1}\neq \pi \in 
\widehat{G}}{\max }\left\vert \frac{\chi _{\pi }(T)}{\dim (\pi )}\right\vert
.$
\end{enumerate}
\end{corollary}

Part 2 of Corollary \ref{C-DS} is immediate from (\ref{PCl}), while for Part
1 one might in addition use the fact that the total variation norm is half
of the $L^{1}$-norm, then apply Cauchy--Schwartz inequality, and finally use
Schur's orthogonality of characters.

\begin{figure}[h]\centering
\includegraphics
{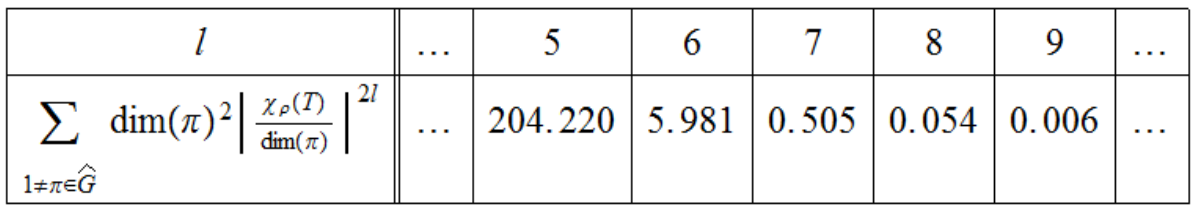}%
\caption{The sum at the right-hand side of (\protect\ref{R_el}) for $%
G=SL_{7}(\mathbb{F}_{3})$. }\label{r_ell}%
\end{figure}%

The numerics appearing in Figure \ref{r_ell} illustrates the possibility
that a good bound on the sum at the right-hand side of (\ref{R_el}) will
give the desired information on the mixing time $l_{M}$.

In order to use Corollary \ref{C-DS} to verify Theorem \ref{H-Thm}, we want
to have a method to get information on the dimensions $\dim (\pi ),$ and
most importantly on the character ratios $\frac{\chi _{\pi }(T)}{\dim (\pi )}
$ of the irreps $\pi $ of $G=SL_{n}(\mathbb{F}_{q})$ at the transvection $T$
(\ref{T}).

Recently, in \cite{Gurevich-Howe15, Gurevich-Howe17}, we have discovered
such a method, that seems to work nicely for all classical groups over
finite fields and probably for character ratios of many other elements of
interest.

\subsection{\textbf{Rank of a Representation}}

Since the 1960s, Harish-Chandra's \textit{philosophy of cusp forms} \cite%
{Harish-Chandra70} is the main organizational principle in representation
theory of reductive groups over finite and local fields. The central objects
in his approach are the cuspidal representations. It turns out that
cuspidality is a generic property, i.e., these irreps constitute a major
part of all irreps, and most of them are, in some sense, the "largest".

The philosophy of cusp forms is well adapted to establishing the Plancherel
formula for reductive groups over local fields, and leads to Lusztig's
classification \cite{Lusztig84} of the irreps of reductive groups over
finite fields.\smallskip

However, analysis of character ratios seems to require a different
approach.\smallskip

With this motivation in mind, we proposed in \cite{Gurevich-Howe15,
Gurevich-Howe17} to turn, in some sense, things upside down, and to have an
organization of the irreps of finite classical groups that is generated by
the very few "smallest" representations. As a result, representations that
may seem to be anomalies from the philosophy of cusp forms viewpoint play a
key role here. This is interesting already in the case of $SL_{2}(\mathbb{F}%
_{q})$, and this example was carried out in \cite{Gurevich-Howe18}. Although
the representations of $SL_{2}(\mathbb{F}_{q})$ have been known for a long
time, we think that the perspective of rank enhances understanding of them.

Our new organization induces several (compatible) invariants of
representations that provide strong information on the character ratios. We
call these invariants collectively \textbf{rank}.

In this note we describe parts of the development that apply to the group $%
GL_{n}(\mathbb{F}_{q}),$ and deduce from it the harmonic analytic
information we requested in Section \ref{S-HA-of-RW} for the group $SL_{n}(%
\mathbb{F}_{q})$.

In particular, for each irreducible representation $\rho $ of $GL_{n}(%
\mathbb{F}_{q})$ we attach an integer $k$ between $0$ and $n$, called its 
\textit{tensor rank}, and show, among other things, that on the transvection 
$T$ (\ref{T}) we have,\smallskip

\textbf{Theorem. }Fix $0\leq k\leq n$. Then for an irrep $\rho $ of $GL_{n}(%
\mathbb{F}_{q})$ of tensor rank $k$, we have an estimate:%
\begin{equation}
\frac{\chi _{\rho }(T)}{\dim (\rho )}=\left\{ 
\begin{array}{c}
\text{ }\frac{1}{q^{k}}+o(...)\text{, \ \ \ \ if \ \ \ }k<\frac{n}{2};\text{
\ \ \ \ \ \ \ \ \ \ \ \ } \\ 
\text{\ \ \ \ \ \ } \\ 
\frac{c_{\rho }}{q^{k}}+o(...)\text{, \ \ \ \ if \ \ }\frac{n}{2}\leq k\leq
n-1;\text{\ \ } \\ 
\text{\ } \\ 
\text{ \ }\frac{-1}{q^{n-1}-1}\text{, \ \ \ \ \ \ \ \ \ if \ \ \ }k=n,\text{
\ \ \ \ \ \ \ \ \ \ \ }%
\end{array}%
\right.  \label{CRT}
\end{equation}%
where $c_{\rho }$ is a certain integer (independent of $q$) combinatorially
associated with $\rho $.

\begin{remark}
For irreps $\rho $ of tensor rank $\frac{n}{2}\leq k\leq n-1,$ the constant $%
c_{\rho }$ in (\ref{CRT}) might be equal to zero. In this case, the estimate
on $\frac{\chi _{\rho }(T)}{\dim (\rho )}$ is simply $o(\frac{1}{q^{k}}).$
However, it is typically non-zero, and in many cases it is $1$.
\end{remark}

The estimates in (\ref{CRT}) seem to give a significant improvement to what
currently appears in the literature, and induce similar results for the
irreps of $SL_{n}(\mathbb{F}_{q})$. In particular, using some additional
analytic information, Hildebrand's Theorem \ref{H-Thm} follows.\medskip

\textbf{Acknowledgements. }The material presented in this note is based upon
work supported in part by the National Science Foundation under Grants No.
DMS-1804992 (S.G.) and DMS-1805004 (R.H.). \ 

We want to thank S. Goldstein and J. Cannon for their help with numerical
aspects of the project, part of which is reported here.

We thank J. Bernstein for sharing some of his thoughts concerning the
organization of representations by small ones.

This note was written during 2018-19 while S.G. was visiting the Math
Department and the College of Education at Texas A\&M University, the Math
Departments at Yale University and Weizmann Institute, the CS Department of
Hebrew University, and the MPI - Bonn, and he would like to thank these
institutions, and to thank personally A. Caldwell and R. Howe at TAMU, D.
Altschuler, A. Davis, Y. Minsky and C. Villano at Yale, G. Kozma, H. Naor,
D. Dvash at WI, O. Schwartz at Hebrew U, and P. Moree and C. Wels at MPI.

\tableofcontents

\section{\textbf{Character Ratios and Tensor Rank\label{S-CR-TR}}}

We start with the problem of estimating the character\footnote{%
In this note, for clarity, we denote irreps of $GL_{n}$ mostly by $\rho $
and of $SL_{n}$ mostly by $\pi .$} ratios (CRs) on the transvection $T$ (\ref%
{T}),

\begin{equation}
\frac{\chi _{\rho }(T)}{\dim (\rho )},\ \ \rho \in \widehat{GL}_{n},
\label{CR-T}
\end{equation}%
for the group $GL_{n}=GL_{n}(\mathbb{F}_{q})$ of $n\times n$ invertible
matrices with entries in $\mathbb{F}_{q}$.

\subsection{\textbf{Dimension\label{Dim-CRs-Sub}}}

\begin{figure}[h]\centering
\includegraphics
{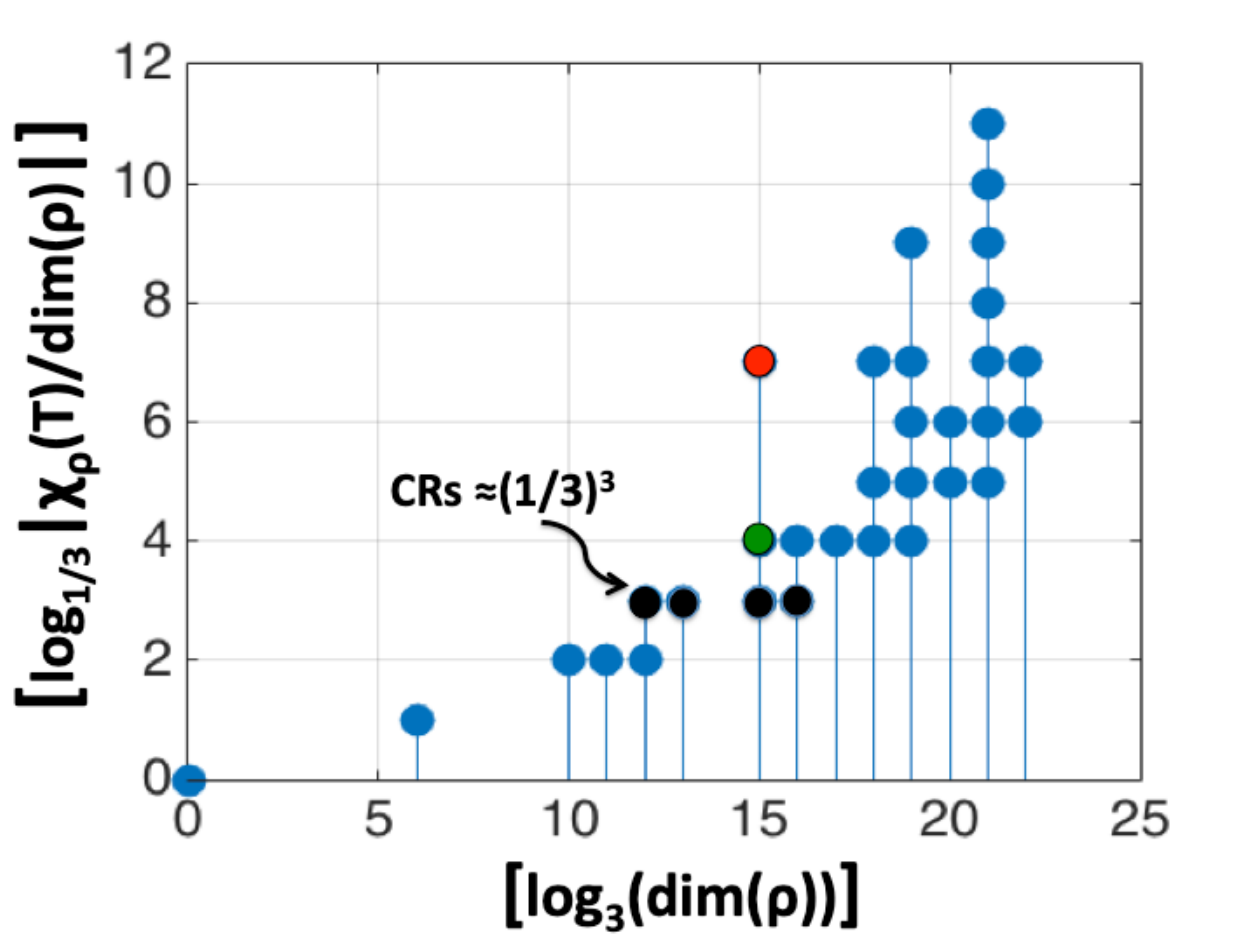}%
\caption{$\log _{1/q}$-scale of CRs\ vs. $\log _{q}$-scale of dimensions for 
$\protect\rho \in \protect\widehat{GL}_{7}(\mathbb{F}_{q}),$ $q=3.$}%
\label{cr-vs-dim-gl7_3}%
\end{figure}

At first sight one might suspect that the size of the character ratio (\ref%
{CR-T}) is to a large extent controlled by the dimension of the
representation (this is how it is usually phrased in the literature - for
example see \cite{Bezrukavnikov-Liebec-Shalev-Tiep18}) since it appears in
the denominator of (\ref{CR-T}). This is in general \textbf{not the case }%
for the transvection $T$---see Figure \ref{cr-vs-dim-gl7_3} for
illustration. In that picture, for each irreducible representation (irrep) $%
\rho $ of $GL_{7}(\mathbb{F}_{3})$ we plot\footnote{%
We denote by $\left[ x\right] $ the nearest integer to the real number $x$.}
the (nearest integer of the) absolute value in $\log _{1/3}$-scale of its
character ratio (\ref{CR-T}) vs. the (nearest integer of the) $\log _{3}$%
-scale of its dimension. In particular, one learns from this numerics that
there are (see the black circles in Figure \ref{cr-vs-dim-gl7_3}) irreps of $%
GL_{n}(\mathbb{F}_{q})$ with dimensions that differ by a multiple of large
power of $q,$ but with the same order of magnitude of CRs, and there are
(see, e.g., the black-green-red circles above 15 in Figure \ref%
{cr-vs-dim-gl7_3}) irreps of the same order of magnitude of dimension but
CRs that differ by multiple of a large power of $\frac{1}{q}$.%

A recent significant discovery \cite{Gurevich-Howe15, Gurevich-Howe17} is
that there is an invariant, different from dimension, that seems to do a
much better job in controlling the CRs (\ref{CR-T})---see Figure \ref%
{cr-vs-rank-gl7_3} for illustration. We proceed to discuss it now.

\subsection{\textbf{Tensor Rank\label{S-TR}}}

An important object attached to any finite group $G$ is its \textit{%
representation (aka Grothendieck)} \textit{ring} \cite{Zelevinsky81} 
\begin{equation*}
R(G)=%
\mathbb{Z}
\lbrack \widehat{G}],
\end{equation*}%
generated from the set $\widehat{G}$ using the operations of addition and
multiplication given, respectively, by direct sum $\oplus $ and tensor
product $\otimes .$

It turns out \cite{Gurevich17, Gurevich-Howe15, Gurevich-Howe17,
Gurevich-Howe18, Howe17-1, Howe17-2} that in the case that $G$ is a finite
classical group the ring $R(G)$ has a natural filtration that we call 
\textit{tensor rank filtration. }In particular, for each irrep we get a
non-negative integer that we call \textit{tensor rank }and might be
considered intuitively as its "size". Most importantly, this invariant seems
to nicely control analytic properties of irreps such as character ratio.

Let us describe the development in the case of $G=GL_{n}.$

Consider the permutation representation\footnote{%
Up to a sign, $\omega $ is the restriction of the oscillator representation
of $Sp_{2n}$ to $GL_{n}$ \cite{Gerardin77, Howe73-2, Weil64}.} $\omega $ of $%
GL_{n}$ on the space $L^{2}(\mathbb{F}_{q}^{n})$ of complex valued functions
on $\mathbb{F}_{q}^{n}$ given by%
\begin{equation}
\left[ \omega (g)f\text{ }\right] (x)=f(g^{-1}x),  \label{omega}
\end{equation}%
for every $g\in GL_{n},$ $f\in L^{2}(\mathbb{F}_{q}^{n}),$ and $x\in \mathbb{%
F}_{q}^{n}.$

Denote by $\widehat{GL}_{n}(\omega ^{\otimes ^{k}})$ the set of irreps of $%
GL_{n}$ that appear in $\omega ^{\otimes ^{k}}$- the $k$-fold tensor product
of $\omega $, and by $\mathbf{1}$ the trivial representation.

\begin{proposition}
\label{P-TRF}We have a sequence of proper containments%
\begin{equation}
\{\mathbf{1}\}\subsetneqq \widehat{GL}_{n}(\omega ^{\otimes
^{1}})\subsetneqq \ldots \subsetneqq \widehat{GL}_{n}(\omega ^{\otimes
^{n}})=\widehat{GL}_{n}.  \label{TRF}
\end{equation}
\end{proposition}

For a proof of \ref{P-TRF} see Appendix \ref{P-P-TRF}.\smallskip

Looking at (\ref{TRF}), we see one natural way to associate a non-negative
integer to an irrep, i.e.,

\begin{definition}[\textbf{Strict tensor rank}]
\label{D-STR}We say that an irrep $\rho $ of $GL_{n}$ is of \underline{%
strict tensor rank} $k$, if in (\ref{TRF}) its 1st occurrence is in $%
\widehat{GL}_{n}(\omega ^{\otimes ^{k}})$.
\end{definition}

We may write $\otimes $-$rank^{\star }(\rho )=k$, $r_{\otimes }^{\star
}(\rho )=k,$ or $rank_{\otimes }^{\star }(\rho )=k$, to indicate that an
irrep $\rho $ of $GL_{n}$ is of strict tensor rank $k$, and denote the set
of all such irreps by $(\widehat{GL}_{n})_{\otimes ,k}^{\star }$.

But, looking at (\ref{TRF}), there is also another way to attach a
non-negative integer to each irrep, taking into account the action of
characters (i.e., $1$-dim representations) on irreps:

\begin{definition}[\textbf{Tensor rank}]
\label{D-TR}We will say that an irrep $\rho $ of $GL_{n}$ is of \underline{%
\textbf{tensor rank}} $k$, if it is a tensor product of a character and an
irrep of strict tensor rank $k$, but not less.
\end{definition}

Again, we may use the notations $\otimes $-$rank(\rho )=k$, or $r_{\otimes
}(\rho )=k,$ or $rank_{\otimes }(\rho )=k$, to indicate that a
representation $\rho $ of $GL_{n}$ has tensor rank $k$, and denote the set
of all such irreps by $(\widehat{GL}_{n})_{\otimes ,k}$.

We extend the definition to arbitrary (not necessarily irreducible)
representation of $GL_{n}$ and say it is of tensor rank $k$ if it contains
irreps of tensor rank $k$ but not of higher tensor rank.

In particular, the \underline{\textit{tensor rank filtration}} mentioned
above is obtained by taking $F_{\otimes ,k}$ to be the elements of $R(G)$
that are sums of irreps of tensor rank less or equal to $k.$ Then, $%
F_{\otimes ,(k-1)}\subset F_{\otimes ,k}$, $F_{\otimes ,i}\otimes F_{\otimes
,j}\subset F_{\otimes ,i+j}$ for every $i,j,k,$ and $F_{\otimes ,n}=R(G).$

Sometime it is also convenient to make the following distinction and to say
that a representation of $GL_{n}$ is of \underline{low tensor rank} if it is
of tensor rank $k<\frac{n}{2}$.

We note that,

\begin{remark}
The two notions of strict tensor rank and tensor rank differ because $GL_{n}$
is not simple, and is (almost) the product of $SL_{n}$ and $\mathbb{F}%
_{q}^{\ast }$. The two notions agree on restriction to $SL_{n}$.
\end{remark}

The following example tells us how the tensor rank one and strict tensor
rank one look like, and will be vastly generalized later in Section \ref%
{S-eta-PCF}.

\begin{example}
\label{Ex-rank-k=1}The irreps of tensor rank $k=1$ of $GL_{n},$ $n\geq 2,$
are (up to twist by a character) the (non-trivial) irreducible components of 
$\omega $ (\ref{omega}). The group $GL_{1}=\mathbb{F}_{q}^{\ast }$ acts on
the space $L^{2}(\mathbb{F}_{q}^{n})$ through its action by homotheties on $%
\mathbb{F}_{q}^{n}$. For every character $\lambda $ of $\mathbb{F}_{q}^{\ast
}$ we have the $\lambda $-isotypic component $\omega _{\lambda }=\{f:\mathbb{%
F}_{q}^{n}\rightarrow 
\mathbb{C}
^{\ast }$; $f(av)=\lambda (a)f(v)$, $\ a\in \mathbb{F}_{q}^{\ast },$ $v\in 
\mathbb{F}_{q}^{n}\}.$ It is not difficult to see using direct calculations
that,

\begin{enumerate}
\item For $\lambda \neq \mathbf{1}$ the space $\omega _{\lambda }$ is
irreducible as a $GL_{n}$-representation, it has dimension $\frac{q^{n}-1}{%
q-1}\approx q^{n-1},$ and its CR on $T$ (\ref{T}) is%
\begin{equation}
\text{\ }\frac{\chi _{\omega _{\lambda }}(T)}{\dim (\omega _{\lambda })}=%
\frac{q^{n-1}-1}{q^{n}-1}\approx \frac{1}{q}  \label{CR-TRone}
\end{equation}

\item The space $\omega _{\mathbf{1}}^{o}=\{f\in \omega _{\mathbf{1}};$ $%
f(0)=0$ and $\tsum\limits_{v\in \mathbb{F}_{q}^{n}}f(v)=0\}$ is irreducible
as a $GL_{n}$-representation, it has dimension $\frac{q^{n}-q}{q-1}\approx
q^{n-1},$ and its CR on $T$ is%
\begin{equation*}
\frac{\chi _{\omega _{\mathbf{1}}^{o}}(T)}{\dim (\omega _{\mathbf{1}}^{o})}=%
\frac{q^{n-2}-1}{q^{n-1}-1}\approx \frac{1}{q}.
\end{equation*}
\end{enumerate}

In particular, one deduces that there are roughly $q^{2}$ irreps of $\otimes 
$-rank $k=1.$
\end{example}

\begin{remark}
In the case of the group $GL_{2},$ using the terminology of the "philosophy
of cusp forms" \ \cite{Harish-Chandra70}, we have,%
\begin{equation}
(\widehat{GL}_{2})_{\otimes ,0}=\text{characters, }(\widehat{GL}%
_{2})_{\otimes ,1}=\text{principal series, }(\widehat{GL}_{2})_{\otimes ,2}=%
\text{cuspidals.}  \label{GL2}
\end{equation}
\end{remark}

\subsection{\textbf{Intrinsic Characterization of Strict Tensor Rank and
Tensor Rank}}

Definitions \ref{D-TR} and \ref{D-STR} of, respectively, tensor rank and
strict tensor rank, are not intrinsic as they use the representation $\omega 
$ (\ref{omega}). At various places of this note, it will be useful for us to
use the following intrinsic characterization (given in \cite{Gurevich-Howe17}%
) of these notions.

For $0\leq k\leq n$, consider the subgroup $H_{k}\subset GL_{n}$ of elements
that pointwise fix the first $k$-coordinates subspace in $\mathbb{F}_{q}^{n}$%
, i.e., 
\begin{equation*}
H_{k}=\left\{ 
\begin{pmatrix}
I_{k} & \ast \\ 
0 & A_{n-k}%
\end{pmatrix}%
;\text{ }A_{n-k}\in GL_{n-k}\right\} .
\end{equation*}%
Note that $H_{0}=GL_{n}$, $H_{n}=\{1\}$, and $H_{k}\subset H_{k-1}$, for
every $k=1,...,n$.

In \cite{Gurevich-Howe17} we observed that,

\begin{proposition}[\textbf{Intrinsic characterisation}]
\label{P-ID-TR}A representation $\rho \in \widehat{GL}_{n}$ is of tensor
rank $k$ (respectively, strict tensor rank $k$) if and only if it admits an
eigenvector (respectively, invariant vector) for $H_{k},$ but not for $%
H_{k-1}.$
\end{proposition}

\subsection{\textbf{Numerics\label{Numerics-CRvsRank-Sub}}}

In this note, we will think on tensor rank as a formal notion of size of a
representation. But, is it going to do a good job in controlling the CRs on
the transvection (\ref{T})?%
\begin{figure}[h]\centering
\includegraphics
{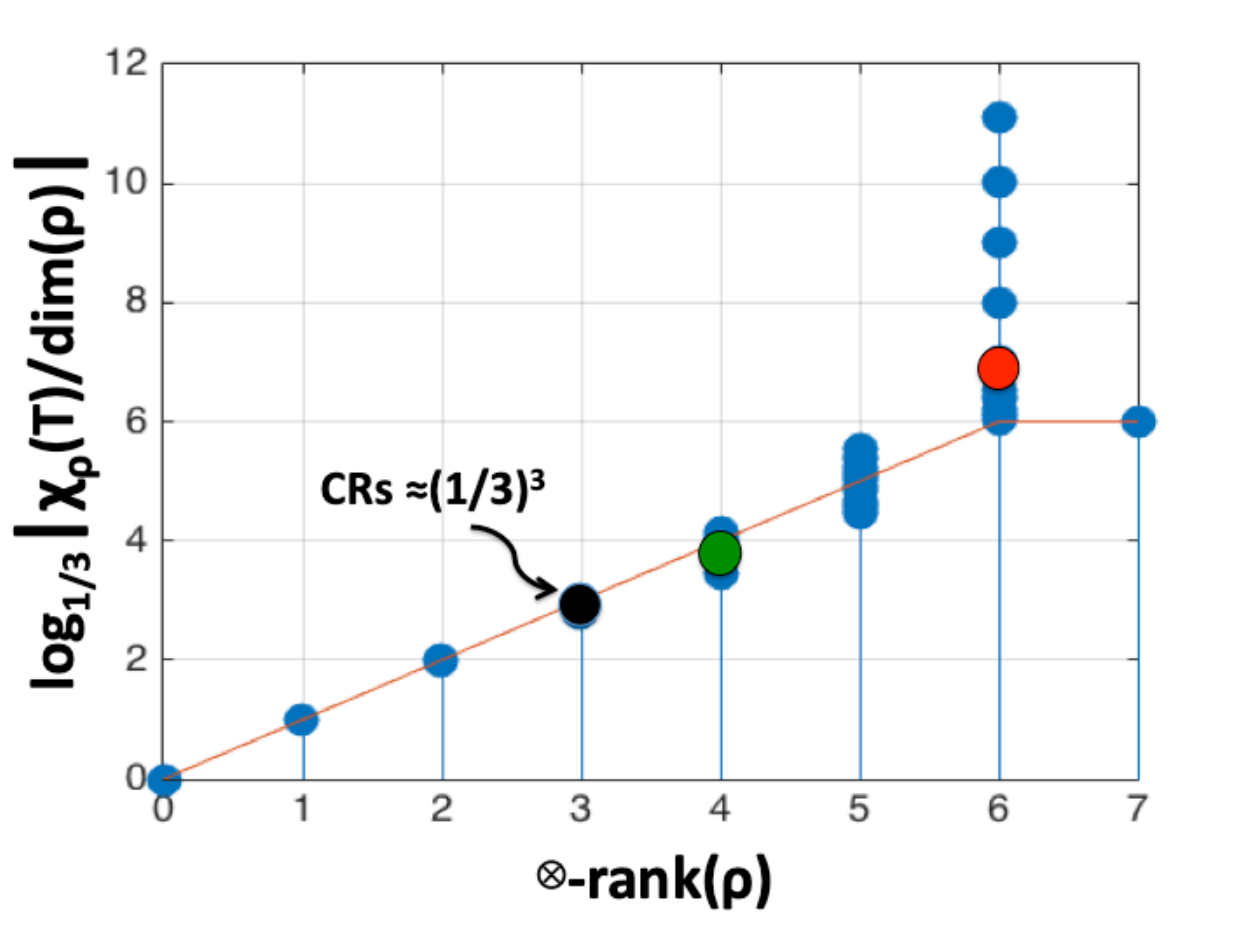}%
\caption{$\log _{1/q}$-scale of CRs\ vs. $\otimes $-rank for $\protect\rho %
\in \protect\widehat{GL}_{7}(\mathbb{F}_{q}),$ $q=3.$}%
\label{cr-vs-rank-gl7_3}%
\end{figure}%

At this stage let us present numerical data collected for the group $GL_{7}(%
\mathbb{F}_{3})$ that hints toward a positive answer to the above question.

Indeed, a comparison of Figures \ref{cr-vs-rank-gl7_3} and \ref%
{cr-vs-dim-gl7_3} indicates that the tensor rank of a representation does a
much better job than dimension in telling what should be expected for the
order of magnitude of the CRs on the transvection $T$. Indeed, Figures \ref%
{cr-vs-rank-gl7_3} show something from the general truth: For tensor rank $k<%
\frac{n}{2}$ (i.e., the low tensor rank) irreps, although the dimensions
might differ by a factor of a large power of $q,$ all the CRs are
essentially of the same size $\frac{1}{q^{k}}$ (compare the black circles in
both figures); Moreover, for higher rank $\frac{n}{2}\leq k\leq n-1,$ the
CRs are of the order of magnitude of $\frac{1}{q^{k}}$ time a constant
(independent of $q$), and it seems that for all tensor rank $n$ irreps the
CRs are exactly $\frac{1}{q^{n-1}}$ in absolute value; Finally, irreps of
the same dimensions can have different character ratios (compare the
black-green-red circles above 15 in Figure \ref{cr-vs-dim-gl7_3} with how
they appear in Figure \ref{cr-vs-rank-gl7_3} ) which are accounted for by
looking at tensor rank.\smallskip

The above numerical results can be quantified precisely and proved. This is
part of what we do next.

\section{\textbf{Analytic Information on Tensor Rank }$k$\textbf{\ Irreps of 
}$GL_{n}$\label{S-AI-GLn}}

In this section we present information concerning the character ratios and
dimensions of the irreps of $\otimes $-rank $k,$ i.e., the members of $(%
\widehat{GL}_{n})_{\otimes ,k},$ including the cardinality of that set.

\subsection{\textbf{Character Ratios on the Transvection}}

For the CRs on the transvection $T$ (\ref{T}) we obtain the following,
essentially sharp, estimate in term of the tensor rank.

\begin{theorem}
Fix $0\leq k\leq n$. Then, for $\rho \in (\widehat{GL}_{n})_{\otimes ,k},$
we have an estimate:

\begin{equation}
\frac{\chi _{\rho }(T)}{\dim (\rho )}=\left\{ 
\begin{array}{c}
\text{ }\frac{1}{q^{k}}+o(...)\text{, \ \ \ \ if \ \ \ }k<\frac{n}{2};\text{
\ \ \ \ \ \ \ \ \ \ \ \ } \\ 
\text{\ \ \ \ \ \ } \\ 
\frac{c_{\rho }}{q^{k}}+o(...)\text{, \ \ \ \ if \ \ }\frac{n}{2}\leq k\leq
n-1;\text{\ \ } \\ 
\text{\ } \\ 
\text{ \ }\frac{-1}{q^{n-1}-1}\text{, \ \ \ \ \ \ \ \ \ \ if \ \ \ }k=n,%
\text{ \ \ \ \ \ \ \ \ \ \ }%
\end{array}%
\right.  \label{CRs-GLn}
\end{equation}%
where $c_{\rho }$ is a certain integer (independent of $q$) combinatorially
associated with $\rho $.
\end{theorem}

\begin{remark}
For irreps $\rho $ of tensor rank $\frac{n}{2}\leq k\leq n-1,$ the constant $%
c_{\rho }$ in (\ref{CRs-GLn}) might be equal to zero. In this case, the
estimate on $\frac{\chi _{\rho }(T)}{\dim (\rho )}$ is simply $o(\frac{1}{%
q^{k}}).$ However, the possibility of $c_{\rho }=0$ is fairly rare, and (at
least for $k\neq n-1$) we are not sure if it happens at all.
\end{remark}

For a derivation of Estimates (\ref{CRs-GLn}), see Section \ref{S-Der-CR-T}%
.\medskip

Note that (\ref{CRs-GLn}) is a formal validation to some of the phenomena
that Figure \ref{cr-vs-rank-gl7_3} illustrates.

\subsection{\textbf{Dimensions}}

\begin{figure}[h]\centering
\includegraphics
{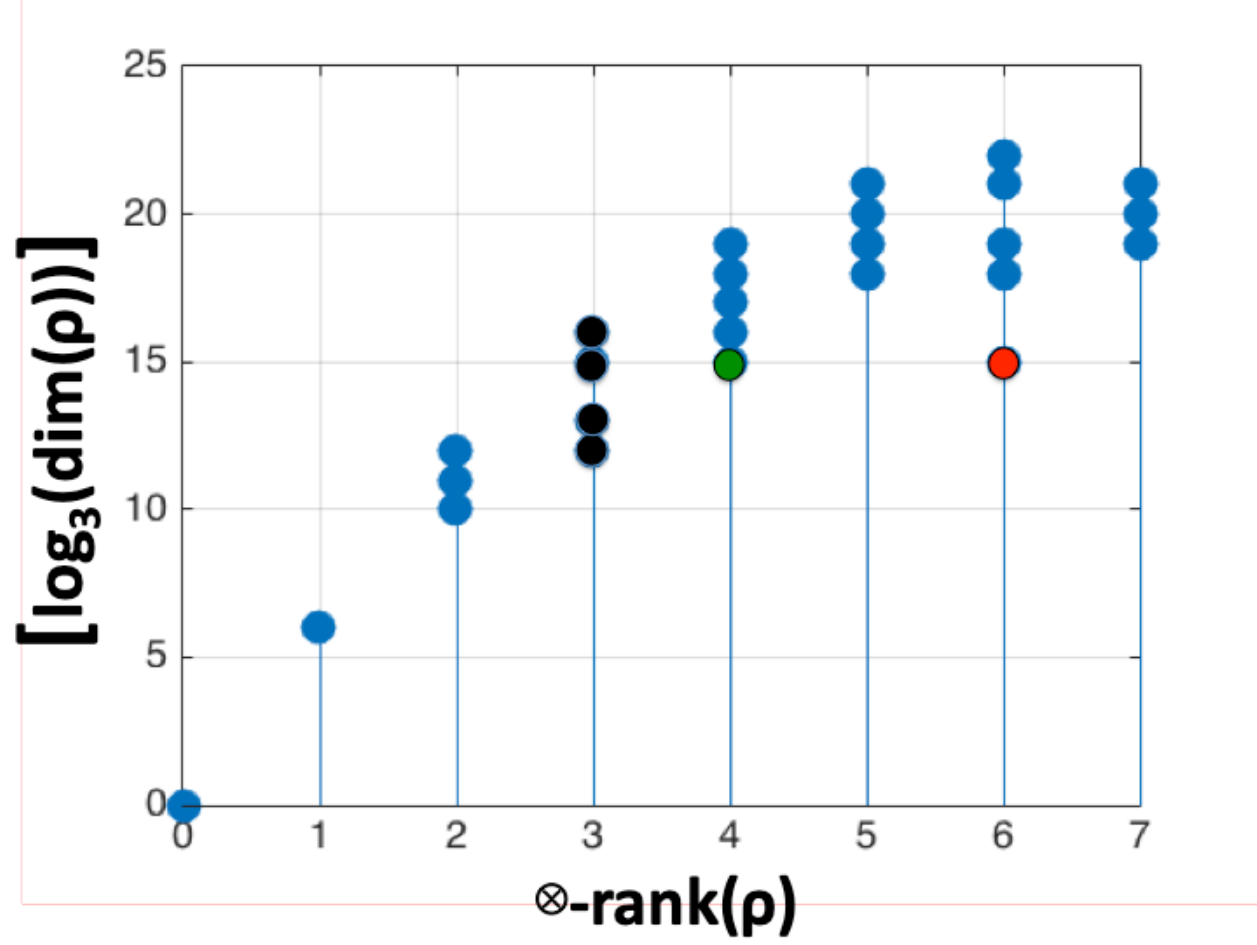}%
\caption{$\log _{q}$-scale of dimension vs. $\otimes $-rank for irreps $%
\protect\rho $ of $GL_{7}(\mathbb{F}_{q}),$ $q=3.$}\label{dim-vs-rank-gl7_3}%
\end{figure}

We proceed to present information on the dimensions of the irreps of tensor
rank $k$. Figure \ref{dim-vs-rank-gl7_3} gives a numerical illustration for
the distribution of the dimensions of the irreps of $GL_{7}(\mathbb{F}_{3})$
within each given tensor rank.%

In this note we obtain sharp lower and upper bounds (that formally explain
Figure \ref{dim-vs-rank-gl7_3}; the black-green-red dots were discussed in
Section \ref{Dim-CRs-Sub}) on the dimensions of the $\otimes $-rank $k$
irreps. Indeed, we have,

\begin{theorem}
\label{T-dim-GLn}Fix $0\leq k\leq n$. Then, for $\rho \in (\widehat{GL}%
_{n})_{\otimes ,k},$ we have an estimate:%
\begin{equation}
q^{k(n-k)+\frac{k(k-1)}{2}}+o(...)\geq \dim (\rho )\geq \left\{ 
\begin{array}{c}
q^{k(n-k)}+o(\ldots ),\text{ \ \ \ \ \ \ \ \ \ if\ \ }k<\frac{n}{2};\text{ \
\ \ \ \ \ \ \ \ \ } \\ 
\\ 
\text{ }q^{(n-k)(3k-n)}+o(\ldots ),\text{ \ \ \ if \ }\frac{n}{2}\leq k<%
\frac{2n}{3};\text{ \ \ \ } \\ 
\\ 
\text{ \ \ }q^{k(n-k)+\frac{k^{2}}{4}}+o(\ldots ),\text{ \ \ \ \ \ if \ }%
\frac{2n}{3}\leq k\leq n,\text{ even}; \\ 
\\ 
\text{\ \ \ \ \ }q^{k(n-k)+\frac{(k-3)^{2}}{4}+3(k-2)}+o(...),\text{ \ if }%
\frac{2n}{3}\leq k\leq n,\text{ odd;\ \ }%
\end{array}%
\right.  \label{Dim-GLn}
\end{equation}%
Moreover, the upper and lower bounds in (\ref{Dim-GLn}) are attained.
\end{theorem}

For a proof of Theorem \ref{T-dim-GLn} see Section \ref{Der-Est-Dim}%
.\smallskip

In \cite{Guralnick-Larsen-Tiep17}, the authors give bounds on the dimensions
of irreps of $GL_{n}$ of tensor rank $k$. However, the estimates (\ref%
{Dim-GLn}) are optimal for each tensor rank $k,$ and in general stronger
than those given in the cited paper.

\subsection{\textbf{The Number of Irreps of Tensor Rank }$k$\textbf{\ of }$%
GL_{n}$}

Finally, we present information concerning the cardinality of the set of
irreps of $\otimes $-rank $k$---see Figure \ref{log-irr(g)_k-gl7_3} for
illustration. \textbf{%
\begin{figure}[h]\centering
\includegraphics
{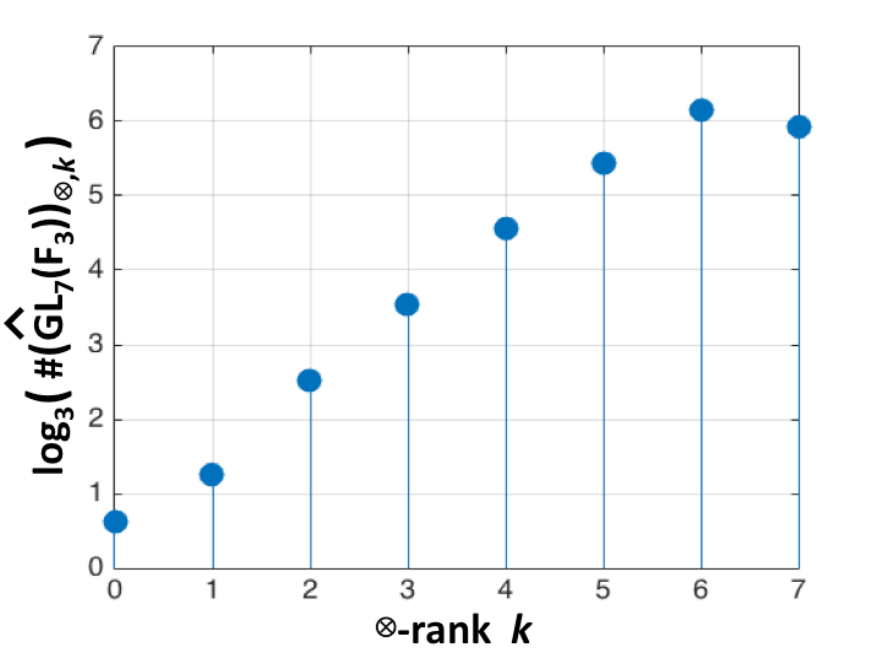}%
\caption{$\log _{q}$-scale of the number of $\otimes $-rank k irreps of $%
GL_{7}(\mathbb{F}_{q}),$ $q=3.$}\label{log-irr(g)_k-gl7_3}%
\end{figure}%
}

In this aspect, we have the following sharp estimate:

\begin{theorem}
\label{T-Card-trank-k}Fix $0\leq k\leq n$. Then, we have an estimate: 
\begin{equation}
\#((\widehat{GL}_{n})_{\otimes ,k})=\left\{ 
\begin{array}{c}
q^{k+1}+o(...)\text{, \ if \ }k\leq n-2; \\ 
c_{k}q^{n}+o(...)\text{, \ if \ }n-2<k,%
\end{array}%
\right.  \label{Card-k-GLn}
\end{equation}%
where $0<c_{n-1},$ $c_{n}<1,$ $c_{n-1}+c_{n}=1$.
\end{theorem}

For a proof of Theorem \ref{T-Card-trank-k} see Section \ref%
{S-Der-Card-trank-k}.

\subsection{\textbf{Perspective}}

We would like to make several remarks concerning the analytic information
announced just above, and to put it in some perspective to our storyline,
and to what seems to be the best known estimates in the literature on
character ratios at the transvection.

\subsubsection{\textbf{Tensor Rank vs. Dimension as Indicator for Size of
Character Ratio\label{S-TRvsDim-for-CR}}}

Looking on the analytic information presented in the sections just above, we
observe the following:\smallskip

\textbf{(A) For irreps in a given tensor rank.\smallskip }

\begin{figure}[h]\centering
\includegraphics
{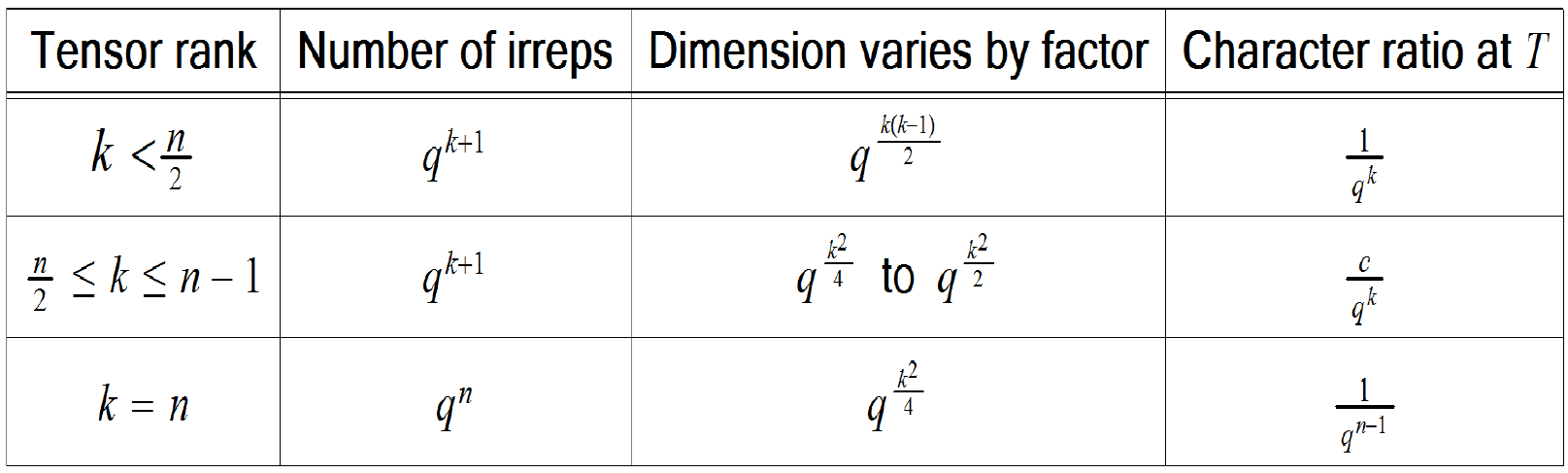}%
\caption{CRs vs. variation in dimensions (in order of magnitude of power of $%
q$) for $\protect\widehat{GL}_{n}$. }\label{cr-vs-dim-trank-k-gl7_3}%
\end{figure}

A comparison of (\ref{Dim-GLn}) and (\ref{CRs-GLn}) demonstrates---see
Figure \ref{cr-vs-dim-trank-k-gl7_3} for a summary---what we illustrated in
Sections \ref{Dim-CRs-Sub} and \ref{Numerics-CRvsRank-Sub}: Within a given
tensor rank $k$ the 
dimensions may vary by a large factor (around $q^{\frac{k(k-1)}{2}}$ for
rank $k<\frac{n}{2}$, and between $q^{\frac{k^{2}}{4}}$ to $q^{\frac{k^{2}}{2%
}}$ for $\frac{n}{2}\leq k$ - quantities are given in approximate order of
magnitude of power of $q$) but the CRs are practically the same, of size
around $\frac{1}{q^{k}}$ (for $\frac{n}{2}\leq k\leq n-1$ a multiple of $%
\frac{1}{q^{k}}$ by a constant independent of $q$).

\textbf{(B) For irreps of different tensor ranks.\smallskip }

Looking on (\ref{Dim-GLn}) we notice that:

\begin{itemize}
\item for $n>\frac{(k+1)(k+2)}{2}$, the upper bound for the dimension of $%
\otimes $-rank $k$ irreps is (for sufficiently large $q$) smaller than the
lower bound for rank $k+1$.

But,

\item when $n<\frac{(k+1)(k+2)}{2}$, the range of dimensions for $\otimes $%
-rank $k$ irreps overlaps (for large enough $q$) the range for $k+1$, and
the overlap grows with $k$. For $k$ in this range, representations of the
same dimension can have different character ratios, which are accounted for
by looking at rank.\smallskip
\end{itemize}

In conclusion, it seems that tensor rank of a representation is a better
indicator than dimension for the size of its character ratio, at least on
elements such as the transvection.

\subsubsection{\textbf{Comparison with Existing Formulations in the
Literature}}

In most of the literature on character ratios that we have seen (see, e.g., 
\cite{Bezrukavnikov-Liebec-Shalev-Tiep18} or \cite{Guralnick-Larsen-Tiep17},
and the references there), estimates on character ratios are given in terms
of the dimension of representations.

Although the dimension is a standard invariant of representations, as we
have seen in Parts (A) and (B) of Section \ref{S-TRvsDim-for-CR}, the
dimensions of representations with a given tensor rank can vary
substantially (i.e., by large powers of $q$), while the character ratio
stays more or less constant (at least for $k<\frac{n}{2}$). Thus, using only
dimension to bound character ratio will often lead to non-optimal estimates.

In particular, the estimates in this note for the character ratio on the
transvection are optimal (in term of the tensor rank), and are, in general,
stronger than the corresponding estimates in the papers cited above. For
example, for $k<\frac{n}{2}$, rather than the bound of $\frac{1}{q^{k}}$,
the paper \cite{Bezrukavnikov-Liebec-Shalev-Tiep18} gives bounds of the
order of magnitude of $\frac{q^{\frac{k(k-1)}{n-1}}}{q^{k}},$ and the
exponent $\frac{k(k-1)}{n-1}$ can be fairly large when $n$ is large and $k$
is near $\frac{n}{2}$ (the second cited paper obtained slightly weaker
bounds, on the transvection, from the first, and also formulated the result
only for irreps of tensor rank $k<\sqrt{n}$)$.$ The table in Figure \ref%
{compcrs-gln} gives some examples of the relationship between the results of
this note, and of the literature cited above.

\begin{figure}[h]\centering
\includegraphics
{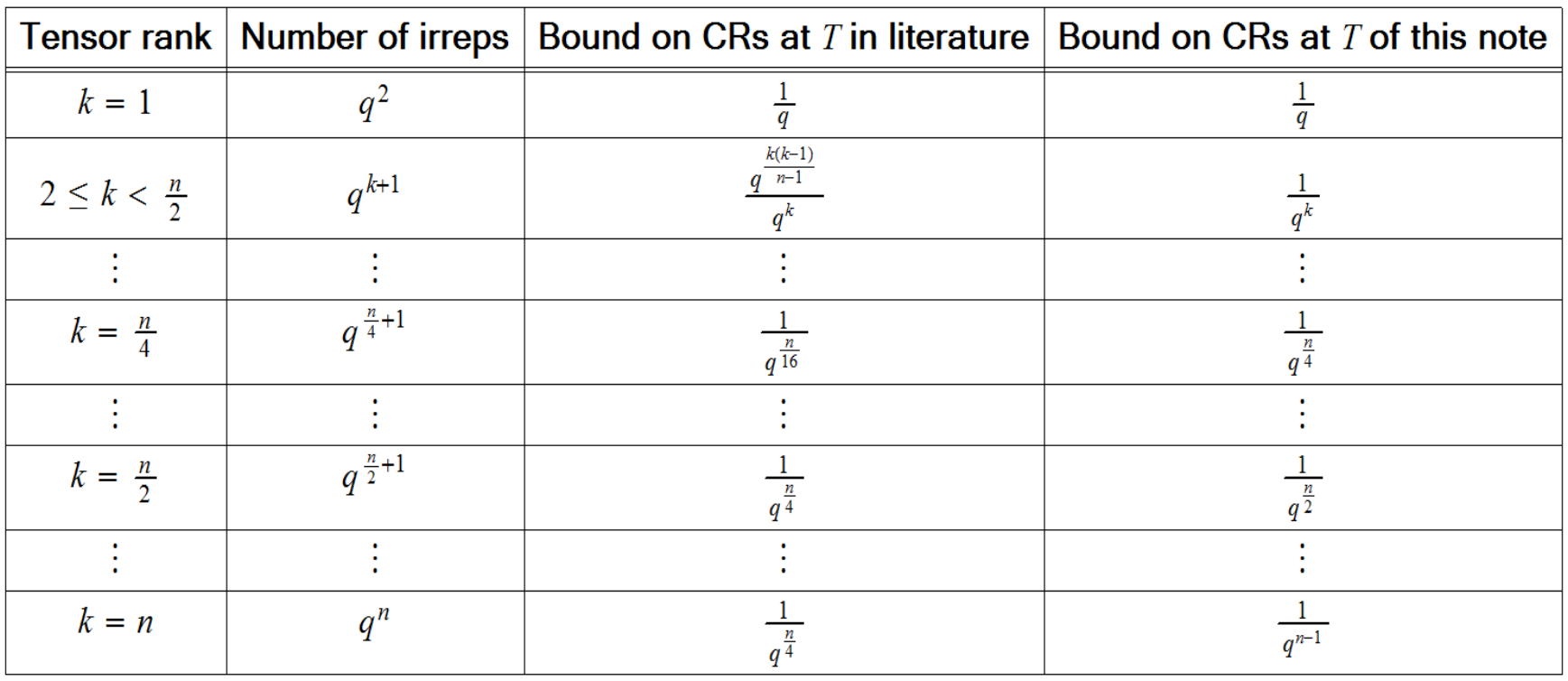}%
\caption{Bounds on CRs: Current literature vs. this note (in order of
magnitude). }\label{compcrs-gln}%
\end{figure}%

We proceed to deduce information on irreps of $SL_{n}$.

\section{\textbf{Analytic Information on Tensor Rank }$k$\textbf{\ Irreps of 
}$SL_{n}$\textbf{\label{S-AI-SLn}}}

In this section we describe analytic results for the irreps of $SL_{n},$ $%
n\geq 3,$ of a given tensor rank $k$. In some cases these estimates can be
derived as an immediate corollary of the corresponding results for $GL_{n}$,
and sometime we will need more information on the irreps of $GL_{n}$ and in
such cases we will postpone the proofs to Section \ref{S-Der-AI-SLn}. The
case of $SL_{2}$ is somewhat special - see Remark \ref{SL2} below.

\subsection{\textbf{Tensor Rank for Representations of }$SL_{n}$}

First we introduce the following terminology. We assume $n\geq 3$.

\begin{definition}
\label{D-TR-SLn}We will say that an irreducible representation $\pi $ of $%
SL_{n}$ has \underline{tensor rank} $k$ if it appears in the restriction of
a tensor rank $k$ (and not less) irrep of $GL_{n}.$
\end{definition}

As before, we denote by $(\widehat{SL}_{n})_{\otimes ,k}$ the set of irreps
of $SL_{n}$ of $\otimes $-rank $k$.\smallskip

\begin{remark}
\label{R-TR-STR-SLn}Note that the condition that $\pi $ should satisfy in
Definition \ref{D-TR-SLn} is equivalent to the requirement that (replacing $%
GL_{n}$ by $SL_{n}$) in (\ref{TRF}) it will appear in the set $\widehat{SL}%
_{n}(\omega ^{\otimes k})$ but not at earlier stage. In particular, the two
notions of strict tensor rank and tensor rank for irreps of $GL_{n}$, agree
on restriction to $SL_{n}$.
\end{remark}

Our technique to get information on irreps of $SL_{n}$ is through the way
they appear inside irreps of $GL_{n}$. Let us start with some information on
this relation.

\subsection{\textbf{Some Properties of the Restriction of Irreps from }$%
GL_{n}$\textbf{\ to }$SL_{n},$\textbf{\ }$n\geq 3$}

Take a representation $\rho $ of $GL_{n}$ and consider its restriction to $%
SL_{n}.$ We will call the set of irreps that appear in this way the $SL_{n}$%
\textit{-spectrum }of $\rho $.\textit{\ }The group $GL_{n}$ acts on $%
\widehat{SL}_{n}$ through its action by conjugation on $SL_{n}$. This in
turn induces an action of $GL_{n}$ on the $SL_{n}$-spectrum of any
representation of $GL_{n}.$

Irreducibility implies that,

\begin{claim}
The $SL_{n}$-spectrum of an irrep of $GL_{n}$ consists of a single $GL_{n}$%
-orbit.
\end{claim}

It is helpful to know that irreps of $GL_{n}$ that share the same $SL_{n}$%
-spectrum have the following simple relation:

\begin{fact}
\label{F-SLn-Spec}Two irreps of $GL_{n}$ have the same $SL_{n}$-spectrum
(equivalently share any representation of $SL_{n}$) iff they differ by a
twist by a character of $GL_{n}$.
\end{fact}

The restriction can be described more precisely as follows. For each $a\in C=%
\mathbb{F}_{q}^{\ast }$ and each $\pi \in $ $\widehat{SL}_{n},$ denote by $%
\pi _{a}$ the representation $\pi _{a}(g)=\pi (s(a)gs(a)^{-1})$, $g\in
SL_{n} $, where $s(a)\in GL_{n}$ is the diagonal matrix with $a$ in the
first entry and all other diagonal entries equal to $1$. Then,

\begin{fact}
\label{F-SLn-Mult}The restriction of an irreducible representation $\rho $
to $SL_{n}$ is multiplicity free. Moreover, for any $\pi $ in the $SL_{n}$%
-spectrum of $\rho $ we have, 
\begin{equation*}
\rho _{|SL_{n}}=\sum_{a\in C/C_{\pi }}\pi _{a},
\end{equation*}%
where $C_{\pi }$ is the stabilizer of $\pi $ in $C$.
\end{fact}

Facts \ref{F-SLn-Spec} and \ref{F-SLn-Mult} are special cases of general
results (see Corollary \ref{C-Res-N}) in Clifford-Mackey's theory \cite%
{Clifford37, Mackey49} (that we recall in Appendix \ref{A-MLGM}) on
restriction of representations from a group to general normal subgroups
(first fact), and normal subgroup with cyclic quotient (second
fact).\smallskip

We proceed to derive the estimates on the character ratios.

\subsubsection{\textbf{Character Ratios on the Transvection}}

We have the following useful Lemma:

\begin{lemma}
\label{L-CR-SLn}Any element of $SL_{n}$ whose centralizer in $GL_{n}$ maps
onto $\mathbb{F}_{q}^{\ast }$ under determinant will have the same character
ratios on any irrep of $GL_{n}$ and any irrep appearing in its restriction
to $SL_{n}.$
\end{lemma}

For a proof of Lemma \ref{L-CR-SLn} see Appendix \ref{P-L-CR-SLn}. \smallskip

Since, in the case $n\geq 3,$ the transvection $T$ (\ref{T}) meets the
conditions of Lemma \ref{L-CR-SLn}, we have, using result (\ref{CRs-GLn}),
the following sharp estimates:

\begin{corollary}
\label{C-CRs-SLn}Fix $n\geq 3$, and $0\leq k\leq n$. Then, for $\pi \in (%
\widehat{SL}_{n})_{\otimes ,k},$ we have an estimate:%
\begin{equation}
\frac{\chi _{\pi }(T)}{\dim (\pi )}=\left\{ 
\begin{array}{c}
\text{ }\frac{1}{q^{k}}+o(...)\text{, \ \ \ \ if \ \ \ }k<\frac{n}{2};\text{
\ \ \ \ \ \ \ \ \ \ \ \ } \\ 
\text{\ \ \ \ \ \ } \\ 
\frac{c_{\pi }}{q^{k}}+o(...)\text{, \ \ \ \ if \ \ }\frac{n}{2}\leq k\leq
n-1;\text{\ \ } \\ 
\text{\ } \\ 
\text{ \ }\frac{-1}{q^{n-1}-1}\text{, \ \ \ \ \ \ \ \ \ \ if \ \ \ }k=n,%
\text{ \ \ \ \ \ \ \ \ \ \ }%
\end{array}%
\right.  \label{CRs-SLn}
\end{equation}%
where $c_{\pi }$ is a certain integer (independent of $q$) combinatorially
associated with $\pi $.
\end{corollary}

\begin{remark}
For irreps $\pi $ of tensor rank $\frac{n}{2}\leq k\leq n-1,$ the constant $%
c_{\pi }$ in (\ref{CRs-SLn}) might be equal to zero. In this case, the
estimate on $\frac{\chi _{\pi }(T)}{\dim (\pi )}$ is simply $o(\frac{1}{q^{k}%
}).$
\end{remark}

\subsubsection{\textbf{Lower and Upper Bounds on Dimensions of Tensor rank }$%
k$\textbf{\ Irreps of }$SL_{n}$\textbf{\ }}

It turns out that, most irreps of $GL_{n}$ stay irreducible after
restriction to $SL_{n}$, among them all the irreps that give the lower
bounds and most of those that give the upper bounds on dimensions of tensor
rank $k$ irreps. As a consequence, from the corresponding results for $%
GL_{n} $, we obtain,

\begin{corollary}
\label{C-Dim-SLn}Fix $n\geq 3$, and $0\leq k\leq n$. Then, for $\pi \in (%
\widehat{SL}_{n})_{\otimes ,k}$, we have an estimate:%
\begin{equation}
q^{k(n-k)+\frac{k(k-1)}{2}}+o(...)\geq \dim (\pi )\geq \left\{ 
\begin{array}{c}
q^{k(n-k)}+o(\ldots ),\text{ \ \ \ \ \ \ \ \ \ if\ \ }k<\frac{n}{2};\text{ \
\ \ \ \ \ \ \ \ \ } \\ 
\\ 
\text{ }q^{(n-k)(3k-n)}+o(\ldots ),\text{ \ \ \ if \ }\frac{n}{2}\leq k<%
\frac{2n}{3};\text{ \ \ \ } \\ 
\\ 
\text{ \ \ }q^{k(n-k)+\frac{k^{2}}{4}}+o(\ldots ),\text{ \ \ \ \ \ if \ }%
\frac{2n}{3}\leq k\leq n,\text{ even}; \\ 
\\ 
\text{\ \ \ \ \ }q^{k(n-k)+\frac{(k-3)^{2}}{4}+3(k-2)}+o(...),\text{ \ if }%
\frac{2n}{3}\leq k\leq n,\text{ odd;\ \ }%
\end{array}%
\right.  \label{Dim-SLn}
\end{equation}%
Moreover, the upper and lower bounds in (\ref{Dim-SLn}) are attained.\bigskip
\end{corollary}

The detailed derivation of Corollary \ref{C-Dim-SLn} can be found in Section %
\ref{Der-Est-Dim-SLn}.

\begin{remark}
\label{SL2}Corollary \ref{C-Dim-SLn} fails for $n=2$. In that case, there
are one split principal series and one cuspidal representation of $GL_{2}$
that when restricted to $SL_{2}$ decompose, respectively, into two pieces of
dimension $\frac{q+1}{2}$, $\frac{q-1}{2}$. Moreover, for these
representations, the character ratio is of order $\frac{1}{\sqrt{q}}$. This
case, discussed in \cite{Gurevich-Howe18}, arises from the \textquotedblleft
accidental" isomorphism $SL_{2}\simeq Sp_{2},$ and these representations
should be thought of as constituents of the Weil/oscillator representation
for $Sp_{2}$, and to be the representations of tensor rank $k=1$ of this
group, while both the rest of the split principal series (dimension $=q+1$)
and the \textquotedblleft discrete series" (dimension $=q-1$) of $Sp_{2}$,
should be considered to have tensor rank $k=2$.
\end{remark}

\subsection{\textbf{The Number of Irreps of Tensor Rank }$k$\textbf{\ of }$%
SL_{n}$}

The fact, mentioned earlier, that most (in a quantified way) tensor rank $k$
irreps of $GL_{n}$ stay irreducible after restricting them to $SL_{n}$,
implies (see estimates (\ref{Card-k-GLn})) the following:

\begin{proposition}
\label{P-Num-Irr-k-SLn}Fix $n\geq 3$, and $0\leq k\leq n$. Then, we have an
estimate:%
\begin{equation}
\#((\widehat{SL}_{n})_{\otimes ,k})=\left\{ 
\begin{array}{c}
\text{ \ }q^{k}+o(...)\text{, \ \ \ if \ }k\leq n-2; \\ 
c_{k}q^{k-1}+o(...)\text{, \ if \ }n-2<k,%
\end{array}%
\right.  \label{Card-SLn}
\end{equation}%
where $0<c_{n-1},$ $c_{n}<1,$ with $c_{n-1}+c_{n}=1$.
\end{proposition}

For a detailed derivation of (\ref{Card-SLn}) see Section \ref%
{Der-P-Num-Irr-k-SLn}.

\section{\textbf{Back to the Random Walk\label{S-BRW}}}

Having the analytic information on the irreps of $SL_{n}$, $n\geq 3$, we can
address the random walk problem discussed in the introduction.

\subsection{\textbf{Setting}}

Recall (see the introduction) that we consider the conjugacy class $C\subset
SL_{n}$ of the transvection $T$ (\ref{T}), and use it, as a generating set,
to do a random walk on $SL_{n}$. We denote by $P_{C}^{\ast l}(g)$ the
probability that after $l$ steps we arrive to a given element $g\in SL_{n}$.

We know that the difference of $P_{C}^{\ast l}$ from the uniform
distribution $U$ is, in total variation, 
\begin{equation*}
\left\Vert P_{C}^{\ast l}-U\right\Vert _{TV}\approx 1,\text{ \ \ \ for \ }%
l<n,
\end{equation*}%
and want to show that the mixing time $l_{M}$ is $n,$ i.e., there is a
dramatic change at the $n$-th step where suddenly the two distributions
become close, and an exponential rate of decay---called mixing rate and
denoted $r_{M}$---kicks in.

\subsection{\textbf{The Mixing Time and Mixing Rate}}

We can derive the following sharp estimates for $l_{M}$ and $r_{M}$:

\begin{theorem}
\label{T-RW}For the random walk on $SL_{n}$ using $C$, as a generating set,
we have, for sufficiently large $q$,

\begin{enumerate}
\item The mixing time $l_{M}=n.$

\item The mixing rate $r_{M}=\frac{1}{q}+O(\frac{1}{q^{n}}).$
\end{enumerate}
\end{theorem}

Part 2 of Theorem \ref{T-RW} follows from the fact that 
\begin{equation*}
r_{M}=\underset{\mathbf{1}\neq \pi \in \widehat{SL}_{n}}{\max }\left\vert 
\frac{\chi _{\pi }(T)}{\dim (\pi )}\right\vert ,
\end{equation*}%
and then use the estimates (\ref{CR-TRone}) and (\ref{CRs-SLn}).

For a proof of Part 1 of Theorem \ref{T-RW}, first we recall that Formula (%
\ref{n}) implies that $l_{M}$ can not be less than $n,$ and then we use,

\begin{proposition}
\label{P-TVB}For sufficiently large $q$ we have, 
\begin{equation*}
\left\Vert P_{C}^{\ast l}-U\right\Vert _{TV}\leq \frac{1}{2\sqrt{q}}\left( 
\frac{1}{q}\right) ^{l-n}+o(...).
\end{equation*}
\end{proposition}

For a verification of Proposition \ref{P-TVB} see Appendix \ref{P-P-TVB}.

\section{\textbf{The eta Correspondence and the Philosophy of Cusp Forms 
\label{S-eta-PCF}}}

To derive the analytic results that we described in Section \ref{S-AI-GLn},
we need to address the following:\smallskip

\textbf{Question: }How to get information on the $\otimes $-rank $k$ irreps
of $GL_{n}$?\smallskip

In this note we would like to describe a technique which leads to an answer
to the above question and is based on the interplay between two methods: the 
\textit{Philosophy of Cusp Forms (P-of-CF)}; and the\textit{\ eta
Correspondence}.\medskip

The P-of-CF was put forward in the 60s by Harish-Chandra \cite%
{Harish-Chandra70}. It is one of the main organizing principles in
representation theory of reductive groups over local \cite%
{Berenstein-Zelevinsky76} and finite fields \cite{Deligne-Lusztig76,
Lusztig84}.\medskip

The eta correspondence was implicitly discovered in the manuscript \cite%
{Howe73-1}. This method can be applied in order to investigate irreps of
classical groups over local and finite fields. It is based on the notion of
dual pair \cite{Howe89} of subgroups in a finite symplectic group, and a
special correspondence between certain subsets of their irreps, which is
induced by restricting the oscillator (aka Weil) representation \cite%
{Gerardin77, Howe73-2, Weil64} of the relevant symplectic group to the given
dual pair.\ In recent years we have been developing this theory much further
in order to support our theory of "size" of a representation for finite
classical groups \cite{Gurevich17, Gurevich-Howe15, Gurevich-Howe17,
Gurevich-Howe18, Howe17-1, Howe17-2} and, in particular, in order to
estimate the dimensions and character ratios of their irreps.

In this note we have refined and essentially completed the development given
in \cite{Gurevich-Howe17} for the $\eta $-correspondence for the dual pair $%
(GL_{k},GL_{n})$.\medskip

In retrospect, we note that the information we obtain on irreducible
representations of tensor rank $k,$ and strict tensor rank $k$, gives an
essentially explicit description for the set of these representations, and
for the $\eta $-correspondence. The combination of the P-of-CF with this $%
(GL_{k},GL_{n})$-duality provides a simple and effective way to identify the
strict tensor rank and tensor rank $k$ pieces inside the large
representation $\omega ^{\otimes ^{k}}$ of $GL_{n}$ that we used in Section %
\ref{S-TR} to define these notions.

For the rest of this section we assume $0\leq k\leq n$.

\subsection{\textbf{The eta Correspondence - Non-Explicit Form}}

We start with a non-explicit form of the eta correspondence.

Recall (see Section \ref{S-TR}) that the vector space $L^{2}(M_{k,n}),$ of
functions on the set of $k\times n$ matrices over $\mathbb{F}_{q},$ is a
host for all (up to tensoring with characters) irreps of $GL_{n}$ of strict
tensor rank less or equal to $k$. In Section \ref{S-TR} we denoted the
action of $GL_{n}$ on this space by $\omega ^{\otimes ^{k}}$.

Of course we have a larger group of symmetries acting on this space, i.e.,
we have a pair of commuting actions 
\begin{equation}
GL_{k}\overset{\omega _{kn}}{\curvearrowright }L^{2}(M_{k,n})\overset{\omega
_{kn}}{\curvearrowleft }GL_{n},  \label{ome-nk}
\end{equation}%
given by $\left[ \omega _{kn}(h,g)f\right] (m)=f(h^{-1}mg)$, for every $h\in
GL_{k},$ $m\in M_{k,n},$ $g\in GL_{n}$, and $f\in L^{2}(M_{k,n})$.

We will also refer to $\omega _{kn}$ (\ref{ome-nk}) as the \textit{%
oscillator representation }of $GL_{k}\times GL_{n}$.

For $0<k\leq n$, the action of $GL_{k}$ does not generate the full commutant
of $GL_{n}$ in $End(\omega _{kn}),$ and vice versa (see Example \ref%
{Ex-rank-k=1} for $k=1$). Let us look at the action of $GL_{k}\times GL_{n}$
on the smaller space 
\begin{equation}
(\omega _{kn})_{\otimes ,k}^{\star }<\omega _{kn},  \label{k-spec-omega-nk}
\end{equation}%
consisting of the (sums of) components of $\omega _{kn}$ that have strict
tensor rank exactly $k.$ On this space we do have,

\begin{theorem}
\label{T-FC}The groups $GL_{k}$ and $GL_{n}$ generate each other's full
commutant in $End((\omega _{kn})_{\otimes ,k}^{\star })$.
\end{theorem}

Let us write the decomposition of $\omega _{kn}$ into a direct sum of
isotypic components for the irreps of $GL_{k}$ as follows%
\begin{equation}
\omega _{kn}=\tsum\limits_{\tau \in \widehat{GL}_{k}}\tau \otimes \Theta
(\tau ),  \label{Iso-Dec}
\end{equation}%
where the multiplicity space $\Theta (\tau )$ is a representation of $GL_{n}$%
.

Now, the Burnside's double commutant theorem \cite{Weyl46} together with
Theorem \ref{T-FC} implies that

\begin{corollary}[\textbf{eta correspondence - non explicit form}]
Each $\Theta (\tau )$ contains at most one irreducible component $\eta (\tau
)$ of strict tensor rank $k$ (and then, it appears with multiplicity one),
in addition to irreps of lower strict tensor rank. In particular, we have a
natural bijective mapping 
\begin{equation}
\tau \longmapsto \eta (\tau )\text{,}  \label{eta1}
\end{equation}%
from a subcollection of $\widehat{GL}_{k}$ onto the set $(\widehat{GL}%
_{n})_{\otimes ,k}^{\star })$ of strict tensor \textit{rank }$k$ irreps of $%
GL_{n}$. \ 
\end{corollary}

We call the mapping (\ref{eta1}) the \underline{\textbf{eta correspondence}}
or $(GL_{k},GL_{n})$\textit{-duality}.

\begin{conclusion}
Up to twist by a character of $GL_{n},$ \textbf{all }$\otimes $\textit{-rank 
}$k$ irreps of $GL_{n}$ appear in the image of the $(GL_{l},GL_{n})$-duality
for $l=k$, and not before.
\end{conclusion}

\begin{remark}
A\ proof of Theorem \ref{T-FC} first appeared in \cite{Howe73-1} (there, the
tensor rank $k$ representations were called the "new spectrum" in some
relevant "Witt tower" associated with corresponding "oscillator
representations" of symplectic groups), and a similar treatment was given in 
\cite{Gurevich-Howe17}. The outcome in both papers is the $\eta $%
-correspondence for general dual pairs in finite symplectic groups.
\end{remark}

In this note, we will prove Theorem \ref{T-FC} as a by-product of making the
description of the correspondence (\ref{eta1}) explicit.

\subsection{\textbf{The eta Correspondence - Explicit Form}}

We want to get a good formula for $\eta (\tau )$ (\ref{eta1}), including an
explicit description of its domain in $\widehat{GL}_{k}$. In this section we
get an approximate one (see Formula (\ref{For-eta}) in Theorem \ref{T-EC}
below) showing that $\eta (\tau )$ is essentially a certain simple to write
down \textit{parabolic induction} which, in addition, can be effectively
analyzed. In particular, this description will give us in Section \ref%
{S-BEDeC}, using the P-of-CF, an exact formula for $\eta (\tau )$ (see
Equation (\ref{Form-eta-tauUSD}) of Theorem \ref{EF-TR-k}).

We fix $0\leq k\leq n$ and consider inside $GL_{n}$ the maximal parabolic
subgroup 
\begin{equation}
P_{k,n-k}=Stab_{GL_{n}}(V_{k}),  \label{Pkn}
\end{equation}%
stabilizing the $k$-dimensional subspace 
\begin{equation}
V_{k}=\left\{ 
\begin{pmatrix}
x_{1} \\ 
\vdots \\ 
x_{k} \\ 
0 \\ 
\vdots \\ 
0%
\end{pmatrix}%
;\text{ }x_{1},\ldots ,x_{k}\in \mathbb{F}_{q}\right\} \subset \mathbb{F}%
_{q}^{n}.  \label{Vk}
\end{equation}%
The group $P_{k,n-k}$ has \textit{Levi decomposition }(see Appendix \ref%
{A-PCF}), i.e., it can be written as a semi-direct product of subgroups\ 
\begin{equation*}
P_{k,n-k}=U_{k,n-k}\cdot L_{k,n-k},
\end{equation*}%
where $U_{k,n-k}$ and $L_{k,n-k\text{ }},$ called, respectively, the \textit{%
unipotent radical} and the \textit{Levi component of }$P_{k,n-k}$, are given
by 
\begin{eqnarray}
\text{(1) \ \ }U_{k,n-k} &=&\left\{ 
\begin{pmatrix}
I_{k} & B \\ 
0 & I_{n-k}%
\end{pmatrix}%
;\text{ }B\in M_{k,n-k}\right\} ,  \label{LDPkn-k} \\
\text{(2) \ \ }L_{k,n-k} &=&\left\{ 
\begin{pmatrix}
A & 0 \\ 
0 & C%
\end{pmatrix}%
;\text{ }A\in GL_{k},\text{ }C\in GL_{n-k}\right\} ,  \notag
\end{eqnarray}%
where $I_{k},I_{n-k}$, are the corresponding identity matrices.

In particular, we have a surjective homomorphism 
\begin{equation}
P_{k,n-k}\overset{pr}{\twoheadrightarrow }P_{k,n-k}/U_{k,n-k}=L_{k,n-k}%
\simeq GL_{k}\times GL_{n-k}.  \label{pr}
\end{equation}%
Now, take $\tau \in \widehat{GL}_{k}$, tensor it with the trivial
representation $\mathbf{1}_{n-k}$ of $GL_{n-k},$ and form the \textit{%
parabolic induction }(see Appendix \ref{A-PCF}),%
\begin{equation}
I_{\tau }=Ind_{P_{k,n-k}}^{GL_{n}}(\tau \otimes \mathbf{1}_{n-k}),
\label{PIkn-k}
\end{equation}%
namely, the induced representation from $P_{k,n-k}$ to $GL_{n}$ of the
pullback of $\tau \otimes \mathbf{1}_{n-k}$ from $GL_{k}\times GL_{n-k}$ via
(\ref{pr}).

It turns out that,

\begin{proposition}[\textbf{Mutiplicity one}]
\label{P-MO}Consider $\tau \in \widehat{GL}_{k}$. Then, $I_{\tau }$ (\ref%
{PIkn-k}) is multiplicity free.
\end{proposition}

We give a more informative version of this result in Section \ref%
{S-Comp-I-tau}.\smallskip

It is easy to see that $I_{\tau }<$\textbf{\ }$\Theta (\tau )$ for every $%
\tau \in \widehat{GL}_{k}$ (see Part (\ref{Part2-C-OO}) of Claim \ref{C-OO}).

The representation $I_{\tau }$ gives the "approximate formula" (this is the
meaning of Equation (\ref{For-eta}) below) for $\eta (\tau )$, that we
mentioned at the beginning of this section. More precisely,

\begin{theorem}[\textbf{eta} \textbf{correspondence - explicit form}]
\smallskip \label{T-EC}Take $\tau \in \widehat{GL}_{k},$ $k\leq n,$ and look
at the decomposition (\ref{Iso-Dec}) of $\omega _{kn}$. We have,

\begin{enumerate}
\item \label{E}\textbf{Existence. }The representation $\Theta (\tau )$
contains a strict tensor rank $k$ component if and only if $\tau $ is of
strict tensor rank \ $\geq 2k-n.$

Moreover, if the condition of Part (\ref{E}) is satisfied, then,

\item \label{U}\textbf{Uniqueness.}\textit{\ } The representation $\Theta
(\tau )$ has a unique constituent $\eta (\tau )$ of strict tensor rank $k$,
and it appears with multiplicity one.\smallskip

and,\smallskip

\item \textbf{Formula. }The constituent $\eta (\tau )$ satisfies $\eta (\tau
)<I_{\tau }<\Theta (\tau )$, and we get%
\begin{equation}
\ 
\begin{tabular}{||ccc||}
\hline\hline
&  &  \\ 
& $I_{\tau }=\eta (\tau )+\dsum\limits_{\rho }\rho ,\ $ &  \\ 
&  &  \\ \hline\hline
\end{tabular}
\label{For-eta}
\end{equation}%
where the sum is multiplicity free, and over certain irreps $\rho $ which
are of strict tensor rank less then $k$ and dimension smaller than $\eta
(\tau )$.

Finally, the mapping 
\begin{equation}
\tau \longmapsto \eta (\tau )\text{,}  \label{eta2}
\end{equation}%
gives an explicit bijective correspondence between the collection $(\widehat{%
GL}_{k})_{\otimes ,\geq 2k-n}^{\star }$ of irreps of $GL_{k}$ of strict
tensor rank $\geq 2k-n$, and the set $(\widehat{GL}_{n})_{\otimes ,k}^{\star
}$ of strict tensor rank\textit{\ }$k$ irreps of $GL_{n}$.
\end{enumerate}
\end{theorem}

Note that, indeed, Theorem \ref{T-EC} gives an explicit description of the
eta correspondence (\ref{eta1}) and hence of all members of $(\widehat{GL}%
_{n})_{\otimes ,k}^{\star }$ and (up to twist by a character) of $(\widehat{%
GL}_{n})_{\otimes ,k}$ the tensor rank $k$ irreps$.$ In particular, it
implies Theorem \ref{T-FC}.

The rest of this section is devoted to formulations, and proofs, of more
informative versions of Proposition \ref{P-MO} and Theorem \ref{T-EC}.

\begin{remark}
Parts \ref{E} and \ref{U} were also formulated and proved, using extensive
character theoretic techniques, in \cite{Guralnick-Larsen-Tiep17}. The
techniques we use in this note, which among other things produce a proof of
Theorem \ref{T-EC}, are different. They are based on the philosophy of cusp
forms and, in particular, are spectral theoretic in nature.
\end{remark}

\subsection{\textbf{Decomposing }$I_{\protect\tau }=Ind_{P_{k,n-k}}^{GL_{n}}(%
\protect\tau \otimes \mathbf{1}_{n-k})$ \textbf{and the Philosophy of Cusp
Forms\label{S-DPI}}}

Denote by $(M_{k,n})_{k}\subset M_{k,n},$ the $GL_{k}\times GL_{n}$ "open"
orbit consisting of matrices of rank equal to $k.$ We observe that

\begin{claim}
\label{C-OO}The following hold:

\begin{enumerate}
\item \label{Part1-C-OO}The strict tensor rank $k$ part $(\omega
_{kn})_{\otimes ,k}^{\star }$ (\ref{k-spec-omega-nk}) of $\omega _{kn}$ is
contained in $L^{2}((M_{k,n})_{k}).$

\item \label{Part2-C-OO}We have%
\begin{equation*}
L^{2}((M_{k,n})_{k})\simeq \sum_{\tau \in \widehat{GL}_{k}}\tau \otimes
I_{\tau },
\end{equation*}%
as a representation of $GL_{k}\times GL_{n},$ where $I_{\tau }$ is given by (%
\ref{PIkn-k}).
\end{enumerate}
\end{claim}

For a proof of Claim \ref{C-OO} see Appendix \ref{P-C-OO}. \smallskip

From Claim 5.3.1, we see that the proofs of Proposition \ref{P-MO} and
Theorem \ref{T-EC} come down to learning the decomposition of $I_{\tau }.$
Our main tool for doing this involves the description of representations
coming from the philosophy of cusp forms (P-of-CF) \cite{Harish-Chandra70}.

\subsubsection{\textbf{Recollection from the Philosophy of Cusp Forms\label%
{S-R-from-P-of-CF}}}

We recall some of the basics of the P-of-CF, that are relevant for us,
leaving a more detailed account of this theory, including relevant
references, for Appendix \ref{A-PCF}.

A representation $\kappa $ of $GL_{n}$ is called \textit{cuspidal} if it
does not contain a non-trivial fixed vector for the unipotent radical of any
parabolic subgroup stabilizing a flag in $\mathbb{F}_{q}^{n}$. Given this
definition, it is easy to show that any irrep is contained in a
representation induced from a representation of a parabolic subgroup $P$ that

\begin{itemize}
\item is trivial on the unipotent radical $U_{P}$ of $P$; and

\item is a cuspidal representation of the quotient $L_{P}=P/U_{P}$.
\end{itemize}

Note that in the case of $GL_{n},$ the group $L_{P}=P/U_{P},$ called the 
\textit{Levi component} of $P$, is a product of $GL_{m}$'s for $m\leq n$, so
that the P-of-CF provides an inductive construction of all irreps. The main
ingredients needed to carry out this construction explicitly are

\begin{enumerate}
\item knowledge of the cuspidal representations of the $GL_{m}$, $m\leq n$;
and

\item \label{Part-2}decomposing the representations induced from cuspidal
representations.
\end{enumerate}

With regard to (\ref{Part-2}), there is a very general result due to
Harish-Chandra that relates different induced-from-cuspidal representations
(see also Theorem \ref{T-HC} and Corollary \ref{C-HC}):

\begin{fact}
\label{F-HC}Two such representations $Ind_{P_{j}}^{GL_{n}}(\kappa _{j})$,
for $j=1,2,$ are either equivalent, or they are completely disjoint - they
have no irreducible constituents in common. For them to be equivalent, two
conditions must be satisfied. First, the inducing parabolics must be \textit{%
associate}, meaning that their Levi components $L_{P_{j}}$ must be
conjugate. Secondly, there must be an element $g\in GL_{n}$ that both
conjugates $L_{P_{1}}$ to $L_{P_{2}}$, and at the same time, conjugates the
representations $\kappa _{1}$ and $\kappa _{2}$ to each other. In other
words, $g(L_{P_{1}})g^{-1}=L_{2}$, and moreover, the representation $g^{\ast
}\kappa _{2}$, which sends $h$ in $L_{P_{1}}$ to $\kappa _{2}(ghg^{-1})$, is
equivalent to $\kappa _{1}$. In this way, association classes of cuspidal
representations of parabolic subgroups define a partition of the unitary
dual $\widehat{GL}_{n}$ into disjoint subsets.
\end{fact}

\begin{remark}
\label{R-T-EC}Fact \ref{F-HC} can assist in the demonstration of the
necessity statement in Part (1) of Theorem \ref{T-EC}. See Appendix \ref%
{P-NP1-T-EC} for a detailed proof.\smallskip
\end{remark}

\paragraph{\textbf{The Split and Spherical Principal Series}}

As noted already, for $GL_{n}$ any Levi component is a product of copies of
groups $GL_{m}$ for $m\leq n$. The collection $D_{P}=\{m_{j}\}$ of the sizes
of the $GL_{m_{j}}$'s factors of a Levi component define a partition of $n$.
Up to conjugation, we can assume that $P$ consists of block upper triangular
matrices, and, given this, we let $m_{j}$ be the size of the $j$-th block
from the upper left corner of the matrices. We will refer to this as the $P$%
-partition. Also, we will say that a cuspidal representation of $GL_{m}$ has 
\textit{cuspidal size} $m$. Thus, the cuspidal sizes of a cuspidal
representation of a parabolic subgroup $P$ also define a partition, the same
as the $P$-partition. Up to association in $GL_{n}$, we can arrange that the
block sizes $m_{j}$ of $P$, equivalently, cuspidal sizes of a cuspidal irrep 
$\kappa $ of $P$, decrease as $j$ increases. We then also associate to the
partition, and to $P$, a Young diagram \cite{Fulton97}, whose $j$-th row has
length $m_{j}$.

If the block sizes are all equal to $1$, then the parabolic is (conjugate
to) the Borel subgroup $B$ of upper triangular matrices. The Levi component
of $B$ is $(GL_{1})^{n}$. Since this group is abelian, all of its
irreducible representations are characters - one dimensional
representations, specified by homomorphisms into $%
\mathbb{C}
^{\ast }$. We will refer to constituents of representations induced from
characters of $B$ as the \textit{split principal series}. Constituents of
the representation induced from the trivial character $\mathbf{1}$ of $B$
will be referred to as \textit{spherical principal series} (or \textit{SPS}
for short). Any representation induced from a (one dimensional) character of
a parabolic subgroup will have constituents all belonging to the split
principal series.

There is a second partition, that permits a more refined understanding of
the split principal series, essentially reducing it to understanding the
spherical principal series. A character $\chi $ of the Borel subgroup $B$ is
given by a collection $\chi _{j}$, $1\leq j\leq n$, of characters of $GL_{1}=%
\mathbb{F}_{q}^{\ast }$, where $\chi _{j}$ is the restriction of the
character $\chi $ to the $j$-th diagonal entry of an element of $B$. Up to
association, these characters can be reordered as desired. Thus, we may
assume that all diagonal entries $j$ for which the $\chi _{j}$ are equal to
a given character of $GL_{1}$ are consecutive. Given this, we can consider a
block upper triangular parabolic subgroup such that, in each diagonal block,
the characters $\chi _{j}$ are equal, and the $\chi _{j}$'s contained in
different diagonal blocks are different. We may also assume when convenient
that the sizes of these blocks decrease from top to bottom. This associates
a well-defined partition to a given character $\chi $ of $B$.

Consider the parabolic subgroup $P_{\chi }$ defined by the blocks associated
to the character $\chi $ of $B$, as in the preceding paragraph. Let the $i$%
-th block from the top of $P_{\chi }$ be $GL_{m_{i}}$. Let $\rho $ be a
constituent (i.e., an irreducible sub-rep) of the representation of $P_{\chi
}$ induced from the character $\chi $ of $B$. Then, the philosophy of cusp
forms tells us that:\smallskip

(a) the representation $Ind_{P_{\chi }}^{GL_{n}}(\rho )$ will be
irreducible;\smallskip

(b) $\rho \simeq \otimes \rho _{i}$, where $\rho _{i}$ is a constituent of
the representation of $GL_{m_{i}}$ induced from $B\cap GL_{m_{i}}$;
and\smallskip

(c) this process gives a bijection from the constituents of $%
Ind_{B}^{P_{\chi }}(\chi )$ to the constituents of $Ind_{B}^{GL_{n}}(\chi )$%
, and this last set is the product (in the natural sense) of the sets of
constituents of the $Ind_{B\cap GL_{m_{i}}}^{GL_{m_{i}}}(\chi )$.\smallskip

Moreover, because of the way $P_{\chi }$ was defined, each representation $%
Ind_{B\cap GL_{m_{i}}}^{GL_{m_{i}}}(\chi )$ has the form%
\begin{equation*}
(\chi _{i}\circ \det )\otimes Ind_{B\cap GL_{m_{i}}}^{GL_{m_{i}}}(\mathbf{1}%
),
\end{equation*}%
where $\chi _{i}$ indicates the common character of $GL_{1}$ assigned to the
diagonal entries of $B\cap GL_{m_{i}}$, and $\det $ is the determinant
homomorphism from $GL_{m_{i}}$ to $GL_{1}$. This means that the constituents
of each $Ind_{B\cap GL_{m_{i}}}^{GL_{m_{i}}}(\chi )$ has the form $(\chi
_{i}\circ \det )\otimes \rho _{i}$, where $\rho _{i}$ is a member of the
spherical principal series for $GL_{m_{i}}$. This leads us to focus on
understanding the spherical principal series.\smallskip\ 

\paragraph{\textbf{Spherical Principal Series}}

The spherical principal series of $GL_{n}$ have been studied extensively
(see Appendices \ref{S-Spherical} and \ref{A-IrrS_l-SPS}). They can be
helpfully studied through the family of induced representations $%
Ind_{P}^{GL_{n}}(\mathbf{1)}$, for all parabolic subgroups. Up to
conjugation, it is enough to consider the parabolics that contain the Borel
subgroup $B$, i.e., the block upper triangular parabolic subgroups. Also, it
is standard that if $P$ and $P^{\prime }$ are associate parabolics, then the
representations $Ind_{P}^{GL_{n}}(\mathbf{1)}$ and $Ind_{P^{\prime
}}^{GL_{n}}(\mathbf{1)}$ are equivalent. Thus, we can select a
representative from each association class of parabolics. We do this in the
usual way, by requiring that the block sizes $m_{i}$ of the diagonal blocks $%
GL_{m_{i}}$ of $P$, listed from top to bottom, are decreasing with
increasing $i$. This again gives us a partition (our third partition) of $n$%
, with an associated Young diagram $D_{P}.$ (We note that, if all of the
above discussion on the P-of-CF is referenced to a fixed original $n$, then
the successive partitions we have been describing are partitions of parts of
the preceding partition).\smallskip

\textbf{Notation: }For the rest of this note, let us denote the set of
partitions of $n$ by $\mathcal{P}_{n},$ and the corresponding set of Young
diagrams by $\mathcal{Y}_{n}$.\smallskip

Let $P=P_{D}$ be a parabolic as above, with blocks whose sizes decrease down
the diagonal, associated to the Young diagram $D\in \mathcal{Y}_{n}.$
Consider the induced representation

\begin{equation}
I_{D}=Ind_{P_{D}}^{GL_{n}}(\mathbf{1).}  \label{I_D}
\end{equation}

All the constituents of $I_{D}$ are spherical principal series
representations. We can be somewhat more precise (see Appendix \ref%
{A-IrrS_l-SPS} for a more detailed account). Recall that the set of
isomorphism classes of representations of a group $G$ form a free abelian
semi-group (monoid) on the irreducible representations, and as such, has a
natural order structure $\leq $ given by the notion of sub-representation
(or equivalently given by dominance of all coefficients in the expression of
a given representation as a sum of irreducibles). The set of
partitions/Young diagrams also has a well-known order structure $\preceq $,
the \textit{dominance order }(see Definition \ref{D-Dom}).

We know the following facts \cite{Gurevich-Howe19, Howe-Moy86} (see also
Proposition \ref{P-Mon-Max-IL} in Appendix \ref{A-IrrS_l-SPS}):

\begin{Facts}
\label{F-SPS}Consider the representations $I_{D}$ (\ref{I_D}). We have,

\begin{enumerate}
\item The map $D\mapsto I_{D}$ is order preserving from the set $\mathcal{P}%
_{n}$ of partitions of $n$, with its reverse dominance order, to the
semigroup of spherical representations (i.e., contains non-trivial $B$%
-invariant vectors) of $GL_{n}$.

\item The representation $I_{D}$ contains a constituent $\rho _{D}$ with
multiplicity one, and with the property that it is not contained in any $%
I_{D^{\prime }}$ with $D^{\prime }\succneqq D$ in the dominance order.
\end{enumerate}
\end{Facts}

\begin{remark}
\label{rho_D-dist-by-dim}The representation $\rho _{D}$ can also be
distinguished by its dimension: it is the only constituent of $I_{D}$ whose
dimension, as a polynomial in $q$, has the same degree as the cardinality of 
$GL_{n}/P_{D}$ (see Corollary \ref{C-dim} in Appendix \ref{A-est-dim-SPS})
\end{remark}

Facts \ref{F-SPS} are parallels of similar facts for the symmetric group $%
S_{n}$ \cite{Howe-Moy86}. For a given partition $D$ of $n$, let $S_{D}$
denote the stabilizer of $D$ in $S_{n}$, and let $Y_{D}$ denote the\textit{%
Young module \cite{Ceccherini-Silberstein-Scarabotti-Tolli10}}%
\begin{equation*}
Y_{D}=Ind_{S_{D}}^{S_{n}}(\mathbf{1).}
\end{equation*}%
Also let $\sigma _{D}$ be the irrep of $S_{n}$ associated to the partition $%
D $. Then the analog of Facts \ref{F-SPS} are valid. In particular, $\sigma
_{D}$ is contained in $Y_{D}$ with multiplicity one. Moreover, for any two
partitions $D_{1}$ and $D_{2}$ of $n$, the Bruhat decomposition for $GL_{n}$ 
\cite{Borel69, Bruhat56}, i.e., that 
\begin{equation}
P_{D_{1}}\diagdown GL_{n}\diagup P_{D_{2}}\simeq S_{D_{1}}\diagdown
S_{n}\diagup S_{D_{1}},  \label{BD}
\end{equation}%
implies \cite{Howe-Moy86} that we have an equality of intertwining numbers%
\begin{equation}
\left\langle I_{D_{1}},I_{D_{2}}\right\rangle =\left\langle
Y_{D_{1}},Y_{D_{2}}\right\rangle .  \label{INI-BD}
\end{equation}%
As a consequence of the facts just mentioned above, one can show (see
Appendices \ref{A-IrrS_l-SPS} and \ref{S-Spherical}, and the reference \cite%
{Gurevich-Howe19} for more precise statement) that the description of the
spherical principal series representations of $GL_{n}$ is essentially the
same as the representation theory of the symmetric group.\smallskip\ 

\paragraph{\textbf{Split, Unsplit, and a P-of-CF\ Formula for General Irreps}%
}

In contrast to the split principal series irreps, we have the irreps that we
call \textit{unsplit}. These are the components of representations induced
from cuspidal representations of parabolics with block sizes of $2$ or
larger (i.e., no blocks of size $1$).

The general representation is gotten by combining unsplit representations
and split principal series. More precisely, given a parabolic $P_{D}\subset
GL_{n}$, let $P_{u,s}$ be the maximal parabolic subgroup, with Levi
component $GL_{u}\times GL_{s}$, where $GL_{u}$ contains all the blocks of $%
P_{D}$ of size greater than $1$, and $GL_{s}$ contains all the blocks of $%
P_{D}$ of size $1$. (We remind the reader of the convention that $P_{D}$ is
block upper triangular, with the block sizes decreasing down the diagonal).
Then a constituent $\rho _{U}$ of a representation of $GL_{u}$ induced by a
cuspidal representation of $GL_{u}\cap P_{D}$ will be an unsplit
representation of $GL_{u}$. On the other hand, $P_{D}\cap GL_{s}$ will be a
Borel subgroup of $GL_{s}$, and a constituent $\rho _{S}$ of the
representation of $GL_{s}$ induced from a cuspidal representation of $%
GL_{s}\cap P_{D}$ will be a split principal series of $GL_{s}$. Since $%
GL_{u}\times GL_{s}$ is a quotient of $P_{u,s}$, the tensor product $\rho
_{U}\otimes \rho _{S}$ will also define a representation of $P_{u,s}$, and
this representation will be irreducible. Now if we look at the induced
representation%
\begin{equation}
\rho _{U,S}=Ind_{P_{u,s}}^{GL_{n}}(\rho _{U}\otimes \rho _{S}),
\label{rho_US-1}
\end{equation}%
then, the philosophy of cusp forms tells us that

(a) $\rho _{U,S}$ is irreducible; and

(b) the map $(\rho _{U},\rho _{S})\mapsto \rho _{U,S}$ is an injection from
the relevant subsets of the unitary dual of $GL_{u}\times GL_{s}$ into the
unitary dual of $GL_{n}$; and

(c) all irreducible representations of $GL_{n}$ arise in this way (including
the cuspidal representations, which are included in the situation when $%
P_{D}=GL_{n}$).

\begin{remark}[\textbf{Uniqueness and the P-of-CF formula}]
\label{R-SS}As discussed above in Section \ref{S-R-from-P-of-CF}, the
P-of-CF tells us that the split part $\rho _{S}$, appearing in (\ref%
{rho_US-1}), is induced irreducibly from a (standard, upper triangular)
parabolic subgroup $P_{S}$ of $GL_{s}$ (corresponding to a partition $%
S=\{s_{1}\geq ...\geq s_{l}\}$ of $s$) and representation of it such that,
on each diagonal block of the parabolic the constituent representation has
the form $\rho _{S_{i}}=(\chi _{i}\circ \det )\otimes \rho _{D_{i}}$, where $%
\rho _{D_{i}}$ is a spherical principal series of the relevant $GL_{s_{i}}$%
-block of $P_{S}$, and the characters $\chi _{i}$ of $GL_{1}$ are distinct
for different blocks of $P_{S}$. Moreover, the association class of $P_{S}$,
and the inducing representations, are uniquely determined. Overall, Formula (%
\ref{rho_US-1}) can be replaced by the following more precise formula:%
\begin{equation}
\begin{tabular}{||c||}
\hline\hline
\\ 
$\rho _{U,S}=Ind_{P_{u,s_{1},...,s_{l}}}^{GL_{n}}\left( \rho _{U}\otimes %
\left[ \tbigotimes\limits_{i=1}^{l}(\chi _{i}\circ \det )\otimes \rho
_{D_{i}}\right] \right) ,$ \\ 
\\ \hline\hline
\end{tabular}
\label{rho_US-2}
\end{equation}%
where $P_{u,s_{1},...,s_{l}}$ is the standard upper triangular parabolic
with blocks of sizes $u,s_{1},...,s_{l}.$
\end{remark}

Let us call (\ref{rho_US-2}) the \textit{unsplit-split P-of-CF formula }(or%
\textit{\ parametrization}).

\subsubsection{\textbf{Decomposing }$I_{\protect\tau %
}=Ind_{P_{k,n-k}}^{GL_{n}}(\protect\tau \otimes \mathbf{1}_{n-k})$\textbf{\ 
\label{S-C-of-I}}}

We are ready to describe the components of the induced representation $%
I_{\tau }=Ind_{P_{k,n-k}}^{GL_{n}}(\tau \otimes 1_{n-k})$ given in (\ref%
{PIkn-k}), using their P-of-CF formulas. Let us start with the situation
where the representation $\tau $ on the $GL_{k}$ block is a SPS
representation.\ In this case the Pieri rule produces such a
description.\smallskip\ 

\paragraph{\textbf{The Pieri Rule}}

Consider the induced representation%
\begin{equation}
I_{\rho _{D}}=Ind_{P_{k,n-k}}^{GL_{n}}(\rho _{D}\otimes \mathbf{1}_{n-k}),
\label{I_rho_D}
\end{equation}%
where $\rho _{D}$ is the SPS of $GL_{k}$ associated to the partition $D$ of $%
n$.

The parallelism between the spherical principal series and the
representations of the symmetric group implies that,

\begin{claim}
\label{C-Mul}Consider two partitions $E$ of $n$ and $D$ of $k$. Denote by $%
\rho _{E}$ and $\rho _{D}$ the corresponding SPS representations of $GL_{n}$
and $GL_{k}$, respectively. Also, denote by $\sigma _{E}$ and $\sigma _{D}$
the corresponding irreps of $S_{n}$ and $S_{k}$, respectively. Then, we have
an equality%
\begin{equation}
\left\langle \rho _{E},I_{\rho _{D}}\right\rangle =\left\langle \sigma
_{E},I_{\sigma _{D}}\right\rangle ,  \label{Mul-L-D}
\end{equation}%
where $\left\langle \bullet ,\bullet \right\rangle =\dim (Hom(\bullet
,\bullet ))$ is the standard intertwining number, and $I_{\sigma _{D}}$
denotes the induced representation 
\begin{equation}
I_{\sigma _{D}}=Ind_{S_{k}\times S_{n-k}}^{S_{n}}(\sigma _{D}\otimes \mathbf{%
1}_{n-k}),  \label{I-sigma-D}
\end{equation}%
where the subgroup $S_{k}\times S_{n-k}$ is contained in the symmetric group 
$S_{n}$ in the standard way, and $\mathbf{1}_{n-k}$ is the trivial
representation of $S_{n-k}$.
\end{claim}

For a proof of Claim \ref{C-Mul}, see Appendix \ref{P-C-Mul}.\smallskip

In conclusion, we can replace the spectral analysis of $I_{\rho _{D}}$ (\ref%
{I_rho_D}) by that of $I_{\sigma _{D}}$ (\ref{I-sigma-D}). On the latter
representation we have a complete understanding. To spell it out, let us
recall \cite{Fulton97} that if we have Young diagrams $\widetilde{D}\in 
\mathcal{Y}_{n}$ and $D\in \mathcal{Y}_{k}$ such that $\widetilde{D}$
contains $D$, denoted $\widetilde{D}\supset D$, then by removing from $%
\widetilde{D}$ all the boxes belonging to $D$, we obtain a configuration,
denoted $\widetilde{D}-D,$ called \textit{skew-diagram}. If, in addition,
each column of\ \ $\widetilde{D}$ is at most one box longer than the
corresponding column of $D$, then we call $\widetilde{D}-D$ a \textit{%
skew-row}. With this terminology, we have \textit{\cite%
{Ceccherini-Silberstein-Scarabotti-Tolli10}},

\begin{theorem}[\textbf{Pieri rule}]
\label{T-I-tauD-Dec}Let $D\in \mathcal{Y}_{k}$. Then, the induced
representation $I_{\sigma _{D}}$ (\ref{I-sigma-D}) is a multiplicity-free
sum of irreps $\sigma _{\widetilde{D}}$ of $S_{n},$ where, the Young diagram 
$\widetilde{D}\in \mathcal{Y}_{n}$ satisfies:

\begin{enumerate}
\item $\widetilde{D}\supset D$; and

\item $\widetilde{D}-D$ is a skew-row.
\end{enumerate}
\end{theorem}

In fact, the Pieri Rule can be understood geometrically as a result about
tensor products of representations of the complex general linear groups $%
GL_{m}(%
\mathbb{C}
)$ \cite{Howe92}. In particular, in Appendix \ref{P-T-I-tauD-Dec} we give a
seemingly new proof of Theorem \ref{T-I-tauD-Dec}, by translating that
result from the $GL_{m}(%
\mathbb{C}
)$ case to the $S_{n}$ case, using the classical $S_{n}$-$GL_{m}(%
\mathbb{C}
)$ \textit{Schur }(aka\textit{\ Schur-Weyl})\textit{\ duality} \cite{Howe92,
Schur27, Weyl46}. Our approach was inspired by a remark of \textbf{Nolan
Wallach}.

Finally, we would like to remark that, nowadays the Pieri rule can be
understood as a particular case of the celebrated Littlewood-Richardson rule 
\cite{Littlewood-Richardson34, Macdonald79}, but was known \textit{\cite%
{Pieri1893} }a long time before this general result.\smallskip

As noted before, Theorem \ref{T-I-tauD-Dec} together with Identity (\ref%
{Mul-L-D}), implies the analogous "Pieri rule" for the decomposition of the
induced representation $I_{\rho _{D}}$ (\ref{I_rho_D}) for a spherical
principal series $\rho _{D}$ of $GL_{k}$.

\paragraph{\textbf{The Components of }$I_{\protect\tau }$\label{S-Comp-I-tau}%
}

Next, the components of $I_{\tau }=Ind_{P_{k,n-k}}^{GL_{n}}(\tau \otimes 
\mathbf{1}_{n-k}),$ for a general $\tau \in \widehat{GL}_{k},$ can be
obtained using the P-of-CF Formula (\ref{rho_US-2}). Indeed, take
non-negative integers $u,s_{1},...,s_{l},d,$ such that $%
u+s_{1}+...+s_{l}+d=k.$ Applying Formula (\ref{rho_US-2}) with $%
S=\{s_{1},...,s_{l},d\},$ we obtain an irreducible representation 
\begin{equation}
\tau _{U,S}=Ind_{P_{u,s_{1},...,s_{l},d}}^{GL_{k}}\left( \rho _{U}\otimes 
\left[ \tbigotimes\limits_{i=1}^{l}(\chi _{i}\circ \det )\otimes \rho
_{D_{i}}\right] \otimes \rho _{D}\right) ,  \label{tau_US}
\end{equation}%
where, $\rho _{U}$ is unsplit irrep of $GL_{u},$ the $\rho _{D_{i}}$ and $%
\rho _{D}$ are, respectively, SPS irreps of the corresponding $GL_{s_{i}}$
and $GL_{d}$ blocks. Moreover, the P-of-CF ensures that the $\tau _{U,S}$'s (%
\ref{tau_US}) produced in this way are all the irreps of $GL_{k}.$ With the
realization (\ref{tau_US}) of irreps, the Pieri rule implies that $I_{\tau
_{U,S}}=Ind_{P_{k,n-k}}^{GL_{n}}(\tau _{U,S}\otimes \mathbf{1}_{n-k})$ has
the following multiplicity free decomposition into sum of irreps%
\begin{equation}
I_{\tau _{U,S}}\simeq \dsum\limits_{\widetilde{D}%
}Ind_{P_{u,s_{1},...,s_{l},d+n-k}}^{GL_{n}}\left( \rho _{U}\otimes \left[
\tbigotimes\limits_{i=1}^{l}(\chi _{i}\circ \det )\otimes \rho _{D_{i}}%
\right] \otimes \rho _{\widetilde{D}}\right) ,  \label{I_tau_US}
\end{equation}%
where $\widetilde{D}$ runs over all Young diagrams in $\mathcal{Y}_{d+n-k}$
that satisfies conditions (1) and (2) of Theorem \ref{T-I-tauD-Dec}%
.\smallskip\ 

Note that, in particular, we obtained a proof of a more precise version of
Proposition \ref{P-MO}.

\subsection{\textbf{Computing Tensor Rank using the P-of-CF Formula}}

The Pieri rule implies, as a corollary, that one can compute the strict
tensor rank and tensor rank of a representation from its P-of-CF Formula (%
\ref{rho_US-2}), more precisely directly from its split principal series
component. To state this, and similar results, it is convenient to use the
notions of \underline{\textit{tensor co-rank}} and \underline{\textit{strict
tensor co-rank}}, by which we mean, respectively, $n$ minus the tensor rank
and $n$ minus the strict tensor rank.

\begin{corollary}
\label{C-TR-PofCF}We have,

\begin{enumerate}
\item \label{C-TR-PofCF-SPS}For a partition $D=\{d_{1}\geq ...\}$ of $n$,
the tensor co-rank of the SPS representation $\rho _{D}$ is the same as its
strict tensor co-rank and is equal to $d_{1}$, that is, the longest row of
the associated Young diagram.

\item \label{C-TR-PofCF-general}The tensor co-rank of the representation $%
\rho _{U,S}$ of $GL_{n}$ described by Formula (\ref{rho_US-2}), is the
maximum of the tensor co-ranks of the SPS representations $\rho _{D_{i}},$ $%
i=1,...,l$, that appear in description of the split part of $\rho _{U,S}$.
The strict tensor co-rank of $\rho _{U,S}$ is the strict tensor rank of the
SPS representation $\rho _{D_{i}}$ that is twisted in (\ref{rho_US-2}) by
the trivial character.
\end{enumerate}
\end{corollary}

For a proof of Corollary \ref{C-TR-PofCF} see Appendix \ref{P-C-TR-PofCF}.

\subsection{\textbf{Back to the Explicit Description of the eta
Correspondence}\label{S-BEDeC}}

The above development implies a more informative version of Theorem \ref%
{T-EC}.

Indeed, looking at Formula (\ref{I_tau_US}), we see that condition (2) of
Theorem \ref{T-I-tauD-Dec} (that the difference $\widetilde{D}-D$ be a skew
row) means that the $n-k$ boxes of $\widetilde{D}-D$ all live in different
columns of $\widetilde{D}$. In particular, $\widetilde{D}$ must contain at
least $n-k$ columns, and the only way that it can contain only $n-k$ columns
is for the boxes of $\widetilde{D}-D$ to belong to the first $n-k$ columns,
consecutively, of $\widetilde{D}$. Thus, for all representations $\tau
_{U,S} $ (\ref{tau_US}) of $GL_{k}$, of strict tensor co-rank not exceeding $%
n-k$ (which is equivalent to strict tensor rank exceeding $k-(n-k)=2k-n$),
there is exactly one constituent---let us denote it by $\eta (\tau _{U,S})$%
---of the corresponding induced representation $I_{\tau _{U,S}}$ (\ref%
{I_tau_US}) of co-rank not exceeding $n-k$. Moreover, for such $\tau _{U,S}$%
's, by the philosophy of cusp forms, $\tau _{U,S}\mapsto \eta (\tau _{U,S})$
is one-to-one correspondence.

Concluding, we have completed the proof of Theorem \ref{T-EC} (see also
Remark \ref{R-T-EC}). Moreover, we obtained an explicit form of the eta
correspondence. To write it in a pleasant way, let us express the
representation $\tau _{U,S}$ (\ref{tau_US}) of $GL_{k}$ as,%
\begin{equation}
\tau _{U,S,D}=Ind_{P_{u,s,d}}^{GL_{k}}(\rho _{U}\otimes \rho _{S}\otimes
\rho _{D}),  \label{tau_USD}
\end{equation}%
where $u,s,d$ are non-negative integers such that $u+s+d=k,$ $P_{u,s,d}$ is
the corresponding standard block diagonal parabolic with blocks of sizes $%
u,s,d$, $\rho _{U}$ is an unsplit irrep of $GL_{u}$, $S$ is a partition $%
\{s_{1}\geq ...\geq s_{l}\}$ of $s$, $\rho _{S}$ is the split representation
of $GL_{s}$ induced from $(\chi _{i}\circ \det )\otimes \rho _{D_{i}}$ on
the $GL_{s_{i}}$-block of the parabolic $P_{S}$, $i=1,...,l$ (see Formula (%
\ref{tau_US})), where $\chi _{i}$ are non-trivial and distinct, and,
finally, $\rho _{D}$ is the SPS\ representation of $GL_{d}$ associated with
a Young diagram $D\in \mathcal{Y}_{d}$.

\begin{theorem}[\textbf{eta correspondence - explicit formula}]
\label{EF-TR-k}Take $\tau _{U,S,D}\in \widehat{GL}_{k}$ of the form (\ref%
{tau_USD}), and of strict tensor rank greater or equal to $2k-n$. Then, the
unique constituent $\eta (\tau _{U,S,D})$ of strict tensor rank $k$ of $%
I_{\tau _{U,S,D}}=Ind_{P_{k,n-k}}^{GL_{n}}(\tau _{U,S,D}\otimes \mathbf{1}%
_{n-k})$, satisfies,%
\begin{equation}
\begin{tabular}{||c||}
\hline\hline
\\ 
$\eta (\tau _{U,S,D})\simeq Ind_{P_{u,s,d+n-k}}^{GL_{n}}\left( \rho
_{U}\otimes \rho _{S}\otimes \rho _{\widetilde{D}}\right) ,$ \\ 
\\ \hline\hline
\end{tabular}
\label{Form-eta-tauUSD}
\end{equation}%
where $\rho _{\widetilde{D}}$ is the SPS representation of $GL_{d+n-k}$
associated with the unique Young diagram $\widetilde{D}\in \mathcal{Y}%
_{d+n-k}$\ whose first row has length $n-k$, and the rest of its rows make $%
D $. Moreover, the mapping $\tau _{U,S,D}\longmapsto $ $\eta (\tau _{U,S,D})$
agrees with the eta correspondence (\ref{eta2}).
\end{theorem}

\section{\textbf{Deriving the Analytic Information on }$\otimes $\textbf{%
-rank }$k$\textbf{\ Irreps of }$GL_{n}$\textbf{\label{S-Der-AI-GLn}}}

The concrete descriptions of the tensor rank $k$ irreps of $GL_{n}$, given
in the previous section, will be used now to derive the analytic information
stated in Section \ref{S-AI-GLn}.

\subsection{\textbf{Deriving the Character Ratios on the Transvection\label%
{S-Der-CR-T}}}

Our analysis of the character ratios (CRs) $\frac{\chi _{\rho }(T)}{\dim
(\rho )}$, where $\rho \in \widehat{GL}_{n}$ and $T$ is the transvection (%
\ref{T}), proceeds in three steps: We start with the case when $\rho $ is of
tensor rank $n$, then we treat the case when it is a spherical principal
series (SPS) representation, and finally we combine the first two cases to
derive the CRs estimate for general tensor rank $k$ irrep.

\subsubsection{\textbf{CRs for Tensor Rank }$n$\textbf{\ Irreps\label%
{S-CR-TR-n}}}

Let $U_{n-1,1}$ be the unipotent radical of the parabolic $P_{n-1,1}$ (aka
"mirabolic" \cite{Gelfand-Kazhdan71}). This group is isomorphic to $\mathbb{F%
}_{q}^{n-1},$ and any non-identity element in this group is a transvection.
Denote by $reg_{_{U_{n-1,1}}}^{\circ }=reg_{_{U_{n-1,1}}}-$ $\mathbf{1}%
_{U_{n-1,1}}$ the regular representation of $U_{n-1,1}$ minus its trivial
representation.

We have,

\begin{proposition}
\label{P-res-n-U_1,n-1}The restriction of a tensor rank $n$ irrep $\rho $ of 
$GL_{n}$ to $U_{n-1,1}$ is a multiple of $reg_{_{U_{n-1,1}}}^{\circ }$. In
particular, the character ratio of such irrep on the transvection is equal
to the CR of $reg_{_{U_{n-1,1}}}^{\circ }$ on this element, namely, 
\begin{equation}
\frac{\chi _{\rho }(T)}{\dim (\rho )}=\frac{-1}{q^{n-1}-1}.  \label{t-rank-n}
\end{equation}
\end{proposition}

For a proof of Proposition \ref{P-res-n-U_1,n-1} see Appendix \ref%
{P-P-res-n-U_1,n-1}.

\subsubsection{\textbf{CRs for Spherical Principal Series Irreps\label%
{S-CR-SPS}}}

Consider a SPS representation $\rho _{D}$ of $GL_{n},$ where $D=\{d_{1}\geq
...\geq d_{s}\}$ is a partition of $n$. The tensor rank of $\rho _{D}$ is
equal to $n-d_{1}$ (see Corollary \ref{C-TR-PofCF}). We will show that,%
\begin{equation}
\frac{\chi _{\rho _{D}}(T)}{\dim (\rho _{D})}=\left\{ 
\begin{array}{c}
\frac{1}{q^{n-d_{1}}}+o(...)\text{, \ \ if \ }d_{1}>d_{2}\text{;} \\ 
\\ 
\frac{c_{D}}{q^{n-d_{1}}}+o(...)\text{, \ \ \ \ otherwise,}%
\end{array}%
\right. \text{,}  \label{CR-rhoN}
\end{equation}%
where $c_{D}$ is a certain integer depending only on $D$ (and not on $q$).

An effective description of the representation $\rho _{D}$ is given by
Proposition \ref{P-Mon-Max-IL} in Appendix \ref{A-SM}. In particular, it
tells us that there are integers $m_{D^{\prime }D\text{,}}$ independent of $%
q $, such that, 
\begin{equation}
\rho _{D}=I_{D}+\dsum\limits_{D^{\prime }\succneqq D}m_{D^{\prime
}D}I_{D^{\prime }},  \label{rho_N-expansion}
\end{equation}%
where $\succ $ denotes the dominance relation on partitions/Young diagrams
(see Definition \ref{D-Dom}), and $I_{D}$ (respectively $I_{D^{\prime }}$)
is the natural induced module (\ref{I_D}) attached to $D$.

\begin{remark}
The integers $m_{D^{\prime }D}$ will certainly sometimes take negative
values.
\end{remark}

We will show that the CR at $T$ of the representation $I_{D}$ on the
right-hand side of (\ref{rho_N-expansion}), can be seen as the numerical
quantity that implies estimate (\ref{CR-rhoN}). To justify this assertion,
first, recall (see Remark \ref{rho_D-dist-by-dim}) that, 
\begin{equation}
\dim (\rho _{D})=\dim (I_{D})+o(...).  \label{dim-rhoN}
\end{equation}

Second, we have,

\begin{proposition}
\label{P-CR-I_N}The character ratio of the induced representation $I_{D}$ at
the transvection satisfies,%
\begin{equation}
\frac{\chi _{I_{D}}(T)}{\dim (I_{D})}=\frac{m_{d_{1}}}{q^{n-d_{1}}}+o(...),
\label{CR-I_N}
\end{equation}%
where $m_{d_{1}}$\ is the number of times the quantity $d_{1}$ appears in $D$%
.
\end{proposition}

For a proof of Proposition \ref{P-CR-I_N} see Appendix \ref{P-P-CR-I_N}%
.\medskip

And third, we use,

\begin{proposition}
\label{P-Rel-CR-I_D}Suppose $D=\{d_{1}\geq ...\geq d_{s}\}$ is a partition
of $n$ which is strictly dominated by another partition $D^{\prime }$. Then,%
\begin{equation}
\frac{\chi _{I_{D^{\prime }}}(T)}{\dim (I_{D})}=\left\{ 
\begin{array}{c}
\text{ \ \ \ \ }o(\frac{\chi _{I_{D}}(T)}{\dim (I_{D})}),\text{ \ \ \ \ \ \
\ \ \ \ \ if \ }d_{1}>d_{2}; \\ 
\\ 
c_{D,D^{\prime }}\cdot \frac{\chi _{I_{D}}(T)}{\dim (I_{D})}+o(...),\text{ \
\ \ otherwise,}%
\end{array}%
\right.  \label{Rel-CR-I_D}
\end{equation}%
where $c_{D,D^{\prime }}$ is an explicit non-negative integer depending only
on $D$ and $D^{\prime }$ (and not on $q$).
\end{proposition}

For a proof of Proposition \ref{P-Rel-CR-I_D} see Appendix \ref%
{P-P-Rel-CR-I_D}.\smallskip

Now, we can derive (\ref{CR-rhoN}). Indeed, we have,%
\begin{eqnarray}
\frac{\chi _{\rho _{D}}(T)}{\dim (\rho _{D})} &=&\frac{\chi
_{I_{D}}(T)+\dsum\limits_{D^{\prime }\succneqq D}m_{D^{\prime }D}\chi
_{I_{D^{\prime }}}(T)}{\dim (I_{D})}+o(...)  \label{Comp-CR-SPS} \\
&&  \notag \\
&=&\left\{ 
\begin{array}{c}
\frac{1}{q^{n-d_{1}}}+o(...)\text{, \ \ if \ }d_{1}>d_{2}\text{;} \\ 
\\ 
\frac{c_{D}}{q^{n-d_{1}}}+o(...)\text{, \ \ \ \ otherwise;}%
\end{array}%
\right.  \notag
\end{eqnarray}%
where $c_{D}$ is an integer depending only on $D$ (and not on $q$), the
first equality follows from Formulas (\ref{rho_N-expansion}) and (\ref%
{dim-rhoN}), and the second equality is due to (\ref{CR-I_N}) and (\ref%
{Rel-CR-I_D}).

\subsubsection{\textbf{CRs for Tensor Rank }$k$\textbf{\ Irreps - General
Case}}

We now treat the general case.

We can reduce the estimation task to the specific cases discussed in
Sections \ref{S-CR-TR-n} and \ref{S-CR-SPS}. Indeed, for the SPS\
representations, the computation (\ref{Comp-CR-SPS}) of the character ratios
is reduced to the case of induced representations of the form $%
I_{D}=Ind_{P_{D}}^{GL_{n}}(\mathbf{1)}$ (\ref{I_D}) for a Young diagram $D$.
In the same way, the computation of character ratio on the transvection, for
general split principal series representation is reduced to the case of an
induced representation from a (one-dimensional) character of a standard
parabolic, i.e., $I_{D,\chi }=Ind_{P_{D}}^{GL_{n}}(\chi )$ where $\chi $ is
a character of $P_{D}$. But on the transvection $T$ the character $\chi $ is
trivial, so we are back in the case of $I_{D}$. In particular, using Formula
(\ref{rho_US-2}) and the standard formula \cite{Fulton-Harris91} for
computing the character of induced representations, we see that it is enough
to estimate the character ratios of irreps of the form

\begin{equation}
\rho _{U,D}=Ind_{P_{u,d}}^{GL_{n}}\left( \rho _{U}\otimes \rho _{D}\right) ,
\label{rho_UD}
\end{equation}%
where

\begin{itemize}
\item $\rho _{U}$ is an unsplit irrep of $GL_{u}$;

\item $D\in \mathcal{Y}_{d}$ is a Young diagram with longest row of length $%
d_{1}=n-k;$

and

\item there are $m_{d_{1}}$ rows in $D$ of that size, and $\rho _{D}$ the
associated SPS representation.
\end{itemize}

Denote by $G(u,n)=GL_{n}/P_{u,d}$ the Grassmannian of subspaces of dimension 
$u$ in $V=\mathbb{F}_{q}^{n}$. Its cardinality is \cite{Artin57}%
\begin{equation}
\#G(u,n)=\frac{\prod\limits_{j=u+1}^{n}(q^{j}-1)}{\prod%
\limits_{j=1}^{n-u}(q^{j}-1)}.  \label{Card-Gr}
\end{equation}%
In particular, using (\ref{Card-Gr}), we observe that, 
\begin{equation}
\frac{\#\left( G(u,n-1)\right) }{\#\left( G(u,n)\right) }=\frac{q^{n-u}-1}{%
q^{n}-1}=\frac{1}{q^{u}}+o(...).  \label{Ratio-Gs}
\end{equation}

The group $GL_{n}$ acts on the set $G(u,n),$ and the collection $%
(G(u,n))^{T},$ of elements fixed by the transvection $T,$ decomposes into a
union of two sets: 
\begin{equation}
(G(u,n))^{T}=G(u,n-1)\cup G(u-1,n-1)\text{.}  \label{G^T}
\end{equation}%
The first set in the union consists of subspaces of dimension $u$ that live
inside the kernel of $T-I,$ so we can identify it with the Grassmannian $%
G(u,n-1)$; while the second set consists with those subspaces of dimension $%
u $ containing the line $L=\func{Im}(T-I),$ so we can identify it with the
Grassmannian $G(u-1,n-1)$. Note that the two sets at the right-hand side of (%
\ref{G^T}) overlap on the set $G(u-1,n-2)$ of those $V_{u}$'s that contain $%
L $, and live inside $\ker (T-I).$

For each subspace $V_{u}$ of dimension $u$, that is fixed by $T$, we
identify $GL(V_{u})\simeq GL_{u}$ and $GL(V/V_{u})\simeq GL_{d}$. In this
way, we can think of the induced actions of $T$ on $V_{u}$ and on $V/V_{u},$
as elements $T_{u}\in GL_{u}$ and $T_{d}\in GL_{d}$, respectively. In
conclusion, we obtain that

\begin{eqnarray}
\frac{\chi _{\rho _{U,D}}(T)}{\dim (\rho _{U,D})} &=&\frac{\#\left(
G(u,n-1)\right) }{\#\left( G(u,n)\right) }\cdot \frac{\chi _{\rho
_{D}}(T_{d})}{\dim (\rho _{D})}+\frac{\#(G(u-1,n-1))}{\#(G(u,n))}\cdot \frac{%
\chi _{\rho _{U}}(T_{u})}{\dim (\rho _{U})}+o(...)  \label{Comp-CR} \\
&&  \notag \\
&=&\left\{ 
\begin{array}{c}
\text{ }\frac{1}{q^{k}}+o(...)\text{, \ \ \ \ \ if \ \ }k<\frac{n}{2};\text{
\ \ \ \ \ \ \ \ \ } \\ 
\text{\ \ \ \ \ \ } \\ 
\frac{c_{D}}{q^{k}}+o(...)\text{, \ \ if \ \ }\frac{n}{2}\leq k<n-1; \\ 
\\ 
\frac{c_{D}-1}{q^{n-1}}+o(...)\text{, \ \ if \ \ }k=n-1;\text{ \ \ \ \ \ }
\\ 
\\ 
\frac{-1}{q^{n-1}}+o(...)\text{, \ \ if \ \ }k=n,\text{ \ \ \ \ \ \ \ \ }%
\end{array}%
\right.  \notag
\end{eqnarray}%
where: The first equality is a consequence of the standard formula for
computing the character of induced reps; The little-$o$ in the first row of
Formula (\ref{Comp-CR}) comes from the overlap of the two sets on the
right-hand side of (\ref{G^T}); The second equality incorporates results (%
\ref{t-rank-n}), (\ref{CR-rhoN}) - specialized to the case $n-d_{1}=k,$ and (%
\ref{Card-Gr}), (\ref{Ratio-Gs}); Finally, note that appearance in (\ref%
{Comp-CR}) of the constant $c_{D}$, an integer depending only on $D$ (and
not on $q$) that comes from its appearance in (\ref{CR-rhoN}).

This completes the derivation of the result (\ref{CRs-GLn}) on the CRs of
the irreps of $GL_{n}$ on the transvection $T$.

\subsection{\textbf{Deriving the Estimates on Dimensions\label{Der-Est-Dim}}}

It is enough, as was in the case of the computations of the CRs just above,
to compute the dimensions of the irreps of the form (\ref{rho_UD}) where $%
\rho _{U}$ is an unsplit representation of $GL_{u}$ attached to cuspidal
datum associated with Young diagram $U\in \mathcal{Y}_{u},$ with (in case $%
u\neq 0$) rows all of which are of length greater or equal $2$, and $\rho
_{D}$ an SPS representation attached to a Young diagram $D\in \mathcal{Y}%
_{d} $.

To compute the dimension of $\rho _{U},$ we use the following \cite%
{Gel'fand70, Green55} crude approximation to the dimension of a cuspidal
representation:

\begin{proposition}
\label{P-dim-cusp}The dimension of a cuspidal representation of $GL_{u}$ is $%
q^{\frac{u(u-1)}{2}}+o(...)$.
\end{proposition}

Using Proposition \ref{P-dim-cusp}, the standard formula for dimension of
induced representation, and the explicit expression for the dimension of $%
\rho _{D}$ (see Corollary \ref{C-dim} in Appendix \ref{A-est-dim-SPS}), we
can obtain a sharp estimate on the dimension of $\rho _{U,D}$. In
particular, we can compute sharp upper and lower bounds for the dimensions
of the tensor rank $k$ irreps.

\subsubsection{\textbf{Upper Bound for the Dimensions of the Tensor Rank }$k$%
\textbf{\ Irreps}}

Let us start with tensor rank $n$ irreps.\medskip

\paragraph{\textbf{The Tensor Rank }$n$\textbf{\ Case}}

From Part (2) of Corollary \ref{C-TR-PofCF}, we learn that an irrep $\rho $
of $GL_{n}$ is of tensor rank $n$ if and only if it is an unsplit
representation of the form $\rho =\rho _{U}$, associated to some Young
diagram $U\in \mathcal{Y}_{n}$ with rows all of which are of length $\geq 2$%
, and a corresponding cuspidal datum. The philosophy of cusp forms (in
particular, Part (\ref{C-param-isobaric-rep-P2}) of Corollary \ref%
{C-param-isobaric-rep} in Appendix \ref{A-UCD}) tells us that the $\rho _{U}$%
's of maximal dimension are those where the cuspidal datum consists of
non-isomorphic cuspidal representation on the various blocks of the
corresponding Levi component. In particular, Proposition \ref{P-dim-cusp}
implies that all these $\rho _{U}$'s are of dimension $\dim (\rho _{U})=q^{%
\frac{n(n-1)}{2}}+o(...)$.

We conclude,

\begin{proposition}
\label{P-largest-dim-rank-n}The largest possible dimension of a tensor rank $%
n$ irrep of $GL_{n}$ is $q^{\frac{n(n-1)}{2}}+o(...).$\medskip
\end{proposition}

\paragraph{\textbf{The Tensor Rank }$k<n$\textbf{\ Case\label%
{Dim-Trank-k-les-n-GLn}}}

Fix $0\leq k<n$, and consider the SPS representation $\rho _{D}$ where $D$
is the partition $\{n-k,1,...,1\}$ of $n$. Formula (\ref{dim-rho_L}) implies
that 
\begin{equation}
\dim (\rho _{D})=q^{k(n-k)+\frac{k(k-1)}{2}}+o(...).  \label{L-dim}
\end{equation}

The dimension appearing in (\ref{L-dim}) is the largest possible for a
tensor rank $k$ irrep. There are several ways to justify this assertion, and
we choose to proceed with a simple construction of all the tensor rank $k$
irreps of the form $\rho _{U,D}=Ind_{P_{u,d}}^{GL_{n}}\left( \rho
_{U}\otimes \rho _{D}\right) $ that are candidates for winning the "maximal
dimension competition", and then observe that they all have dimension as in (%
\ref{L-dim}).

First, the rank constraint implies that the Young diagram $D,$ that defines
the $\rho _{D}$ datum of $\rho _{U,D},$ must have longest row of size $n-k$,
so we can assume that $d=n-k+l,$ for some $0\leq l\leq k$. In particular, if
we want to maximize the dimension of $\rho _{D},$ under this constraint, the
dominance relation on partitions tells us that, the partition $D$ must be of
the form%
\begin{equation*}
D=\{n-k,\overset{l\text{ times}}{\overbrace{1,...,1}}\}\text{.}
\end{equation*}%
Second, to maximize the dimension of the unsplit part $\rho _{U}$ of $\rho
_{U,D},$ we take it to be one of the tensor rank $u$ representation of $%
GL_{u}$ of maximal dimension constructed in the Section just above (see \ref%
{P-largest-dim-rank-n}). In particular, $\dim (\rho _{U})=q^{\frac{u(u-1)}{2}%
}+o(...)$.

Finally, with any inducing $\rho _{U}$ and $\rho _{D}$ such as these we just
described, the dimension of the corresponding induced representation $\rho
_{U,D}$ (recall that $u=k-l$) is 
\begin{eqnarray*}
\dim (\rho _{U,D}) &=&q^{(k-l)(n-k+l)}\cdot q^{\frac{(k-l)(k-l-1)}{2}}\cdot
q^{l(n-k)+\frac{l(l-1)}{2}}+o(...) \\
&=&q^{k(n-k)+\frac{k(k-1)}{2}}+o(...).
\end{eqnarray*}

In conclusion, we obtain,

\begin{proposition}
\label{P-Largest-dim}The largest possible dimension of a tensor rank $k$
irreducible representation of $GL_{n}$ is $q^{k(n-k)+\frac{k(k-1)}{2}%
}+o(...).$
\end{proposition}

Overall, this completes the verification of the assertion on upper bound on
dimensions, appearing in Theorem \ref{T-dim-GLn}.

\subsubsection{\textbf{Lower Bound for the Dimensions of the Tensor Rank }$k$%
\textbf{\ Irreps}}

Let us start with tensor rank $n$ irreps.\medskip

\paragraph{\textbf{The Tensor Rank }$n$\textbf{\ Case}}

Let as assume that $n=l\lambda $ for some $\lambda \geq 1$. Consider the
standard parabolic $P_{l\times \lambda }$ with Levi component an $l$-fold
product of $GL_{\lambda }$'s, and on each $GL_{\lambda }$ the same cuspidal
representation $\kappa _{\lambda }$. Recall (see Appendix \ref{A-UCD}) that,
such a cuspidal datum is called \textit{isobaric} and the induced
representation $Ind_{P_{_{l\times \lambda }}}^{GL_{n}}(\kappa _{\lambda
}^{\otimes ^{l}})$ is, in general (e.g., for $l>1$), reducible. Moreover, it
has a unique component $\rho _{_{l\times \lambda }}$ of smallest dimension
(see Formula (\ref{dim-isobaric-rep})) with%
\begin{equation}
\dim (\rho _{_{l\times \lambda }})=\frac{\dim (\kappa _{\lambda })^{l}\cdot
\#(GL_{n}/P_{_{l\times \lambda }})}{\#(GL_{l}(\mathbb{F}_{q^{\lambda
}})/B_{l}(\mathbb{F}_{q^{\lambda }}))},  \label{dim-isob}
\end{equation}%
where $B_{l}$ is the standard Borel subgroup in $GL_{l}$.

We look at two cases:

\begin{itemize}
\item $n$\textit{\ even:} Consider the tensor rank $n$ representation $\rho
_{_{\frac{n}{2}\times 2}}$ of $GL_{n},$ given by the recipe described above
with $\lambda =2$ and $l=\frac{n}{2}$, i.e., the constituent of $Ind_{P_{%
\frac{n}{2}\times 2}}^{GL_{n}}(\kappa _{2}^{\otimes ^{\frac{n}{2}}})$ of
smallest dimension. Then, a direct calculation, using Formula (\ref{dim-isob}%
), gives%
\begin{equation}
\dim (\rho _{\frac{n}{2}\times 2})=q^{\frac{n^{2}}{4}}+o(...).
\label{dim-trank-n-even}
\end{equation}

\item $n$\textit{\ odd:} Consider the $\otimes $-rank $n$ representation $%
\rho _{3,n-3}$ of $GL_{n},$ given by $Ind_{P_{3,n-3}}^{GL_{n}}(\kappa
_{3}\otimes \rho _{\frac{n-3}{2}\times 2}),$ where $P_{3,n-3}$ is the
standard parabolic with Levi blocks $GL_{3}$ and $GL_{n-3}$, the $\kappa
_{3} $ is a cuspidal representation of $GL_{3}$, and finally, the $\rho _{%
\frac{n-3}{2}\times 2}$ is the irrep of $GL_{n-3}$ defined in the same way
as $\rho _{\frac{n}{2}\times 2}$ above. Then, using (\ref{dim-trank-n-even})
we get,%
\begin{equation}
\dim (\rho _{3,n-3})=q^{\frac{(n-3)^{2}}{4}+3(n-2)}+o(...)\text{.}
\label{dim-trank-n-odd}
\end{equation}
\end{itemize}

In fact, optimizing using the philosophy of cusp forms and Formula (\ref%
{dim-isob}), we see that Examples (\ref{dim-trank-n-even}) and (\ref%
{dim-trank-n-odd}) give the minimizers, in the dimension aspect, among the
tensor rank $n$ irreps.

In conclusion,

\begin{proposition}
\label{P-Smallest-dim-rank-n}The smallest possible dimension of a tensor
rank $n$ irreducible representation of $GL_{n}$ is $q^{\frac{n^{2}}{4}%
}+o(...)$ if $n$ is even, and $q^{\frac{(n-3)^{2}}{4}+3(n-2)}+o(...)$ if $n$
is odd.
\end{proposition}

\paragraph{\textbf{The Tensor Rank }$k<n$\textbf{\ Cases\label%
{Dim-lower-trank-k-less-n}}}

We fix $0\leq k<n,$ and consider an irrep $\rho _{U,D}$ of $GL_{n}$ of the
form (\ref{rho_UD}), i.e., $\rho _{U,D}=Ind_{P_{u,d}}^{GL_{n}}\left( \rho
_{U}\otimes \rho _{D}\right) $, which is in addition of tensor rank $k$,
namely, the Young Diagram $D\in \mathcal{Y}_{d}$ must contains a row of
length $n-k$ and this is its longest one. Optimizing to obtain the lowest
possible dimension of such irreps, we just need to decide what to do with
the "left over" $k$ boxes. Moreover, because there is no interaction between
the unsplit and split inducing data, we just need to decide if $k$ goes to
the diagram $U\in \mathcal{Y}_{u}$ or to $D\in \mathcal{Y}_{d}$.

We divide the discussion to several cases, depending on the size and,
sometime, also the parity of $k$. \medskip

\textbf{Case }$k<\frac{n}{2}$\textbf{:} Here applying (replace $n$ by $k$
there) the numerical results (\ref{dim-trank-n-even}) and (\ref%
{dim-trank-n-odd}), we see that the winner is the SPS representation $\rho
_{D}$ of $GL_{n}$ with $D$ corresponds to the partition $\{n-k,k\}$. The
dimension is of course 
\begin{equation*}
\dim (\rho _{D})=q^{k(n-k)}+o(...).
\end{equation*}%
\medskip

\textbf{Case }$\frac{n}{2}\leq k<\frac{2n}{3}$\textbf{:} Here, again, by a
direct comparison using the numerical results (\ref{dim-trank-n-even}) and (%
\ref{dim-trank-n-odd}), we see that the lowest possible dimension is of an
SPS representation, this time $\rho _{D}$ with $D$ which is associated to
the partition $\{n-k,n-k,$ $2k-n\}$. The dimension is 
\begin{equation*}
\dim (\rho _{D})=q^{(n-k)(3k-n)}+o(\ldots ),
\end{equation*}%
using the formula in Corollary \ref{C-dim}.\medskip

\textbf{Case }$\frac{2n}{3}\leq k<n$\textbf{: }Here, the comparison shows
that the winner is the unsplit side, i.e., the irreps of tensor rank $k$ and
of lowest dimension are of the form $\rho
_{U,D}=Ind_{P_{k,n-k}}^{GL_{n}}\left( \rho _{U}\otimes \rho _{D}\right) $,
where $\rho _{D}$ is the trivial representation of $GL_{n-k}$ and $\rho _{U}$
is the tensor rank $k$ representation of $GL_{k}$ of lowest dimension. In
particular,%
\begin{equation*}
\dim (\rho _{U,D})=\left\{ 
\begin{array}{c}
\text{ \ }q^{k(n-k)+\frac{k^{2}}{4}}+o(...),\text{ \ \ \ \ \ \ \ \ \ if }k%
\text{ is even;} \\ 
q^{k(n-k)+\frac{(k-3)^{2}}{4}+3(k-2)}+o(...),\text{ \ if }k\text{ is odd;}%
\end{array}%
\right.
\end{equation*}%
Overall, this completes the verification of the assertions on lower bounds
on dimensions, appearing in Theorem \ref{T-dim-GLn}.

\subsection{\textbf{Deriving the Cardinality of the Collection of Tensor
Rank }$k$ Irreps\label{S-Der-Card-trank-k}}

To calculate the number of irreps of a given tensor rank, we will use
informations that come from the $\eta $-correspondence and the philosophy of
cusp forms.

Let us start with the largest collections.

\subsubsection{\textbf{\ The Tensor Rank }$n-1$\textbf{\ and }$n$\textbf{\
Cases}}

The cuspidal representations of $GL_{n}$ are of tensor rank $n$, and from
their construction \cite{Gel'fand70} one knows \cite{Bump04, Howe-Moy86}
that there are $aq^{n}+o(...)$ of them, for some $0<a<1$. In addition, from
the direct construction of the generic split principal series irreps, i.e.,
these induced from generic characters of the Borel (or rather its Levi
component - the diagonal torus) we know that there are $bq^{n}+o(..)$ of
them for some $0<b<1.$ But, $\#(\widehat{GL}_{n})=q^{n}+o(...)$, and the eta
correspondence implies that the number of irreps of $GL_{n}$ of tensor rank $%
\leq n-2$, is not more than $q^{n-1}+o(...)$, so we deduce that there
positive constants $c_{n-1},c_{n}$, with $c_{n-1}+c_{n}=1$ such that 
\begin{equation}
\#(\widehat{GL}_{n})_{\otimes ,n-1}=c_{n-1}q^{n}+o(...),\text{ and \ }\#(%
\widehat{GL}_{n})_{\otimes ,n}=c_{n}q^{n}+o(...).  \label{card-tr-n-n-1}
\end{equation}

\begin{remark}
We note that,

\begin{enumerate}
\item The estimates (\ref{card-tr-n-n-1}) hold also for $\#(\widehat{GL}%
_{n})_{\otimes ,n-1}^{\star }$ and $\#(\widehat{GL}_{n})_{\otimes ,n}^{\star
}$

\item It can be shown that the constants $c_{n-1}$ and $c_{n}$ are
independent of $q$.
\end{enumerate}
\end{remark}

We proceed to do the counting in the lower tensor rank cases.

\subsubsection{\textbf{The Tensor Rank }$k\leq n-2$\textbf{\ Case}}

The eta correspondence (see Theorem \ref{T-EC}) gives a bijection between $(%
\widehat{GL}_{k})_{\otimes ,\geq 2k-n}^{\star }$---the irreps of $GL_{k}$ of
strict tensor rank greater or equal to $2k-n$, and $(\widehat{GL}%
_{n})_{\otimes ,k}^{\star }$---the irreps of $GL_{n}$ of strict tensor rank $%
k$. For $k\leq n-2,$ the collection $(\widehat{GL}_{k})_{\otimes ,\geq 2k-n}$
includes the irreps of strict tensor rank $k$ and $k-1$ of $GL_{k}$, so $(%
\widehat{GL}_{n})_{\otimes ,k}^{\star }=q^{k}+o(...)$. Moreover, the
description (\ref{Form-eta-tauUSD}) tells us that these irreps of strict
tensor rank $k$ and $k-1$ of $GL_{k}$ are

\begin{itemize}
\item mapped by the eta correspondence to irreps of tensor rank $k$ of $%
GL_{n}$; and,

\item produce non-isomorphic representations upon twist by any non-trivial
character of $GL_{n}$.
\end{itemize}

As a result, overall we get 
\begin{equation*}
\#(\widehat{GL}_{n})_{\otimes ,k}=q^{k+1}+o(..).
\end{equation*}

This completes the verification of Theorem \ref{T-Card-trank-k}, and the
derivation of all the analytic properties stated in Section \ref{S-AI-GLn}.

\section{\textbf{Deriving the Analytic Information on }$\otimes $\textbf{%
-rank }$k$\textbf{\ Irreps of }$SL_{n}$\textbf{\label{S-Der-AI-SLn}}}

Let us now derive the estimates on dimensions of tensor rank $k$ irreps of $%
SL_{n},$ $n\geq 3,$ and the number of such irreps. In particular, we
complete the verification of the results announced in Section \ref{S-AI-SLn}.

We start with the dimension aspect.

\subsection{\textbf{Deriving the Estimates on Dimensions\label%
{Der-Est-Dim-SLn}}}

As we remarked earlier, the estimates in the dimension aspect are the same
as for $GL_{n}$ for a simple reason: the restrictions to $SL_{n}$ of the
irreps of $GL_{n}$ that give the lowest and the largest dimensions in a
given tensor rank, typically stay irreducible.

The main tool we use to check the irreducibility in question is the
Clifford-Mackey criterion which says (see Corollary \ref{C-Irr-Res} in
Appendix \ref{A-MLGM}) that the restriction to $SL_{n}$ of an irrep $GL_{n}$%
, stays irreducible if and only if it is not fixed by twist of any
non-trivial character of $GL_{n}$.

As was done for $GL_{n}$, we do some case by case computations.

\subsubsection{\textbf{Upper Bound on Dimensions of Tensor rank }$k$\textbf{%
\ Irreps\protect\smallskip }}

\paragraph{\textbf{Tensor Rank }$n$\textbf{\ Case\label{S-UB-Dim-trank-n}}}

The cuspidal representations of $GL_{n}$ have a parametrization by the
complex characters of the maximal anisotropic torus of $GL_{n}$ modulo the
action of the Galois group of the degree $n$ extension of the finite field 
\cite{Gel'fand70, Howe-Moy86}. In particular, there exist cuspidal
representations (in fact most of them have this property) of $GL_{n}$ which
are not fixed by any twist by a character of $GL_{n}$, and so their
restrictions to $SL_{n}$ stay irreducible and have the dimension $q^{\frac{%
n(n-1)}{2}}+o(...)$. This is the sharp upper bound announced in (\ref%
{Dim-SLn}) for tensor rank $n$ irreps.\smallskip

\paragraph{\textbf{Tensor Rank }$k<n$\textbf{\ Case}}

Consider the SPS\ representation $\rho _{D},$ associated to the partition $D$
given by $\{n-k,1,...,1\}$ of $n$. It gives (see Section \ref%
{Dim-Trank-k-les-n-GLn}) the upper bound for the dimensions of tensor rank $%
k<n$ irreps of $GL_{n}$, and stays (for example by the criterion stated just
above) irreducible after restriction to $SL_{n}$. This shows that, indeed,
in the range $k<n$ the upper bound $q^{k(n-k)+\frac{k(k-1)}{2}}+o(...)$
appearing in (\ref{Dim-SLn}) is sharp.

\subsubsection{\textbf{Lower Bound on Dimensions of Tensor rank }$k$\textbf{%
\ Irreps\protect\smallskip }}

\paragraph{\textbf{Tensor Rank }$n$\textbf{\ Case}}

Take a cuspidal representation of $GL_{2}$ which is not fixed by a twist of
any character of $GL_{2}$. In addition, take a cuspidal representation of $%
GL_{3}$ with similar property (see Section \ref{S-UB-Dim-trank-n} for more
detailed discussion). Then, apply the construction of the tensor rank $n$
irreps of $GL_{n}$ of minimal dimension. They will be irreducible after
restriction to $SL_{n},$ and have the dimension given as a sharp lower bound
in (\ref{Dim-SLn}) for tensor rank $n$ irreps.\smallskip

\paragraph{\textbf{Tensor Rank }$k<n$\textbf{\ Case}}

As in the case of $GL_{n}$ (see Section \ref{Dim-lower-trank-k-less-n}) we
go over several cases.\smallskip

\textbf{Case }$k<\frac{n}{2}$\textbf{: }Here, the representation of $GL_{n}$
with smallest possible dimension is (up to tensoring with a character) given
by the SPS\ representation $\rho _{D},$ where $D$ is the partition $%
\{n-k,k\} $ of $n$. By the Clifford-Mackey criterion it stays irreducible
after restriction to $SL_{n}$, This shows that, indeed, in the range $k<%
\frac{n}{2} $ the dimension $q^{k(n-k)}+o(...)$ appearing in (\ref{Dim-SLn})
is a lower bound and, indeed, a sharp one.\smallskip

\textbf{Case }$\frac{n}{2}\leq k<\frac{2n}{3}$\textbf{:} In this interval,
for the $GL_{n}$, the lowest possible dimension $q^{(n-k)(3k-n)}+o(\ldots )$
is (again up to tensoring by a character) of the SPS representation $\rho
_{D}$ with $D$ the partition $\{n-k,n-k,$ $2k-n\}$ of $n$. Again, by the
Clifford-Mackey criterion it stays irreducible after restriction to $SL_{n}$%
, confirming that also in this case what appear in (\ref{Dim-SLn}) is a
lower bound, and a sharp one.\smallskip \smallskip

\textbf{Case }$\frac{2n}{3}\leq k<n$\textbf{: }The same reasoning, using
what we have for $GL_{n}$ (see Section \ref{Dim-lower-trank-k-less-n}),
implies that also for this interval the estimate in (\ref{Dim-SLn}) is a
sharp lower bound.

We proceed to discuss the cardinality aspect.

\subsection{\textbf{Deriving the Number of Irreps of Tensor Rank }$k$\textbf{%
\ of }$SL_{n}$\textbf{\label{Der-P-Num-Irr-k-SLn}}}

We have sharp estimate on the number of tensor rank $k$ irreps of $GL_{n}$ -
see Formula (\ref{Card-k-GLn}). We also know (see Fact \ref{F-SLn-Spec})
that irreps of $GL_{n}$ share any constituent (hence all) after restriction
to $SL_{n}$ if and only if they differ by twist by a character of $GL_{n}$.
Finally, we can show that most tensor rank $k$ irreps of $GL_{n}$ stay
irreducible after restriction to $SL_{n}$. Indeed, we have the following
quantitative result:

\begin{proposition}
\label{P-Irr-Res-Gen}Consider the irreps of $GL_{n}$ of tensor rank $k.$
Then,

\begin{enumerate}
\item For $k<\frac{n}{2}$ all of them stay irreducible after restriction to $%
SL_{n}$.

\item \label{P2}For $\frac{n}{2}\leq k$ the proportion of them which stay
irreducible after restriction to $SL_{n}$ is $1-o(\frac{1}{q}).$
\end{enumerate}
\end{proposition}

For a proof of Proposition \ref{P-Irr-Res-Gen} see Appendix \ref%
{P-P-Irr-Res-Gen}.\smallskip

With the help of Proposition \ref{P-Irr-Res-Gen} we can get the exact
estimates that stated in Formula (\ref{Card-SLn}). We go over two
cases:\smallskip

\textbf{Case }$k<\frac{n}{2}$: Here the conclusion is clear, after
restriction, taking into account Fact \ref{F-SLn-Spec} and Formula (\ref%
{Card-k-GLn}), we get that $\#((\widehat{SL}_{n})_{\otimes
,k})=q^{k}+o(...).\smallskip $

\textbf{Case }$\frac{n}{2}\leq k$: Here, using Part (\ref{P2}) of
Proposition \ref{P-Irr-Res-Gen}, we see, again using Fact \ref{F-SLn-Spec},
that $\#((\widehat{SL}_{n})_{\otimes ,k})=q^{k}+o(...)$ for $k<n-1$, and
there are two positive constants $c_{n-1},c_{n},$ with $c_{n-1}+c_{n}=1,$
such that $\#((\widehat{SL}_{n})_{\otimes ,k})=c_{k}q^{n-1}+o(...)$, for $%
k=n-1,n$.\smallskip

This completes the derivation of estimates (\ref{Card-SLn}), and of all the
analytic properties announced in Section \ref{S-AI-SLn}.

\appendix

\section{\textbf{Clifford-Mackey Theory\label{A-MLGM}}}

We describe some parts from \textit{Clifford theory/Mackey's little group
method }\cite{Clifford37, Mackey49} that are relevant to this note.

\subsection{\textbf{Setting}}

Suppose you have a finite group $G$ which is a semi-direct product 
\begin{equation*}
G=C\ltimes N,
\end{equation*}%
where $N$ is a normal subgroup, and $C$ is cyclic.

A simple version of Clifford-Mackey theory gives the construction of the
irreps of $G$ from the irreps of $N$, and describes how irreps of $G$
decompose under restriction to $N$.

\subsection{\textbf{The Construction}}

Note that the group $C$ acts on $\widehat{N}$, the unitary dual of $N$, by
conjugation. We will call the members of $\widehat{N}$ that appear in, the
restriction to $N$, of a representation $\rho \in \widehat{G},$ the 
\underline{$N$-spectrum} of $\rho $. The irreducibility of $\rho $ implies
that,\ 

\begin{claim}
\label{Single-O}The $N$\textit{-spectrum} of $\rho \in \widehat{G}$ is a
single orbit for the action of $C$ on $\widehat{N}.$
\end{claim}

Let us construct all $\rho \in \widehat{G}$ sharing a given $N$-spectrum.
Take a representation $\pi \in \widehat{N}$, and let $C_{\pi }\subset C$ be
the stabilizer of $\pi $ in $C$. Then for each $c$ in $C_{\pi }$, there is
an operator $\sigma (c)$ on the space of $\pi $ such that%
\begin{equation}
\sigma (c)\pi (n)\sigma (c)^{-1}=\pi (cnc^{-1}),  \label{sigma_c}
\end{equation}%
and this $\sigma (c)$ is determined up to a scalar multiple, by Schur's
Lemma.

\begin{claim}
\label{sigma}We can choose the operators $\sigma (c),$ $c\in C_{\pi },$ from
(\ref{sigma_c}) in such a way that they form a representation of $C_{\pi }$.
\end{claim}

Claim \ref{sigma} follows from the fact that $C_{\pi }$ is cyclic\footnote{%
If $C_{\pi }$ is not cyclic, then it may not happen that the $\sigma (c)$
can be chosen to form a representation. The prime example is when $G$ is the
Heisenberg group, and $N$ is its center.}. Indeed, if $c_{0}$ is a generator
of $C_{\pi },$ then we can choose $\sigma (c_{0}^{k})=\sigma (c_{0})^{k}$
for $0\leq k<\#C_{\pi }=m$. Moreover, equation (\ref{sigma_c}) implies that $%
\sigma (c_{0})^{m}$ is a scalar multiple of the identity. We can multiply $%
\sigma (c_{0})$ by a scalar to arrange that $\sigma (c_{0})^{m}$ is exactly
the identity. Then, with this definition of $\sigma $ we get an extension $%
\widetilde{\pi }$ of $\pi $ to $C_{\pi }\ltimes N$, namely the
representation $\widetilde{\pi }=\sigma \ltimes \pi $ on the space of $\pi $
given by 
\begin{equation}
\widetilde{\pi }(c,n):=\sigma (c)\circ \pi (n),\text{ \ \ }c\in C_{\pi },%
\text{ }n\in N.  \label{pi-tilda}
\end{equation}

We can get other extensions by twisting this with a character of $C_{\pi }$.
\smallskip

Clifford-Mackey's theory \cite{Clifford37, Mackey49} then says,

\begin{theorem}
\label{T-MLGM}We have,

\begin{enumerate}
\item \label{Ext}\medskip The irreps of the form (\ref{pi-tilda}) are (up to
twist by a character of $C_{\pi }$) all the possible extensions of $\pi $
from $N$ to $C_{\pi }\ltimes N$.

\item \label{Class}All irreps $\rho $ of $G$ containing $\pi $ are obtained
by inducing one of these extensions from $C_{\pi }\ltimes N$ to $G$.
\end{enumerate}
\end{theorem}

As a result we obtain,

\begin{corollary}
\label{C-Res-N}We have,

\begin{enumerate}
\item \label{Spec}Irreps of $G$ have the same $N$-spectrum iff they differ
by twist by a character of $C$.

\item \label{Mult}The restriction to $N$ of any member of $\widehat{G}$ is
multiplicity free.
\end{enumerate}
\end{corollary}

For a proof of Corollary \ref{C-Res-N} see Appendix \ref{P-C-Res-N} (Part (%
\ref{Mult}) was proved in \cite{Lehrer72}).\smallskip

Finally, let us rewrite Part (\ref{Spec}) of Corollary \ref{C-Res-N} in a
slightly different and more quantitative way. Denote by $\widehat{G}_{\pi }$
the collection of all irreps of $G$ having $\pi \in \widehat{N}$ in their $N$%
-spectrum. Theorem \ref{T-MLGM} implies,

\begin{corollary}
\label{C-Irr-above-pi}The group of characters of $C_{\pi }$ acts naturally
on $\widehat{G}_{\pi }$ and this action is free and transitive.
\end{corollary}

In particular,

\begin{corollary}
\label{C-Irr-Res}The restriction to $N$ of an irrep $\rho $ of $G$ stays
irreducible iff $\rho $ is not fixed by a twist of any non-trivial character
of $C$.
\end{corollary}

\section{\textbf{Harish-Chandra's "Philosophy of Cusp Forms"\label{A-PCF}}}

In this section we recall several facts from Harish-Chandra's \textit{%
"philosophy of cusp forms"} \ (P-of-CF) \cite{Harish-Chandra70} for the
description/classification of the set of irreps of $GL_{n}$. We follow
closely the exposition of ideas given in \cite{Howe-Moy86} (where the reader
can find more details, including proofs of the various statements). Other
good sources are \cite{Bump04} and the comprehensive study done in \cite%
{Zelevinsky81}.

The upshot of the P-of-CF is a process that exhausts $\widehat{GL}_{n}$ in
three steps:\smallskip

\textbf{Step 1. }Determining the "cuspidal" irreps of the groups $GL_{m},$ $%
m\leq n$.\smallskip

\textbf{Step 2. }Dividing $\widehat{GL}_{n}$ into subsets parametrized by
"cuspidal data".\smallskip

\textbf{Step 3. }Parametrizing the irreps associated with each cuspidal
datum.\smallskip

We will give now more details on Step 2 (see Section \ref{S-CD}) and Step 3
(see Section \ref{S-ICD}), leaving the classification of the irreps to be
given in term of the building blocks - the cuspidal representations of Step
1, which we will not discuss explicitly in this note (see \cite{Bump04,
Gel'fand70, Howe-Moy86, Zelevinsky81}).

\subsection{\textbf{Cuspidal Data Attached to Parabolic Subgroups\label{S-CD}%
}}

Let us denote $V_{n}=\mathbb{F}_{q}^{n}$, and for each $m\leq n$ denote by $%
V_{m}\subset V_{n}$ subspace of $V_{m}$ of vectors having their last $n-m$
coordinates equal to zero, and by $V_{n-m}^{o}$ the complementary subspace
consisting of vectors having their first $m$ coordinates set to zero.

Recall that the \textit{standard flag} associated with an increasing
subsequence of integers 
\begin{equation}
A=\{0=a_{0}<a_{1}<\ldots <a_{l}=n\},  \label{A}
\end{equation}%
is the nested sequence of spaces of $V_{n},$%
\begin{equation}
0=V_{a_{0}}\subset V_{a_{1}}\subset ...\subset V_{a_{l}}=V_{n}.  \label{FA}
\end{equation}

To the flag (\ref{FA}), we associate the following triple of groups:%
\begin{eqnarray}
\text{(1) \ \ }P_{A} &=&\{g\in GL_{n};\text{ }g(V_{a_{i}})=V_{a_{i}}\text{
for every }i\},  \label{LU} \\
\text{(2) \ \ }U_{A} &=&\{g\in P_{A};\text{ }(g-1)(V_{a_{i}})\subset
V_{a_{i-1}}\text{ for every }i\},  \notag \\
\text{(3) \ \ }L_{A} &=&\{g\in P_{A};\text{ }%
g(V_{n-a_{i}}^{o})=V_{n-a_{i}}^{o}\text{ for every }i\}\text{. \ \ \ \ \ } 
\notag
\end{eqnarray}%
The group $P_{A}$ is called the \textit{standard parabolic }subgroup
associated with the flag (\ref{FA}), and the groups $U_{A},L_{A},$ are,
respectively, the \textit{unipotent radical }and \textit{Levi component} of $%
P_{A}$. We have, 
\begin{equation}
\left\{ 
\begin{array}{c}
P_{A}=L_{A}U_{A}; \\ 
L_{A}\simeq \tprod\limits_{i}GL_{\lambda _{i}},%
\end{array}%
\right.  \label{PALA}
\end{equation}%
where $\lambda _{i}=a_{i}-a_{i-1},$ form a partition of $n.$

Now, we can illustrate a recursive process leading to the P-of-CF.

Take $\rho \in \widehat{GL}_{n}$ and consider a standard parabolic subgroup $%
P_{A}\subset GL_{n}.$ Suppose $\rho $ contains a vector invariant under $%
U_{A},$ the unipotent radical of $P_{A}$. Then Frobenius reciprocity \cite%
{Serre77} implies that there is a representation $\kappa $ of $P_{A}$,
trivial on $U_{A},$ such that $\rho $ is contained in the induced
representation of $\kappa $ from $P_{A}$ to $GL_{n},$ 
\begin{equation*}
\rho <Ind_{P_{A}}^{GL_{n}}(\kappa ).
\end{equation*}%
Since a representation of $P_{A},$ trivial on $U_{A},$ is a representation
of $L_{A}=P_{A}/U_{A},$ and $L_{A}$ is a product of $GL_{m}$'s for $m<n,$
the problem of determining the possibilities for $\kappa $ (i.e.,
determining $\widehat{L}_{A}$) is presumably easier than that of determining 
$\widehat{GL}_{n}$. Thus, the problem of determining all $\rho \in \widehat{%
GL}_{n}$ with $U_{A}$ invariant fixed vectors is reduced to the problem of
determining $\widehat{L}_{A}$ and decomposing induced representations.

We described above an inductive procedure for determining $\widehat{GL}_{n}$%
, the building blocks of which are those representations for which no such
reduction is possible, i.e., those irreps $\kappa $ of $GL_{m}$, $m\leq n,$
which contain no $U_{A}$-invariant vectors for any $A\neq \{0,n\}$.
Harish-Chandra called these irreps \textit{cuspidal,} a term suggested by
the theory of automorphic forms.

To make the P-of-CF\ description of $\widehat{GL}_{n}$ more precise, one
introduces the following key definition:

\begin{definition}
\label{D-CD}A \underline{cuspidal datum} is a pair $(P_{A},\kappa )$ where $%
P_{A}\subset GL_{n}$ is a standard parabolic, and $\kappa $ is a cuspidal
irrep of its Levi subgroup $L_{A}$. Two cuspidal data $(P_{A},\kappa )$ and $%
(P_{A^{\prime }},\kappa ^{\prime })$ are \underline{associate} if there is a 
$g\in GL_{n}$ that conjugates the Levi subgroups $L_{A}\subset P_{A}$ to $%
L_{A^{\prime }}\subset P_{A^{\prime }}$ and the corresponding cuspidal
representations $\kappa $ to $\kappa ^{\prime }.$
\end{definition}

A main result of the P-of-CF is

\begin{theorem}[Harish-Chandra]
\label{T-HC}Suppose $\rho \in \widehat{GL}_{n}.$ Up to association, there
exists a unique cuspidal datum $(P_{A},\kappa )$ with%
\begin{equation}
\rho <Ind_{P_{A}}^{GL_{n}}(\kappa ).  \label{Ind-PG-k}
\end{equation}
\end{theorem}

As a consequence we obtain

\begin{corollary}
\label{C-HC}The induced representations $Ind_{P_{A}}^{GL_{n}}(\kappa )$ and $%
Ind_{P_{A^{\prime }}}^{GL_{n}}(\kappa ^{\prime })$ have components in common
if and only if the cuspidal data $(P_{A},\kappa )$ and $(P_{A^{\prime
}},\kappa ^{\prime })$ are associated. Moreover, in this case the induced
representations are equivalent.
\end{corollary}

\begin{remark}
We would like to elaborate a bit more on the structure of a cuspidal datum.

Note that since $L_{A}$ is a product of $GL_{\lambda _{i}}$ (see (\ref{PALA}%
)) then any irrep $\kappa $ of $L_{A}$ will be a tensor product 
\begin{equation}
\kappa =\tbigotimes\limits_{i}\kappa _{\lambda _{i}},  \label{kapa-TP}
\end{equation}%
of representations $\kappa _{\lambda _{i}}$ of the $GL_{\lambda _{i}}$. In
particular, $\kappa $ will be cuspidal if and only if every factor $\kappa
_{\lambda _{i}}$ in (\ref{kapa-TP}) is cuspidal.

Moreover, since up to association the factors $GL_{\lambda _{i}}$ of $L_{A}$
can be permuted arbitrarily, it will be useful for us to have certain
"standard" organization of \ the cuspidal datum:

First, we will call the $A$ (see (\ref{A})) for which the differences $%
\lambda _{i}=a_{i}-a_{i-1}$ are monotonically (weakly) decreasing the 
\underline{\textit{standard representative}} of its association class.
Moreover, in this case we will call a cuspidal representation of $L_{A}$ of
the form (\ref{kapa-TP}) \underline{\textit{standard}} iff, for any cuspidal
representation $\kappa _{\lambda }$ of $GL_{\lambda }$, the set of indices $%
i $ such that $\lambda _{i}=\lambda $ and $\kappa _{\lambda _{i}}\simeq
\kappa _{\lambda }$ is a consecutive set. When these conditions hold, we
will also say that $(P_{A},\kappa )$ is a \underline{standard cuspidal datum}%
.

Second, \ we define the \underline{\textit{decomposition parabolic}} $P_{%
\widetilde{A}}$ attached to a given standard cuspidal datum with parabolic $%
P_{A}$. This is a parabolic $P_{\widetilde{A}}$ that contains $P_{A}$ and
defined by an increasing sequence $\widetilde{A}=\{\widetilde{a}%
_{i}\}\subset A$, where the $\widetilde{a}_{i}$ are such that, two blocks of 
$L_{A}$ belong to the same block of $L_{\widetilde{A}}$ if and only if the
cuspidal representations attached to the two blocks by the given cuspidal
datum are isomorphic.
\end{remark}

In conclusion, Theorem \ref{T-HC} implies that the process of forming
induced representations from parabolic subgroups using cuspidal
representations of Levi subgroups, partitions $\widehat{GL}_{n}$ into
disjoint subsets parametrized by association classes of standard cuspidal
data.

\subsection{\textbf{Parametrizing the Irreps Associated with a Cuspidal
datum \label{S-ICD}}}

The next step in the philosophy of cups forms is to parametrize the
irreducible components of the induced representations $Ind_{P_{A}}^{GL_{n}}(%
\kappa )$.

We first observe that, for a standard cuspidal datum $(P_{A},\kappa )$ with
decomposition parabolic $P_{\widetilde{A}}$, \ we have

\begin{proposition}
\label{P-Irr}If $\sigma $ is any irreducible component of $Ind_{P_{A}}^{P_{%
\widetilde{A}}}(\kappa )$ then $Ind_{P_{\widetilde{A}}}^{GL_{n}}(\sigma )$
is irreducible.
\end{proposition}

Proposition \ref{P-Irr} follows from Theorem \ref{T-HC} and Mackey
irreducibility criteria \cite{Mackey51, Serre77}.

\subsubsection{\textbf{Parametrizing the Irreducible Components of }$%
Ind_{P_{A}}^{P_{\protect\widetilde{A}}}(\protect\kappa )$}

\begin{figure}[h]\centering
\includegraphics
{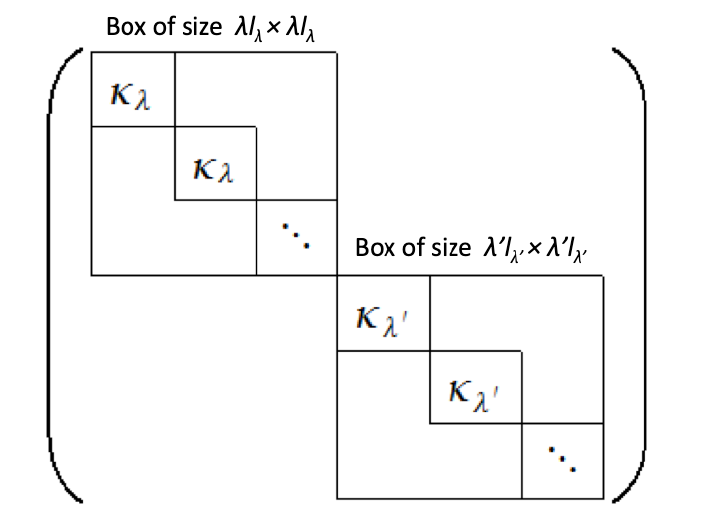}%
\caption{Cuspidal datum consisting of two blocks of an isobaric cuspidal
data.}\label{cd}%
\end{figure}

Thus, by Proposition \ref{P-Irr}, all the reducibility of $%
Ind_{P_{A}}^{GL_{n}}(\kappa )$ happens in the blocks of the Levi component $%
L_{\widetilde{A}}$ and analysis of reducibility reduces to analysis of
cuspidal datum for which all the blocks have equivalent
representations---see Figure \ref{cd} for illustration. A datum of this kind
will be called an \underline{isobaric} cuspidal datum.%

Let us elaborate a bit more on the reduction to an isobaric case.

Indeed, the standard cuspidal datum $(P_{A},\kappa )$ is given by a
representation 
\begin{equation*}
\kappa =\tbigotimes\limits_{\lambda }\kappa _{\lambda }^{\otimes
^{l_{\lambda }}}\text{ \ of the Levi component \ }L_{A}=\tprod\limits_{%
\lambda }(GL_{\lambda })^{l_{\lambda }},
\end{equation*}%
and the decomposition parabolic $P_{\widetilde{A}}\supset P_{A}$ has Levi
component%
\begin{equation*}
L_{\widetilde{A}}=\tprod\limits_{\lambda }GL_{\mu _{\lambda }},\text{ \ }\mu
_{\lambda }=\lambda l_{\lambda }.
\end{equation*}

Moreover, since $\kappa $ is trivial on $U_{A}\supset U_{\widetilde{A}},$
the representation $Ind_{P_{A}}^{P_{\widetilde{A}}}(\kappa )$ will be
effectively the representation $Ind_{L_{\widetilde{A}}\cap P_{A}}^{L_{%
\widetilde{A}}}(\kappa )$ of $L_{\widetilde{A}}$ induced from $\kappa $
considered as a representation of the parabolic subgroup%
\begin{equation*}
L_{\widetilde{A}}\cap P_{A}=\tprod\limits_{\lambda }\left( GL_{_{\mu
_{\lambda }}}\cap P_{A}\right) .
\end{equation*}%
In particular, we have 
\begin{equation*}
Ind_{L_{\widetilde{A}}\cap P_{A}}^{L\widetilde{_{A}}}(\kappa )\simeq
\tbigotimes\limits_{\lambda }Ind_{GL_{\mu _{\lambda }}\cap P_{A}}^{GL_{\mu
_{\lambda }}}(\kappa _{\lambda }^{\otimes ^{l_{\lambda }}}).
\end{equation*}%
The Levi component of $GL_{\mu _{\lambda }}\cap P_{A}$ is a product of $%
l_{\lambda }$ copies of $GL_{\lambda }$.

The conclusion is that, indeed, in order to parametrize the irreducible
components of $Ind_{P_{A}}^{P_{\widetilde{A}}}(\kappa )$ it suffices to
analyze the case of parabolic induction attached to an isobaric cuspidal
datum.

\subsubsection{\textbf{Parametrizing the Irreps Attached to Isobaric
Cuspidal Datum\label{A-UCD}}}

We take $n=l\lambda $ and consider a parabolic $P_{l\times \lambda }\subset
GL_{n}$ with Levi component $L_{l\times \lambda }=(GL_{\lambda })^{l}$
equipped with a representation of the form $\kappa _{\lambda }^{\otimes
^{l}} $, where $\kappa _{\lambda }$ is an irreducible cuspidal
representation of $GL_{\lambda }$.

We would like to parametrize the irreducible components of $Ind_{P_{l\times
\lambda }}^{GL_{n}}(\kappa _{\lambda }^{\otimes ^{l}})$.

Consider the intertwining algebra $End_{GL_{n}}(Ind_{P_{l\times \lambda
}}^{GL_{n}}(\kappa _{\lambda }^{\otimes ^{l}}))$. From general theory (e.g.,
from Burnside's double commutant theorem \cite{Weyl46}) follows that

\begin{proposition}
The joint action of $\ GL_{n}$ and $End_{GL_{n}}(Ind_{P_{l\times \lambda
}}^{GL_{n}}(\kappa _{\lambda }^{\otimes ^{l}}))$ on $Ind_{P_{l\times \lambda
}}^{GL_{n}}(\kappa _{\lambda }^{\otimes ^{l}})$ induces a canonical
bijection between the irreps of $End_{GL_{n}}(Ind_{P_{l\times \lambda
}}^{GL_{n}}(\kappa _{\lambda }^{\otimes ^{l}}))$ and the irreducible
components of $Ind_{P_{l\times \lambda }}^{GL_{n}}(\kappa _{\lambda
}^{\otimes ^{l}}).$
\end{proposition}

\begin{example}
The standard parabolic attached to the set $\{0,1,...,n\}$ is the Borel
subgroup $B\subset GL_{n}$ of upper triangular matrices. The irreducible
representations that appear in $Ind_{B}^{GL_{n}}(\mathbf{1})$ are called 
\underline{spherical principal series} (SPS) representations. Consider the
algebra (under convolution) $H(GL_{n}//B)$ of functions on $GL_{n}$ which
are bi-invariant with respect to $B.$ This algebra is called the \underline{%
spherical Hecke algebra}. Realizing $Ind_{B}^{GL_{n}}(\mathbf{1})$ on the
space of functions on $G/B$ we obtain an identification 
\begin{equation}
End_{GL_{n}}(Ind_{B}^{GL_{n}}(\mathbf{1}))=\mathcal{H}(GL_{n}//B).
\label{SHA}
\end{equation}
\end{example}

Identity (\ref{SHA}) has an important generalization as follows. Consider
the spherical Hecke algebra $\mathcal{H}(GL_{l}(\mathbb{F}_{q^{\lambda
}})//B(\mathbb{F}_{q^{\lambda }}))$. It turns out that, $%
End_{GL_{n}}(Ind_{P_{l\times \lambda }}^{GL_{n}}(\kappa _{\lambda }^{\otimes
^{l}}))$ has a presentation in term of generators and relations that is
identical to the presentation using standard generators and relations for $%
\mathcal{H}(GL_{l}(\mathbb{F}_{q^{\lambda }})//B(\mathbb{F}_{q^{\lambda }}))$
\cite{Howlett-Lehrer80, Howe-Moy86, Lusztig84}. In particular,

\begin{theorem}
There is an explicit isomorphism%
\begin{equation*}
End_{GL_{n}}(Ind_{P_{l\times \lambda }}^{GL_{n}}(\kappa _{\lambda }^{\otimes
^{l}}))\simeq \mathcal{H}(GL_{l}(\mathbb{F}_{q^{\lambda }})//B(\mathbb{F}%
_{q^{\lambda }})),
\end{equation*}%
that preserves the natural $L^{2}$-structures on the two algebras up to
multiples.
\end{theorem}

As a corollary we get

\begin{corollary}
\label{C-param-isobaric-rep}There is a canonical bijection 
\begin{equation}
\beta :\widehat{GL}_{n}(Ind_{P_{l\times \lambda }}^{GL_{n}}(\kappa _{\lambda
}^{\otimes ^{l}}))\longleftrightarrow \widehat{GL}_{l}(Ind_{B(\mathbb{F}%
_{q^{\lambda }})}^{GL_{l}(\mathbb{F}_{q^{\lambda }})}(\mathbf{1})),
\label{C-beta}
\end{equation}%
with the following properties:

\begin{enumerate}
\item The multiplicity of $\rho $ in $Ind_{P_{l\times \lambda
}}^{GL_{n}}(\kappa _{\lambda }^{\otimes ^{l}})$ is equal the multiplicity of 
$\beta (\rho )$ in $Ind_{B(\mathbb{F}_{q^{\lambda }})}^{GL_{l}(\mathbb{F}%
_{q^{\lambda }})}(\mathbf{1}).$

\item \label{C-param-isobaric-rep-P2}For every $\ \rho \in \widehat{GL}%
_{n}(Ind_{P_{l\times \lambda }}^{GL_{n}}(\kappa _{\lambda }^{\otimes ^{l}}))$
we have%
\begin{equation}
\frac{\dim (\rho )}{\dim (\beta (\rho ))}=\dim (\kappa _{\lambda }^{\otimes
^{l}})\frac{\#((GL_{n}(\mathbb{F}_{q})/P_{l\times \lambda }(\mathbb{F}_{q}))%
}{\#(GL_{l}(\mathbb{F}_{q^{\lambda }})/B(\mathbb{F}_{q^{\lambda }}))}.
\label{dim-isobaric-rep}
\end{equation}
\end{enumerate}
\end{corollary}

According to (\ref{C-beta}), in order to parametrize irreps attached to
isobaric cuspidal datum, it is enough to decompose\textbf{\ }the space $%
Ind_{B}^{GL_{l}}(\mathbf{1}).$

\subsubsection{\textbf{Parametrizing the Spherical Principal Series
Representations}\label{S-Spherical}}

We want to parametrize the SPS representations, i.e., the irreps that appear
in $Ind_{B}^{GL_{l}}(\mathbf{1})$.

Let us denote by $W$ the standard Weyl group (i.e., the permutation
matrices) in $GL_{l}$. In addition, for a standard parabolic $P_{A}\subset
GL_{l},$ we consider the induced representation $I_{A}=Ind_{P_{A}}^{GL_{l}}(%
\mathbf{1})$, the subgroup $W_{A}=W\cap P_{A}<$ $W$, and the induced
representation $Y_{A}=Ind_{W_{A}}^{W}(\mathbf{1}))$.

The following theorem gives an effective parametrization of the SPS\
representations:

\begin{theorem}
\label{T-SHA-Sn}There is a unique bijection 
\begin{equation}
\alpha :\widehat{GL}_{l}(Ind_{B}^{GL_{l}}(\mathbf{1}))\longleftrightarrow 
\widehat{W},  \label{C-alpha}
\end{equation}%
such that for every standard parabolic subgroup $P_{A}$ we have $\rho \in 
\widehat{GL}_{l}(I_{A})$ if and only if $\alpha (\rho )\in \widehat{W}%
(Y_{A}).$ Moreover, in that case we have 
\begin{equation}
\dim (Hom_{GL_{l}}(\rho ,I_{A}))=\dim (Hom_{W}(\alpha (\rho ),Y_{A})).
\label{Mult-alpha}
\end{equation}
\end{theorem}

The standard justification for Theorem \ref{T-SHA-Sn} that we are aware of
(see \cite{Bourbaki68, Howe-Moy86, Iwahori66}), goes by the name "Tits's
deformation argument". However, due to its fundamental rule in the
representation theory of the finite general linear groups, it might be
worthwhile to give other derivations of Theorem \ref{T-SHA-Sn}. In
particular, in Appendix \ref{A-IrrS_l-SPS} we sketch a modified approach to
the proof of Theorem \ref{T-SHA-Sn}, which seems to be more elementary than
the approach currently used in the literature, and might give additional
valuable information on representations of $GL_{l}.\medskip $

In conclusion, we obtain a classification of the irreducible components that
appear in $Ind_{B}^{GL_{l}}(\mathbf{1})$. The parametrization is given in
term of partitions of $l$ as is the case for the irreps of $W=S_{l}$ \cite%
{Ceccherini-Silberstein-Scarabotti-Tolli10, Fulton-Harris91, Weyl46}. As a
consequence, using (\ref{C-beta}) and (\ref{C-alpha}) we get a
parametrization of the irreps that appear inside $Ind_{P_{l\times \lambda
}}^{GL_{n}}(\kappa _{\lambda }^{\otimes ^{l}})$ in terms of partitions of $%
l. $

\subsection{\textbf{Summary}}

We have learned that, according to the philosophy of cusp forms, an irrep of 
$GL_{n}$ is specified by its cuspidal datum, augmented by a collection of
partitions. Precisely, for each cuspidal representation $\kappa _{\lambda }$
that appears in the datum, if $l$ is the number of times that $\kappa
_{\lambda }$ appears, then we augment the datum with a partition of $l$
attached to $\kappa _{\lambda }$.

\section{\textbf{Representations of }$S_{l}$\textbf{\ and the Spherical
Principal Series for }$GL_{l}$\textbf{\label{A-IrrS_l-SPS}}}

We sketch (for a more comprehensive treatment, including proofs of the main
statements, see \cite{Gurevich-Howe19}) a seemingly not so well known
organization of the representation theories of - on the one hand the
symmetric group $S_{l},$ and on the other hand the spherical principal
series (SPS) representations of $GL_{l}$.

The modified perspective, starts by putting at the forefront two naturally
arising structures - the \textit{symmetric }and \textit{spherical }monoids%
\textit{. }Then, as a logical outcome of their intrinsic qualities, one is
able to reproduce, in an elegant way, first the classifications - in terms
of partitions - of the irreps of the symmetric group, and of the SPS
representations of the finite general linear group; and second to recast the
bijection stated in Theorem \ref{T-SHA-Sn}, in terms of an isomorphism - the
only one possible - between these two aforementioned monoids.

As a by-product we get additional valuable information on the SPS
representations of $GL_{l}$ from those of $S_{l}$. For example, this is how
we obtained the Pieri rule for $GL_{l}$ in Section \ref{S-C-of-I}.

\subsection{\textbf{The Symmetric Monoid and the Classification of the
Irreps of }$S_{l}$\label{A-SyM}}

Consider the set $M(S_{l})$ of representations of the symmetric group $S_{l}$
up to equivalence. The direct sum operation $\oplus $ on representations
induces, in a natural way, a structure of a monoid on $M(S_{l})$ with
identity element given by the $0$ representation. We will call it the 
\underline{symmetric monoid}. It is well known \cite{Serre77} that the
symmetric monoid (and the analogous structure for any finite group) is a
free abelian semigroup on the irreducible representations.

The symmetric monoid is equipped naturally with

\begin{itemize}
\item an "inner product", given by the non-degenerate symmetric bilinear
form, 
\begin{equation*}
\left\langle \sigma ,\sigma ^{\prime }\right\rangle =\dim (Hom(\sigma
,\sigma ^{\prime })),\text{ \ \ }\sigma ,\sigma ^{\prime }\in M(S_{l})\text{;%
}
\end{equation*}%
and,

\item a partial order, given by 
\begin{equation*}
\sigma <\sigma ^{\prime }\text{ \ iff \ }\sigma \text{ is a
sub-representation of }\sigma ^{\prime }.
\end{equation*}
\end{itemize}

The monoid $M(S_{l})$ has an easily defined and much-studied collection of
elements, called \textit{Young modules}, parametrized by partitions of $l$ 
\cite{Ceccherini-Silberstein-Scarabotti-Tolli10}. Indeed, for a partition $%
L=\{l_{1}\geq l_{2}\geq ...\geq l_{s}\}$ of $l,$ we define the \textit{Young
module }associated with $L$, to be the induced module%
\begin{equation}
Y_{L}=Ind_{S_{l_{1}}\times ...\times S_{l_{s}}}^{S_{l}}(\mathbf{1}),
\label{YL}
\end{equation}%
where the subgroup $S_{l_{1}}\times ...\times S_{l_{s}}$ is contained in $%
S_{l}$ in the standard way.

\ We would like to point out two properties of the collection (\ref{YL}) of
Young modules. Both involve the \textit{dominance relation} on the set $%
\mathcal{P}_{l}$ of partitions of $l$. Recall that,

\begin{definition}
\label{D-Dom} If, in addition to $L$ as above, we have another partition $%
L^{\prime }=\{l_{1}^{\prime }\geq l_{2}^{\prime }\geq ...\geq l_{r}^{\prime
}\}$ of $l,$ then we say that $L^{\prime }$ \underline{dominates} $L$, and
write $L\prec L^{\prime },$ if $r\leq s$ and%
\begin{equation*}
\sum_{i=1}^{j}l_{i}\leq \sum_{i=1}^{j}l_{i}^{\prime },\text{ \ \ \ for }%
j=1,...,r.
\end{equation*}
\end{definition}

With this terminology we have,

\begin{proposition}
\label{P-Mon-Max-YL}Suppose $L$ is a partition of $l$. Then,

\begin{enumerate}
\item \label{Mon-YL}For any partition $L^{\prime }$ of $l$, we have $%
Y_{L^{\prime }}\lvertneqq Y_{L}$ if and only if $L\precneqq L^{\prime }$.

\item \label{Max-YL}There is a unique irrep%
\begin{equation}
\sigma _{L}<Y_{L},  \label{sigmaL}
\end{equation}%
which is not contained in $Y_{L^{\prime }},$ for any partition $L^{\prime }$
that strictly dominates $L$.\textit{\ } Moreover, the multiplicity of $%
\sigma _{L}$ in $Y_{L}$ is one.
\end{enumerate}
\end{proposition}

An elementary proof of Proposition \ref{P-Mon-Max-YL} can be found in \cite%
{Gurevich-Howe19, Howe-Moy86}.\smallskip

As a corollary of Proposition \ref{P-Mon-Max-YL} we reproduce the well known
classification of irreps of the symmetric group:

\begin{corollary}[\textbf{Classification}]
\label{C-Class-Sl}The irreps $\sigma _{L},$ $L\in \mathcal{P}_{l},$ are
pairwise non-isomorphic, and exhaust $\widehat{S}_{l}$. In particular, the
irreps of $S_{l}$ are naturally parametrized\ by partitions of $l$.
\end{corollary}

\subsection{\textbf{The Spherical Monoid and the Classification of the
Constituents of }$Ind_{B}^{GL_{l}}(\mathbf{1})$\label{A-SM}}

Consider the set $M_{B}(GL_{l})$ of representations (up to equivalence) of
the group $GL_{l}$ for which all subreps have a $B$-invariant vector. 
\textit{Mutatis mutandem, }as in the case of the symmetric monoid, the set $%
M_{B}(GL_{l})$ with the direct sum operation $\oplus $ is a monoid. It is
the free abelian semigroup on the irreducible representations of $GL_{l}$
that appear in $Ind_{B}^{GL_{l}}(\mathbf{1})$, that is, in the permutation
action of $GL_{l}$ on the variety of complete flags in $\mathbb{F}_{q}^{l}$.
It inherits a partial order structure $<$ and inner product $\left\langle
\cdot ,\cdot \right\rangle $ from the monoid of all representations of $%
GL_{l}$. We will call it the \underline{spherical monoid}.

As in the case of the symmetric monoid, also the spherical monoid has a
easily defined and much-studied collection of elements parametrized by
partitions.\ In this case, starting with a partition $L=\{l_{1}\geq
l_{2}\geq ...\geq l_{s}\}$ of $l,$ we consider the (parabolically) induced
module 
\begin{equation}
I_{L}=I_{A(L)}=Ind_{P_{A(L)}}^{GL_{l}}(\mathbf{1})<\text{\textbf{\ }}%
Ind_{B}^{GL_{l}}(\mathbf{1}),  \label{IL}
\end{equation}%
where $P_{A(L)}$ is the standard parabolic (see Section \ref{S-CD})
associated with the set $A(L)=\{0,l_{1},...,l_{1}+...+l_{s}=l\}.$\ 

The collection of induced modules $I_{L}$ (\ref{IL}) satisfies the
properties:

\begin{proposition}
\label{P-Mon-Max-IL}Suppose $L$ is a partition of $l$. Then,

\begin{enumerate}
\item \label{Mon-IL}For any partition $L^{\prime }$ of $l$, we have $%
I_{L^{\prime }}\lvertneqq I_{L}$ if and only if $L\precneqq L^{\prime }$.

\item \label{Max-IL}There is a unique irrep 
\begin{equation}
\rho _{L}<I_{L},  \label{rhoL}
\end{equation}%
which is not contained in $I_{L^{\prime }},$ for any partition $L^{\prime }$
that strictly dominates $L$.\textit{\ } Moreover, the multiplicity of $\rho
_{L}$ in $I_{L}$ is one.
\end{enumerate}
\end{proposition}

An elementary proof of Proposition \ref{P-Mon-Max-IL} can be found in \cite%
{Gurevich-Howe19, Howe-Moy86}.\smallskip

As a corollary of Proposition \ref{P-Mon-Max-IL} we reproduce the well known
classification of the spherical principal series representations of $GL_{l}$:

\begin{corollary}[\textbf{Classification}]
\label{C-Class-SPS}The irreps $\rho _{L},$ $L\in \mathcal{P}_{l},$ are
pairwise non-isomorphic, and exhaust $\widehat{GL}_{l}(Ind_{B}^{GL_{l}}(%
\mathbf{1}))$. In particular, the SPS representations of $GL_{l}$ are
naturally parametrized\ by partitions of $l$.
\end{corollary}

\subsection{\textbf{Correspondence between Irreps of }$S_{l}$\textbf{\ and
the Spherical Principal Series of }$GL_{l}$\label{A-Id-IrrSl-SPSGl}}

Consider the symmetric and spherical monoids, $M(S_{l})$ and $M_{B}(GL_{l})$%
, respectively. We have,

\begin{theorem}
\label{T-Id}The assignment%
\begin{equation*}
I_{L}\longmapsto Y_{L},\text{ \ \ }L\in \mathcal{P}_{l},
\end{equation*}%
extends uniquely to an isomorphism of monoids 
\begin{equation}
\alpha :M_{B}(GL_{l})\widetilde{\longrightarrow }M(S_{l}),  \label{iso-alpha}
\end{equation}%
that satisfies the following (equivalent) conditions:

\begin{description}
\item[C1] $\alpha $ preserves the partial orders $<$ on both monoids.

\item[C2] $\alpha $ preserves the inner products $\left\langle
,\right\rangle $ on both monoids.
\end{description}

Moreover, the aforementioned extension $\alpha $ satisfies $\alpha (\rho
_{L})=\sigma _{L}$, for every $L\in \mathcal{P}_{l}$ (see (\ref{rhoL}) and (%
\ref{sigmaL})).
\end{theorem}

The uniqueness part of Theorem \ref{T-Id} is immediate, while the existence
part is a direct consequence of the Bruhat decomposition \cite{Borel69,
Bruhat56} (for more details see \cite{Gurevich-Howe19, Howe-Moy86}%
).\smallskip

Finally, we note that Theorem \ref{T-SHA-Sn} follows from Theorem \ref{T-Id}.

\subsection{\textbf{Estimating the Dimensions of the Spherical Principal
Series Representations of }$GL_{l}$\label{A-est-dim-SPS}}

Proposition \ref{P-Mon-Max-IL} has the following consequences for the
dimensions of the SPS representations:

\begin{corollary}[\textbf{Dimension}]
\label{C-dim}We have,

\begin{enumerate}
\item \textbf{Formula. }The dimension of the\ SPS representation $\rho _{L}$
(\ref{rhoL}) attached to a partition $L=\{l_{1}\geq l_{2}\geq ...\geq
l_{s}\} $of $l,$ satisfies $\dim (\rho _{L})=\dim (I_{L})+o(...),$ as $%
q\rightarrow \infty $, and, in particular,%
\begin{equation}
\dim (\rho _{L})=q^{d_{L}}+o(...),  \label{dim-rho_L}
\end{equation}%
where $d_{L}=\tsum\limits_{1\leq i<j\leq s}l_{i}l_{j}$.

\item \textbf{Monotonicity}. Suppose $L$ and $L^{\prime }$ are two
partitions of $l,$ with $L$ $\precneqq L^{\prime }$. Then, $d_{L^{\prime
}}<d_{L}.$
\end{enumerate}
\end{corollary}

For a proof of Corollary \ref{C-dim} see Appendix \ref{P-C-dim}.

\section{\textbf{Proofs}}

\subsection{\textbf{Proofs for Section \protect\ref{S-CR-TR}}}

\subsubsection{\textbf{Proof of Proposition \protect\ref{P-TRF}\label%
{P-P-TRF}}}

\begin{proof}
Note that $(L^{2}(\mathbb{F}_{q}^{n}))^{\otimes ^{k}}=L^{2}(M_{k,n})$, where 
$M_{k,n}$ denotes the space of matrices of size $k\times n$ over $\mathbb{F}%
_{q}$. In particular, the space $L^{2}(M_{k,n})$ contains the regular
representation if and only if $k=n,$ and the existence of cuspidal
representations, for example, tells us that the filtration does not
stabilize before that stage. This completes the proof of the Proposition.
\end{proof}

\subsection{\textbf{Proofs for Section \protect\ref{S-AI-SLn}}}

\subsubsection{\textbf{\ Proof of Lemma \protect\ref{L-CR-SLn}\label%
{P-L-CR-SLn}}}

\begin{proof}
Consider the section $s$ for the determinant morphism $GL_{n}\overset{\det }{%
\longrightarrow }\mathbb{F}_{q}^{\ast }$, sending $a\in \mathbb{F}_{q}^{\ast
}$ to the diagonal matrix with diagonal $(a,1,...,1).$

Suppose $\pi \in \widehat{SL}_{n}$ appears in the restriction of $\rho \in 
\widehat{GL}_{n}$ to $SL_{n}.$ Then, by Fact \ref{F-SLn-Mult}\textbf{\ }we
have,%
\begin{equation}
\rho _{|SL_{n}}=\sum_{a\in C/C_{\pi }}\pi _{a},  \label{res}
\end{equation}%
where $C=\mathbb{F}_{q}^{\ast }$, $\pi _{a}\in \widehat{SL}_{n}$ for $a\in
C, $ is given by $\pi _{a}(g)=\pi (s(a)gs(a)^{-1}),$ and $C_{\pi }$ is the
stabilizer of $\pi $ in $C$.

In particular, for an element $g\in SL_{n}$ with centralizer in $GL_{n}$
satisfying our assumption, we have $\chi _{\pi _{a}}(g)=\chi _{\pi }(g)$ for
every $a\in \mathbb{F}_{q}^{\ast }$. It follows that, 
\begin{equation*}
\frac{\chi _{\rho }(g)}{\dim (\rho )}=\frac{\chi _{\pi }(g)}{\dim (\pi )},
\end{equation*}%
as claimed.
\end{proof}

\subsection{\textbf{Proofs for Section \protect\ref{S-BRW}}}

\subsubsection{\textbf{Proof of Proposition \protect\ref{P-TVB}\label%
{P-P-TVB}}}

\begin{proof}
We remark that the tensor rank $k=1$ irreps that appear in example \ref%
{Ex-rank-k=1} stay irreps after restriction to $SL_{n}$, $n\geq 3$ (using
the argument given in Appendix \ref{P-P-Irr-Res-Gen}).

Next, we compute,%
\begin{eqnarray*}
4\left\Vert P_{C}^{\ast l}-U\right\Vert _{TV}^{2} &\leq &\underset{\mathbf{1}%
\neq \rho \in \widehat{SL}_{n}}{\sum }\dim (\rho )^{2}\left\vert \frac{\chi
_{\rho }(T)}{\dim (\rho )}\right\vert ^{2l} \\
&=&\overset{n}{\underset{k=1}{\sum }\text{ }}\underset{\rho \in (\widehat{SL}%
_{n})_{\otimes ,k}}{\sum }\dim (\rho )^{2}\left\vert \frac{\chi _{\rho }(T)}{%
\dim (\rho )}\right\vert ^{2l} \\
&\leq &(q-1)\left( \frac{q^{n}-1}{q-1}\right) ^{2}\left( \frac{q^{n-1}-1}{%
q^{n}-1}\right) ^{2l}+o(...) \\
&\leq &\frac{1}{q}\left( \frac{1}{q^{2}}\right) ^{l-n}+o(...),
\end{eqnarray*}%
where the first inequality is (\ref{R_el}); the second inequality
incorporates Example \ref{Ex-rank-k=1} for tensor rank $k=1$, and Formulas (%
\ref{CRs-SLn}), (\ref{Dim-SLn}), and (\ref{Card-SLn}) for higher ranks. This
completes the verification of the proposition.
\end{proof}

\subsection{\textbf{Proofs for Section \protect\ref{S-eta-PCF}}}

\subsubsection{\textbf{Proof of the Necessity Statement in Part (1) of
Theorem \protect\ref{T-EC}\label{P-NP1-T-EC}}}

\begin{proof}
We make use of characterization of strict tensor rank given by Proposition %
\ref{P-ID-TR}.

If $\tau $ has strict tensor rank $2k-n-a$, for some $k\geq a>0,$ then $\tau 
$ has a fixed vector for $H_{2k-n-a}^{^{\prime }}$, the stabilizer of the
first $2k-n-a$ coordinates subspace in $\mathbb{F}_{q}^{k}$. So, $\tau $ is
contained in $Ind_{H_{2k-n-a}}^{GL_{k}}(\mathbf{1})$.

Consider the parabolic $P_{2k-n-a,\text{ }2(n-k)+a}\subset GL_{n}$, its
unipotent radical $U_{2k-n-a,\text{ }2(n-k)+a},$ and Levi component $%
L_{2k-n-a,\text{ }2(n-k)+a}\simeq GL_{2k-n-a}\times GL_{2(n-k)+a}.$ In
particular, inside this parabolic we have the group $\mathcal{G}%
_{2k-n-a,2(n-k)+a}=U_{2k-n-a,\text{ }2(n-k)+a}\cdot GL_{2(n-k)+a}$.

Next, consider the parabolic $P_{n-k+a,n-k}\subset GL_{2(n-k)+a}$, with Levi 
$L_{n-k+a,n-k}=GL_{n-k+a}\times $ $GL_{n-k}$ and unipotent radical $%
U_{n-k+a,n-k}\simeq M_{n-k+a,n-k},$ namely, 
\begin{equation*}
P_{n-k+a,n-k}=\{%
\begin{pmatrix}
A & C \\ 
0 & B%
\end{pmatrix}%
\in GL_{2(n-k)+a};\text{ \ }A\in GL_{n-k+a},\text{ }B\in GL_{n-k},\text{ }%
C\in M_{n-k+a,n-k}\}.
\end{equation*}

Denote by $\widetilde{\rho }_{n-k+a,n-k}$ the pullback to $\mathcal{G}%
_{2k-n-a,2(n-k)+a}$ of 
\begin{equation*}
\rho _{n-k+a,n-k}=Ind_{P_{n-k+a,n-k}}^{GL_{2(n-k)+a}}(\mathbf{1}),
\end{equation*}%
and observe that, 
\begin{equation*}
Ind_{P_{k,n-k}}^{GL_{n}}(\tau \otimes \mathbf{1}_{n-k})<Ind_{\mathcal{G}%
_{2k-n-a,2(n-k)+a}}^{GL_{n}}(\widetilde{\rho }_{n-k+a,n-k}).
\end{equation*}%
Note that the parabolics $P_{n-k+a,n-k}$ and $P_{n-k,n-k+a}$ are associate.
This implies, using Theorem \ref{T-HC}, that the induced representations $%
\rho _{n-k+a,n-k}$ and (the similarly defined) $\rho _{n-k,n-k+a}$ are
equivalent. In particular, the pullback $\widetilde{\rho }_{n-k,n-k+a}$ of $%
\rho _{n-k,n-k+a}$ to $\mathcal{G}_{2k-n-a,2(n-k)+a}$ satisfies 
\begin{equation*}
Ind_{P_{k,n-k}}^{GL_{n}}(\tau \otimes \mathbf{1}_{n-k})<Ind_{\mathcal{G}%
_{2k-n-a,2(n-k)+a}}^{GL_{n}}(\widetilde{\rho }_{n-k,n-k+a}).
\end{equation*}%
But the group $H_{k-a}\subset GL_{n},$ that fixes the first $k-a$
coordinates subspace in $\mathbb{F}_{q}^{n},$ acts trivially on $Ind_{%
\mathcal{G}_{2k-n-a,2(n-k)+a}}^{GL_{n}}(\widetilde{\rho }_{n-k,n-k+a}),$ so,
using Proposition \ref{P-ID-TR}, we conclude that $Ind_{P_{k,n-k}}^{GL_{n}}(%
\tau \otimes \mathbf{1}_{n-k})$ contains irreps of strict tensor rank at
most $k-a$.

This completes the proof of the necessity statement.
\end{proof}

\subsubsection{\textbf{Proof of Claim \protect\ref{C-OO}\label{P-C-OO}}}

\begin{proof}
Consider the natural $GL_{k}\times GL_{n}$-action on the set of matrices $%
M_{k,n}$. An orbit for this action is described by the rank of its elements.
The rank $r$ can vary from $0$ to $k$ (we assume $k\leq n$), and we denote
the corresponding orbit by $(M_{k,n})_{r}$. In particular, we have a
decomposition into direct sum of $GL_{k}\times GL_{n}$-representations,%
\begin{equation*}
L^{2}(M_{k,n})=\sum_{r=0}^{k}L^{2}((M_{k,n})_{r}).
\end{equation*}%
Note that, if $0\leq r<k$, then as $GL_{n}$-representation $%
L^{2}((M_{k,n})_{r})<L^{2}(M_{r,n})$, because each $GL_{n}$-orbit of
matrices of lower rank is equivalent to the open orbit in the matrices of
that rank. So, we see that, representations supported on matrices of lower
rank are of lower strict tensor rank, and, in particular, the strict tensor
rank $k$ part of $L^{2}(M_{k,n})$ is contained in $L^{2}((M_{k,n})_{k})$.
This proves Part (\ref{Part1-C-OO}) of Claim \ref{C-OO}.\smallskip

Next, we want to compute the isotypic components for the action of $GL_{k}$
on $L^{2}((M_{k,n})_{k})$. Since $GL_{k}$ acts freely on $(M_{k,n})_{k}$,
the space $L^{2}((M_{k,n})_{k})$ contains a copy of the regular
representation of $GL_{k}$. Let us denote \ by $\mathcal{H}_{\tau }$ a space
on which $\tau \in \widehat{GL}_{k}$ is represented, and calculate the
multiplicity space 
\begin{equation}
Hom_{GL_{k}}(L^{2}((M_{k,n})_{k}),\mathcal{H}_{\tau })=?  \label{MS}
\end{equation}%
Considering the matrix $I_{k,n}\in (M_{k,n})_{k}$ whose first $k$ diagonal
elements are $1,$ and all other entries are $0$, we can identify $%
(M_{k,n})_{k}=H_{k}\diagdown GL_{n}$, where $H_{k}=Stab_{GL_{n}}(I_{k,n}).$
\ Now, to an intertwiner $\iota $ in (\ref{MS}) we can associate the
function $f_{\iota }:GL_{n}\rightarrow \mathcal{H}_{\tau },$ given by $%
f_{\iota }(g)=\iota (\delta _{H_{k}g}),$ where $\delta _{H_{k}g}$ is the
delta function at the coset $H_{k}g$. Note that since $GL_{k}\cdot
H_{k}=P_{k,n-k}$, and $GL_{k}$ normalizes $H_{k}$, the assignment $\iota
\mapsto f_{\iota }$, gives a morphism 
\begin{equation}
Hom_{GL_{k}}(L^{2}(H_{k}\diagdown GL_{n}),\mathcal{H}_{\tau
})\longrightarrow \{f:GL_{n}\rightarrow \mathcal{H}_{\tau };\text{ }f(pg)=%
\widetilde{\tau }(p)f(g),\text{ }p\in P_{k,n-k},\text{ }g\in GL_{n}\},
\label{MSD}
\end{equation}%
where $\widetilde{\tau }$ is the composition of $\tau $ with the projection $%
P_{k,n-k}\twoheadrightarrow GL_{k}$.

The right-hand side of (\ref{MSD}) is $Ind_{P_{k,n-k}}^{GL_{n}}(\tau \otimes 
\mathbf{1}_{n-k}),$ and, moreover, the \ mapping (\ref{MSD}) is an
isomorphism. This proves Part (\ref{Part2-C-OO}) of Claim \ref{C-OO}.
\end{proof}

\subsubsection{\textbf{Proof of Claim \protect\ref{C-Mul}\label{P-C-Mul}}}

\begin{proof}
We use the development described in Appendix \ \ref{A-IrrS_l-SPS}. \ 

Consider the isomorphism $\alpha $ (\ref{iso-alpha}) given in Theorem \ref%
{T-Id}, between the spherical monoid $M_{B}(GL_{n})$ and the symmetric
monoid $M(S_{n}).$ It preserves the inner product structures, defined by the
intertwining number pairings, on both sides. But, by the way the SPS
representations of the general linear groups, and the irreps of the
symmetric groups, are assigned to partitions (see Appendices \ref{A-SM} and %
\ref{A-SM}, respectively), we know that $\alpha $ sends $I_{\rho _{D}}$ (\ref%
{I_rho_D}) to $I_{\sigma _{D}}$ (\ref{I-sigma-D}). Concluding, for every
partition $E$ of $n$ and $D$ of $k$, the identity $\left\langle \rho
_{E},I_{\rho _{D}}\right\rangle =\left\langle \sigma _{E},I_{\sigma
_{D}}\right\rangle $ (\ref{Mul-L-D}) holds. This completes the proof of
Claim \ref{C-Mul}.
\end{proof}

\subsubsection{\textbf{Proof of Theorem \protect\ref{T-I-tauD-Dec}\label%
{P-T-I-tauD-Dec}}}

\begin{proof}
Recall \textit{Schur duality} \cite{Schur27}: the groups $GL_{m}(%
\mathbb{C}
)$ and $S_{n}$ both act in an obvious way on the $n$-fold tensor power $(%
\mathbb{C}
^{m})^{\otimes ^{n}}$of $%
\mathbb{C}
^{m}$. The actions of $GL_{m}(%
\mathbb{C}
)$ and $S_{n}$ commute with each other, and moreover, they generate mutual
commutants in the endomorphisms of $(%
\mathbb{C}
^{m})^{\otimes ^{n}}$. The resulting correspondence of representations is
compatible with the parametrizations of the representations of $S_{n}$ and
of $GL_{m}(%
\mathbb{C}
)$ with Young diagrams \cite{Howe92}.

If we look at the action of $S_{k}\times GL_{m}(%
\mathbb{C}
)$ on $(%
\mathbb{C}
^{m})^{\otimes ^{k}}$, then the isotypic subspace for the representation $%
\sigma _{D}$ of $S_{k}$ will be isomorphic to $\sigma _{D}\otimes \rho ^{D}$%
, as a representation of $S_{k}\times GL_{m}(%
\mathbb{C}
)$. Here $\rho ^{D}$ is the representation of $GL_{m}(%
\mathbb{C}
)$ parametrized by the Young diagram $D$. Similarly, if we look at the
action of $S_{n-k}\times GL_{m}(%
\mathbb{C}
)$ on the $S_{n-k}$ invariants, then the action of $S_{n-k}\times GL_{m}(%
\mathbb{C}
)$ is $\mathbf{1}_{n-k}\otimes S^{n-k}(%
\mathbb{C}
^{m})$. Here $S^{a}(%
\mathbb{C}
^{m})$ indicates the $a$-th symmetric power of $%
\mathbb{C}
^{m}$, and is the representation of $GL_{m}(%
\mathbb{C}
)$ corresponding to the diagram with a single row of length $a$.

Taking the tensor product, we conclude that the isotypic component of $(%
\mathbb{C}
^{m})^{\otimes ^{n}}$ for the representation $\sigma _{D}\otimes \mathbf{1}%
_{n-k}$ of $S_{k}\times S_{n-k},$ is the $S_{k}\times S_{n-k}\times GL_{m}(%
\mathbb{C}
)$-module%
\begin{equation*}
(\sigma _{D}\otimes \mathbf{1}_{n-k})\otimes (\rho ^{D}\otimes S^{n-k}(%
\mathbb{C}
^{m})).
\end{equation*}

The Pieri rule for the complex general linear group \cite{Howe92} tells us
that the tensor product $\rho ^{D}\otimes S^{n-k}(%
\mathbb{C}
^{m})$ of $GL_{m}(%
\mathbb{C}
)$-modules decomposes in a multiplicity-free sum of representations $\rho ^{%
\widetilde{D}}$, where $\widetilde{D}$ is as described in the statement of
the theorem: $\widetilde{D}$ has $n$ boxes, contains $D$, and $\widetilde{D}%
-D$ is a skew row.

Now consider the $S_{n}\times $ $GL_{m}(%
\mathbb{C}
)$-module generated by $\rho ^{D}\otimes S^{n-k}(%
\mathbb{C}
^{m})$. From Schur duality, we know that it is the sum $\sigma _{\widetilde{D%
}}\otimes \rho ^{\widetilde{D}}$, where $\widetilde{D}$ runs through the set
of diagrams of the previous paragraph, \textit{mutatis mutandem}, of the
statement of the theorem. On the other hand, this is the $S_{n}\times $ $%
GL_{m}(%
\mathbb{C}
)$-module generated by the $S_{k}\times S_{n-k}$-module $\sigma _{D}\otimes 
\mathbf{1}_{n-k}$. It follows that the restriction of a representation $%
\sigma _{E}$ of $S_{n}$ to $S_{k}\times S_{n-k}$ contains $\sigma
_{D}\otimes \mathbf{1}_{n-k}$ if and only if $E=\widetilde{D},$ as described
above, and then the multiplicity of\ $\sigma _{D}\otimes \mathbf{1}_{n-k}$
in $\sigma _{\widetilde{D}}$ is $1$. The theorem now follows by Frobenius
reciprocity.
\end{proof}

\subsubsection{\textbf{Proof of Corollary \protect\ref{C-TR-PofCF}\label%
{P-C-TR-PofCF}}}

\begin{proof}
It is not difficult to see (e.g., using the intrinsic characterization given
by Proposition \ref{P-ID-TR}) that for SPS representation strict tensor rank
and tensor rank agree. Take a Young diagram $D\in \mathcal{Y}_{n}$ with
longest row of size $d_{1}$, and consider the corresponding SPS
representation $\rho _{D}$. According to the Pieri rule (see Theorem \ref%
{T-I-tauD-Dec}), $k=n-d_{1}$ is the first such that $\rho _{D}$ appears
inside an induced representation of the form $I_{\rho
_{E}}=Ind_{P_{k,n-k}}^{GL_{n}}(\rho _{E}\otimes \mathbf{1}_{n-k})$ where $%
\rho _{E}$ is a SPS representation of the $GL_{k}$-block of the parabolic $%
P_{k,n-k}.$ In fact, $E\in \mathcal{Y}_{k}$ is the Young diagram obtained
from $D$ by deleting its first row. So the strict co-tensor rank of $\rho
_{D}$ is $d_{1}$. This proves Part (\ref{C-TR-PofCF-SPS}).

Part (\ref{C-TR-PofCF-general}), i.e., the case of general irrep, is proved
in a similar manner. The Pieri rule implies Formula (\ref{I_tau_US}), and
using it and all of its twists by characters, we deduce the statement
applying the same argument as in the SPS case.

This completes the proof of Corollary \ref{C-TR-PofCF}.
\end{proof}

\subsection{\textbf{Proofs for Section \protect\ref{S-Der-AI-GLn}}}

\subsubsection{\textbf{Proof of Proposition }\protect\ref{P-res-n-U_1,n-1} 
\label{P-P-res-n-U_1,n-1}}

\begin{proof}
Recall that the unipotent radical $U_{n-1,1}$ of the parabolic $P_{n-1,1}$
is isomorphic to $\mathbb{F}_{q}^{n-1},$ and that any non-identity element
in that group is a transvection. So all non-identity elements are conjugate,
and dually, all non-identity characters are conjugate. So the restriction of
any representation of $GL_{n}$ to $U_{n-1,1}$ is a sum of some copies of the
trivial representation $\mathbf{1}_{U_{n-1,1}}$, and some copies of $%
reg_{_{U_{n-1,1}}}^{\circ }=reg_{_{U_{n-1,1}}}-$ $\mathbf{1}_{U_{n-1,1}}$
the regular representation minus the trivial representation.

Let $\rho $ be an irreducible representation of $GL_{n}$. If the restriction
of $\rho $ to $U_{n-1,1}$ contains $\mathbf{1}_{U_{n-1,1}}$, then by
Frobenius reciprocity, $\rho $ must be contained in a representation induces
from the parabolic $P_{n-1,1}$. This means that its realization in terms of
the philosophy of cusp forms must be induction from a parabolic whose Levi
component contains some $GL_{1}$ factors, which in turn means, by Corollary %
\ref{C-TR-PofCF}, that $\rho $ has tensor rank at most $n-1$.

We conclude that tensor rank $n$ irreps restricted to $U_{n-1,1}$ contain
only multiple of $reg_{_{U_{n-1,1}}}^{\circ }$, and the character ratio of
such a representation on the transvection will be the character ratio, on
that element, of $reg_{_{U_{n-1,1}}}^{\circ }$, which is $\frac{-1}{q^{n-1}-1%
}$. This completes the proof of Proposition \ref{P-res-n-U_1,n-1}.
\end{proof}

\subsubsection{\textbf{Proof of Proposition \protect\ref{P-CR-I_N}\label%
{P-P-CR-I_N}}}

\begin{proof}
The proof is by a direct computation of the ratio between the cardinalities
of an appropriate set of flags of vector spaces and its subset of \textit{%
transvection invariant flags}.

The representation $I_{D}$, where $D=\{d_{1}\geq ...\geq d_{r}\}$ is a
partition of $n$, can be realized on the space of functions on the set $%
X_{F} $ of flags of vector spaces in $V=\mathbb{F}_{q}^{n}$ of the form 
\begin{equation}
F:\text{ \ }0=V_{a_{0}}\subset V_{a_{1}}\subset \ldots \subset
V_{a_{r-1}}\subset V_{a_{r}}=V\text{,}  \label{F}
\end{equation}%
where $\dim (V_{a_{j}})=a_{j},$ and $\dim (V_{a_{j}}/V_{a_{j-1}})=d_{j}$,
for $j=1,...,r$.

Let $T$ be a transvection on $V$, i.e., $T-I$ (here, $I$ stands for the
identity operator) has rank one, and $(T-I)^{2}=0$.

We are interested in knowing what restrictions the flag $F$ (\ref{F}) must
satisfy in order to be invariant under $T$. We treat two extreme cases, and
then the general case.\medskip

\textbf{Case 1. }\textit{Suppose the line }$L=\func{Im}(T-I)$\textit{\ is
not contained in\ }$V_{a_{r-1}}$\textit{. }

In this case for $F$ to be invariant under $T$, all the $V_{a_{j}}$, and in
particular, $V_{a_{r-1}}$, must be contained in $\ker (T-I)$, which is a
hyperplane - a subspace of $V$ of co-dimension $1$. In other words, $F$ is
actually a flag in $(n-1)$-space rather than $n$-space.

The collection of flags with the given subspace dimensions in $(n-1)$-space
rather than $n$-space has relative cardinality $\frac{1}{q^{\dim
(V_{a_{r-1}})}}+o(...)$ with respect to the collection of all such flags in $%
n$-space, as one sees by comparing opposite unipotent radicals in the two
situations.\medskip

\textbf{Case 2. }\textit{Suppose, on another hand, that the line }$L$\textit{%
\ is contained already in }$V_{a_{1}}$\textit{. }

In this case $F$ is guaranteed to be invariant under $T$. How many flags can
satisfy this condition? If $L$ is contained in $V_{a_{1}}$, then the flag $F$
will push down to define a flag in the $(n-1)$-dimensional space $V/L$.
Again comparing opposite unipotent radicals, we see that the relative
cardinality of this collection of $T$-invariant $F$ with respect to the
collection of all possible $F$ is $\frac{1}{q^{n-\dim (V_{a_{1}})}}+o(...)$%
.\medskip

\textbf{Case 3.}\textit{\ Now consider the general situation: suppose }$%
V_{a_{j}}\nsupseteq L\subset V_{a_{j+1}}$\textit{.}

Here, by looking only at the sub-flag%
\begin{equation*}
F_{a_{j}}:\text{ \ }V_{a_{1}}\subset \ldots \subset V_{a_{j}}\subset V,
\end{equation*}%
we conclude that $F_{a_{j}}$ is part of a collection of flags of relative
cardinality $\frac{1}{q^{\dim (V_{a_{j}})}}$ with respect to the set of all
flags with the same dimension set as $F_{a_{j}}$.

On the other hand, consider the flag 
\begin{equation*}
_{a_{j}}F:\text{ \ }V_{a_{j+1}}/V_{a_{j}}\subset ...\subset
V_{a_{r-1}}/V_{a_{j}}\subset V/V_{a_{j}}.
\end{equation*}%
It satisfies the second simplified condition of Case 2. This implies that $%
_{a_{j}}F$ is part of in a collection of flags of relative cardinality $%
\frac{1}{q^{n-\dim (V_{a_{j+1}})}}+o(...)$ with respect to the collection of
all flags with the dimension set of $_{a_{j}}F$.

The mapping $F$ $\mapsto $ $F_{a_{j}}$ defines a surjective map from the set 
$X_{F}$ of flags with dimension set the same as $F$, to the collection of
flags with the dimension set of $F_{a_{j}}$. This map is a fibration, with
fiber equal to the collection of flags with dimension set equal to that of $%
_{a_{j}}F$. Looking at the inverse image of the $T$-invariant set in the
fiber over each point of the $T$-invariant set in $F_{a_{j}}$ satisfying the
condition $\func{Im}(T-I)\subset V_{a_{j+1}}$, we conclude that the set of
of flags with the dimension set of $F$ and such that $\func{Im}(T-I)$ is
contained in $V_{a_{j+1}}$ but not in $V_{a_{j}}$ has relative cardinality%
\begin{equation*}
\frac{1}{q^{n-\dim (V_{a_{j+1}})+\dim (V_{a_{j}})}}+o(...)=\frac{1}{%
q^{n-\dim (V_{a_{j+1}}/V_{a_{j}})}}+o(...),
\end{equation*}%
with respect to the set $X_{F}$.

Taking the minimum of the numbers 
\begin{equation*}
n-\underset{j}{\max }\dim (V_{a_{j+1}}/V_{a_{j}}),
\end{equation*}%
which in our case is $n-d_{1},$ and denote by $m_{d_{1}}$ the number of
times the quantity $d_{1}$ appears in the partition $D$, we conclude that
the relative cardinality of the set $X_{F}^{T},$ of $T$-invariant flags with
the dimension set of $F,$ with respect to the set $X_{F}$ is 
\begin{equation}
\frac{\#X_{F}^{T}}{\#X_{F}}=\frac{m_{d_{1}}}{q^{n-d_{1}}}+o(...).
\label{CR-SPS-I_N}
\end{equation}%
Of course (\ref{CR-SPS-I_N}) is equal to $\frac{\chi _{I_{D}}(T)}{\dim
(I_{D})}$. This completes the proof of Proposition \ref{P-CR-I_N}.
\end{proof}

\subsubsection{\textbf{Proof of Proposition \protect\ref{P-Rel-CR-I_D}\label%
{P-P-Rel-CR-I_D}}}

\begin{proof}
Fix an algebraic closure $\mathbf{k}$ of $\mathbb{F}_{q}.$ Consider the flag
variety $\mathbf{X}_{D}$ of flags in $\mathbf{V}=\mathbf{k}^{n}$ defined by
a partition $D=\{d_{1}\geq ...\geq d_{r}\}$ of $n$ \cite{Fulton97}. It is
irreducible of dimension 
\begin{equation*}
\Delta =\dim (\mathbf{X}_{D})=\sum_{1\leq i<j\leq r}d_{i}d_{j}\text{.}
\end{equation*}

Consider a transvection $T$ acting on $\mathbf{X}_{D}$. It has fixed points
(this was explained in the proof just above) that form a Zariski open subset
of a union of flag varieties These flag varieties have dimensions%
\begin{equation*}
\Delta _{i}=\Delta -n+d_{i}.
\end{equation*}%
Consider a partition $D^{\prime }$ that dominates $D$. Then $D^{\prime }$
can be reached from $D$ by a sequence of transformations of the type $%
d_{a}\mapsto d_{a}+1$; and $d_{b}\mapsto d_{b}-1$, for $a<b,$ and leaving
the other $d_{i}$'s unchanged.

For this $D^{\prime }$, we have%
\begin{eqnarray*}
\Delta ^{\prime } &=&\dim (\mathbf{X}_{D^{\prime }})=\Delta +\sum_{i\neq
a}d_{i}-\sum_{i\neq b}d_{i}-1 \\
&=&\Delta -d_{a}+d_{b}-1.
\end{eqnarray*}%
Then the subvarieties of fixed points for the transvection on this $%
X_{D^{\prime }}$ have dimensions%
\begin{equation*}
\Delta _{i}^{\prime }=\Delta ^{\prime }-n+d_{i}^{\prime }.
\end{equation*}

Since we have $d_{i}^{\prime }=d_{i}$ except for $i=a,b,$ and since $%
d_{b}^{\prime }=d_{b}-1$, and since $\Delta ^{\prime }<\Delta $ , these
dimensions are all less than for the corresponding subvarieties for $D$,
except possibly for $\Delta _{a}^{\prime }=\Delta ^{\prime }-n+d_{a}+1$. For
this to be equal to the largest dimension for $D$, we would need that $%
\Delta ^{\prime }=\Delta -1,$ and $d_{a}=d_{1}$.

The condition $\Delta ^{\prime }=\Delta -1$ in turn implies that $%
d_{b}=d_{a} $. Thus, the trace of the transvection on $I_{D^{\prime
}}=L^{2}(X_{D^{\prime }})$ (where $X_{D^{\prime }}$ denotes the set of $%
\mathbb{F}_{q}$-rational points $X_{D^{\prime }}=\mathbf{X}_{D^{\prime }}(%
\mathbb{F}_{q})$) will be of smaller order of magnitude than the trace on $%
I_{D}=L^{2}(X_{D})$, except when $d_{a}=d_{b}$, in which case, it will be of
the same order of magnitude.

This completes the proof of the Proposition.
\end{proof}

\subsection{\textbf{Proofs for Section \protect\ref{S-Der-AI-SLn}}}

\subsubsection{\textbf{Proof of Proposition \protect\ref{P-Irr-Res-Gen}\label%
{P-P-Irr-Res-Gen}}}

\begin{proof}
According to Corollary (\ref{C-Irr-Res}) the restriction to $SL_{n}$ of an
irrep $\rho $ of $GL_{n}$ is irreducible iff $\rho $ is not fixed by a twist
of any non-trivial character of $GL_{n}.$ Moreover, the eta correspondence
(see Theorem \ref{EF-TR-k}), and the \ description (see Formula (\ref%
{Form-eta-tauUSD})) of its image, tells us that, a twist by a character of a
tensor rank $k$ representation will produce isomorphic one only if the
corresponding representation of $GL_{k}$ has tensor rank 
\begin{equation}
r=2k-n.  \label{Fixed-Crit}
\end{equation}

Now, let us go over various cases:\smallskip

\textbf{The tensor rank }$k<\frac{n}{2}$\textbf{\ irreps of }$GL_{n}$\textbf{%
: }In this domain, no irreps of $GL_{k}$ has tensor rank (\ref{Fixed-Crit}).
So every irrep of $GL_{n},$ in this range, stays irreducible after
restriction to $SL_{n}$. This completes the justification of Part 1 of
Proposition \ref{P-Irr-Res-Gen}.\smallskip

\textbf{The tensor rank }$\frac{n}{2}\leq k\leq n-2$\textbf{\ irreps of }$%
GL_{n}$\textbf{: }In this interval $r=2k-n\leq k-2$, so using the counting
coming from the eta correspondence our knowledge on the cardinality of the
tensor rank $r$ irreps of $GL_{k}$ (see Theorem \ref{T-Card-trank-k}) we get
that at most $q^{k-1}+o(...)$ irreps of tensor rank $k$ of $GL_{n}$ might be
reducible after restriction to $SL_{n}$. This justify Part 2 of Proposition %
\ref{P-Irr-Res-Gen}, for the range under discussion. \smallskip\ 

\textbf{The tensor ranks }$n-1$\textbf{\ and }$n$\textbf{\ irreps of }$%
GL_{n} $\textbf{: }For the irreps of tensor rank $n-1$ of $GL_{n}$, the
"generic" split principal series representation, induced from characters of
the standard Borel subgroup which are not fixed by twist of any character of 
$GL_{n},$ stays irreducible after restriction to $SL_{n}$ (in fact these
contribute the $c_{n-1}q^{n}+o(...)$ to the cardinality of irreps of tensor
rank $k=n-1$ given in Formula (\ref{Card-k-GLn})).

A similar argument applies for the irreps of tensor rank $n$ of $GL_{n}$.
Here we use the cuspidal irreps of $GL_{n}$ parametrized by generic
characters of the torus of $GL_{n}$ defined by the multiplicative group of
the field extension $\mathbb{F}_{q^{n}}$ of degree $n$ of $\mathbb{F}_{q}$.
(See Section 7\ref{S-UB-Dim-trank-n} for more information on this
parametrization). This provides $\frac{1}{n}q^{n}+o(...)$ cuspidal irreps,
which all have tensor rank $n,$ and stays irreducible after restriction to $%
SL_{n}$. Similar arguments can be made for any maximal torus in $GL_{n}$.
Generic characters of a given torus will parametrize representations of rank 
$n$ if all the irreducible factors of the torus are multiplicative groups of
proper extensions of $\mathbb{F}_{q}$. If one or more factors of the torus
is $\mathbb{F}_{q}^{\ast }$, then the corresponding representations will
have rank $n-1$.

Overall, the above counting gives the $c_{n}q^{n}+o(...)$ irreps of tensor
rank $k=n$ given in Formula (\ref{Card-k-GLn}).

This completes the justification of Part 2 of Proposition \ref{P-Irr-Res-Gen}%
.
\end{proof}

\subsection{\textbf{Proofs for Appendix \protect\ref{A-MLGM}}}

\subsubsection{\textbf{Proof of Corollary \protect\ref{C-Res-N}\label%
{P-C-Res-N}}}

We use the notations and definitions given in Section \ref{A-MLGM}.

\begin{proof}
We start with the proof of Part (\ref{Spec}).

Of course two irreps of $G$ that differ by a twist by a character of $C$
have the same $N$-spectrum.

For the other direction. Suppose $\pi $ is an irrep of $N.$ Since $C$ is
cyclic, all possible characters of the stabilizer $C_{\pi }$ of $\pi $ arise
by restriction from some character of $C$. If we take one irreducible $\rho $
of $G$ containing $\pi $, it will be induced from some representation $%
\widetilde{\pi }$ of $C_{\pi }\ltimes N$, as described in Section \ref%
{A-MLGM} (see Part (\ref{Class}) of Theorem \ref{T-MLGM}). If we twist $\rho 
$ with a character of $C$ it will be induced from the representation $%
\widetilde{\pi }$, twisted with this character restricted to $C_{\pi }$. But
we know (see Part (\ref{Ext}) of Theorem \ref{T-MLGM}) this will give all
possible extensions of $\pi $ to $C_{\pi }\ltimes N$, and so will give all
possible representations of $\widehat{G}$ containing $\pi $ when restricted
to $N$. This completes the proof of Part (\ref{Spec}).

Part (\ref{Mult}), i.e., the multiplicity-freeness also follows from the
description of the irreducibles containing $\pi $. They are all induced from
the extension of $\pi $ (which is exactly one copy of $\pi $ on $N$) to the
stabilizer of $\pi $ under action of $C$ by conjugation. So the induced
representation restricted to $N$ consists of one copy of each $Ad^{\ast }C$
transform of $\pi $, one for each coset in $G/(C_{\pi }\ltimes N)$.
\end{proof}

\subsection{\textbf{Proofs for Appendix \protect\ref{A-IrrS_l-SPS}}}

\subsubsection{\textbf{Proof of Corollary \protect\ref{C-dim}\label{P-C-dim} 
}}

\begin{proof}
Suppose $L=\{l_{1}\geq l_{2}\geq ...\geq l_{s}\}$ is a partition of $l$. We
have $\dim (I_{L})=\#(GL_{l}/P_{L})=\#(U_{L})$, where $U_{L}$ is the
unipotent radical of $P_{L}$. An easy direct \ computation gives $%
\#U_{L}=q^{d_{L}}+o(...)$, where $d_{L}=\tsum\limits_{1\leq i<j\leq
s}l_{i}l_{j}$. It follows that 
\begin{equation}
\dim (I_{L})=q^{d_{L}}+o(...).  \label{dim-IL}
\end{equation}%
\ 

We want to show that%
\begin{equation}
d_{L^{\prime }}<d_{L},\text{ if }L^{\prime }\succ L.  \label{dL-dL'}
\end{equation}%
Indeed, suppose that for some $j_{0}<j_{1}$, we have a partition of $l$
given by $L^{\prime }=\{l_{1}\geq $ $...$ $\geq l_{j_{0}-1}\geq
l_{j_{0}}+1\geq l_{j_{0}+1}\geq $ $...$ $\geq l_{j_{1}-1}\geq
l_{j_{1}}-1\geq l_{j_{1}+1}\geq $ $...$ $\geq l_{r}\}.$ Then $L^{\prime
}\succ L$, and, in fact, the dominance order on partitions is generated by
such inequalities \cite{Ceccherini-Silberstein-Scarabotti-Tolli10}. In
particular, it is enough to show that $d_{L^{\prime }}<d_{L}$ for such $%
L,L^{\prime }$. A direct computation implies that $d_{L}>d_{L^{\prime }}$,
iff \ $l_{j_{0}}\cdot l_{j_{1}}>$ $(l_{j_{0}}+1)\cdot (l_{j_{1}}-1),$ and
the latter inequality holds true since $l_{j_{0}}\geq l_{j_{1}}$.

Now, combining Formulas (\ref{dim-IL}) and (\ref{dL-dL'}), with Proposition %
\ref{P-Mon-Max-IL}, we get that $\dim (\rho _{L})=\dim (I_{L})+o(...)=$ $%
q^{d_{L}}+o(...)$.

This completes the proof of both parts of Corollary \ref{C-dim}.
\end{proof}

\end{document}